\documentclass[11pt]{amsart}

\usepackage{amsmath,amssymb,amsthm}
\usepackage[a4paper,margin=1.1in]{geometry}
\usepackage[dvipsnames]{xcolor}
\usepackage{enumitem}
\usepackage[T1]{fontenc}
\usepackage{lmodern}
\usepackage{microtype}
\usepackage[colorlinks=true,linkcolor=blue!50!black,citecolor=green!35!black,urlcolor=blue!50!black]{hyperref}

\numberwithin{equation}{section}

\newtheorem{theorem}{Theorem}[section]
\newtheorem{proposition}[theorem]{Proposition}
\newtheorem{lemma}[theorem]{Lemma}
\newtheorem{corollary}[theorem]{Corollary}
\theoremstyle{definition}
\newtheorem{definition}[theorem]{Definition}
\theoremstyle{remark}
\newtheorem{remark}[theorem]{Remark}

\newcommand{\cO}{\mathcal{O}}
\newcommand{\fo}{\mathcal{O}_K}
\newcommand{\fa}{\mathfrak{a}}
\newcommand{\fb}{\mathfrak{b}}
\newcommand{\fc}{\mathfrak{c}}
\newcommand{\fd}{\mathfrak{d}}
\newcommand{\fe}{\mathfrak{e}}
\newcommand{\ff}{\mathfrak{f}}
\newcommand{\fg}{\mathfrak{g}}
\newcommand{\fp}{\mathfrak{p}}
\newcommand{\fq}{\mathfrak{q}}

\newcommand{\fz}{\mathfrak{z}}
\newcommand{\fs}{\mathfrak{s}}

\newcommand{\fS}{\mathfrak{S}}
\newcommand{\N}{\mathrm{N}}
\newcommand{\Nm}{\mathrm{Nm}}
\newcommand{\Tr}{\mathrm{Tr}}
\newcommand{\rad}{\mathrm{rad}}

\newcommand{\Disc}{\operatorname{Disc}}
\newcommand{\ord}{\operatorname{ord}}
\newcommand{\ee}{\mathbf{e}}
\newcommand{\bx}{\mathbf{x}}
\newcommand{\by}{\mathbf{y}}
\newcommand{\bz}{\mathbf{z}}
\newcommand{\bk}{\mathbf{k}}
\newcommand{\bb}{\mathbf{b}}
\newcommand{\bm}{\mathbf{m}}
\newcommand{\bv}{\mathbf{v}}
\newcommand{\bu}{\mathbf{u}}
\newcommand{\bw}{\mathbf{w}}
\newcommand{\ba}{\mathbf{a}}
\newcommand{\bt}{\mathbf{t}}
\newcommand{\bs}{\mathbf{s}}
\newcommand{\bxi}{\boldsymbol{\xi}}
\newcommand{\bl}{\boldsymbol{\lambda}}
\newcommand{\tv}{\tilde{\mathbf{v}}}
\newcommand{\tm}{\tilde{\mathbf{m}}}
\newcommand{\cH}{\mathcal{H}}
\newcommand{\cW}{\mathcal{W}}
\newcommand{\cZ}{\mathcal{Z}}
\newcommand{\cS}{\mathcal{S}}
\newcommand{\cK}{\mathcal{K}}
\newcommand{\cN}{\mathcal{N}}
\newcommand{\bW}{\mathbf{W}}
\newcommand{\ip}[1]{\langle #1\rangle}
\newcommand{\RR}{\mathbb{R}}
\newcommand{\CC}{\mathbb{C}}
\newcommand{\QQ}{\mathbb{Q}}
\newcommand{\ZZ}{\mathbb{Z}}
\newcommand{\FF}{\mathbb{F}}
\newcommand{\PP}{\mathbb{P}}
\newcommand{\sumstar}{\sideset{}{^*}\sum}

\newcommand{\ep}[1]{e_{\fp^{#1}}}

\begin{document}

\title[Pairs of quadrics over number fields]{Pairs of quadrics in eleven variables\\ over number fields}

\author{Esrafil Ali Molla}
\address{Stat-Math Unit, Indian Statistical Institute, Kolkata, 203 B.T. Road, Baranagar, West Bengal 700108, India}
\email{esrafil.math@gmail.com}


\begin{abstract}
A nested version of the $\delta$-method is developed for number fields $K/\QQ$ and is used to establish an asymptotic formula for the number of $K$-rational points of bounded height on non-singular complete intersections of two quadrics in $\PP^{n-1}_K$, when $n\geqslant11$.
\end{abstract}

\maketitle

\setcounter{tocdepth}{1}
\tableofcontents

\section{Introduction}\label{sec:intro}

\subsection{Background}
Let $K$ be a number field of degree $d$ over $\QQ$, with ring of integers $\fo$. Let $Q_1,Q_2\in\fo[x_1,\dots,x_n]$ be quadratic forms such that
\[
X=\{Q_1=Q_2=0\}\subset\PP^{n-1}_K
\]
is a non-singular complete intersection of codimension two. The arithmetic of such varieties has a long history. Colliot-Th\'el\`ene, Sansuc and Swinnerton-Dyer~\cite{CSS87a,CSS87b} proved that the Hasse principle and weak approximation hold for $X$ over any number field as soon as $n\geqslant 9$, and Heath-Brown~\cite{HB18} has extended the Hasse principle to $n\geqslant8$. For $n\geqslant9$ the variety has points over every non-archimedean completion of $K$, by work of Demyanov~\cite{Dem56}, who treats pairs of quadratic forms over complete discretely valued fields with finite residue field, and of Birch, Lewis and Murphy~\cite{BLM62}; compare the remark in~\cite[\S1]{HB18}.

The circle method gives more: an asymptotic formula for the number of points of bounded height. Birch's theorem~\cite{Bir62}, extended to number fields by Skinner~\cite{Ski97}, provides such a formula for non-singular pairs of quadrics once $n\geqslant14$. Over $\QQ$, Munshi~\cite{Mun15} lowered this to $n\geqslant11$. His method, which he called a \emph{nested $\delta$-method}, was the first application of the $\delta$-method of Duke, Friedlander and Iwaniec~\cite{DFI93} and Heath-Brown~\cite{HB96} to a system of more than one non-linear equation. It builds on earlier work of Browning and Munshi~\cite{BM13} on singular intersections of quadrics. More recently, Li, Rydin Myerson and Vishe~\cite{LRV24} have introduced a two-dimensional $\delta$-symbol over $\QQ$ which handles $n\geqslant10$, and Vishe~\cite{Vis23} has treated $n\geqslant9$ over $\FF_q(t)$; see also Heath-Brown and Pierce~\cite{HBP17}.

Over a number field the main difficulty in setting up the $\delta$-method is that there is no Farey dissection. Browning and Vishe~\cite{BV14} overcame this with a smooth $\delta$-symbol that detects the zero \emph{ideal}, and used it to count points on cubic hypersurfaces over $K$. The purpose of this paper is to carry Munshi's argument over to an arbitrary number field.

\subsection{The result}
Let $V=K\otimes_\QQ\RR$, and let $\Nm$ and $\ee(\cdot)$ be the norm and additive character on $V$ described in \S\ref{sec:notation}. Write $X^\ast(V)$ for the set of $\by\in V^n$ such that, at some archimedean place, the gradients $\nabla Q_1$ and $\nabla Q_2$ are linearly dependent at $\by$. Fix $\bxi\in V^n\setminus X^\ast(V)$, and let $W$ be the weight~\eqref{eq:weight}, a smooth bump around $\bxi$ multiplied by a Gaussian of width $(\log P)^{-2}$. Put
\[
N(P)=\sum_{\substack{\bx\in\fo^n\\ Q_1(\bx)=Q_2(\bx)=0}}W(\bx/P).
\]
Let $\sigma_\infty$ be the singular integral attached to $W$ (see~\eqref{eq:sigmainf}), and let $\fS_K$ be the singular series~\eqref{eq:singularseries}. Our main result is the following.

\begin{theorem}\label{thm:main}
Let $K$ be a number field of degree $d$, with discriminant $D_K$ and $r_2$ complex places. Let $n\geqslant 11$, and let $Q_1,Q_2\in\fo[x_1,\dots,x_n]$ define a non-singular complete intersection $X\subset\PP^{n-1}_K$. There is $\eta=\eta(n,d)>0$ such that
\[
N(P)=\Bigl(\frac{2^{r_2}}{\sqrt{D_K}}\Bigr)^{n-2}\fS_K\,\sigma_\infty\,P^{d(n-4)}+O\bigl(P^{d(n-4)-\eta}\bigr),
\]
uniformly for $\bxi$ in any fixed compact subset of $V^n\setminus X^\ast(V)$. The series $\fS_K$ converges absolutely, it is a product of local densities, and it is positive. If $\bxi\in X(V)$, then $\sigma_\infty\asymp(\log P)^{-2d(n-2)}$.
\end{theorem}

The implied constants depend on $K$, $Q_1$, $Q_2$, $n$ and the compact set. One may take $\eta=\min\{\frac d4(n-10-2\theta_K),\frac d4(n-8)\}-\varepsilon$ for any $\varepsilon>0$, where $\theta_K\in[0,\frac12)$ is defined in~\eqref{eq:theta}. Positivity of $\fS_K$ uses the existence of non-singular points of $X$ over every non-archimedean completion of $K$, which holds for $n\geqslant9$ by the results quoted above.

A standard smoothing argument turns Theorem~\ref{thm:main} into an asymptotic formula for sharp boxes.

\begin{corollary}\label{cor:boxes}
Let $K$, $X$ and $n$ be as in Theorem~\ref{thm:main}. Let $\mathcal B\subset V^n$ be a product, over the archimedean places, of boxes whose closure does not meet $X^\ast(V)$. Then
\[
\#\{\bx\in\fo^n\cap P\mathcal B:Q_1(\bx)=Q_2(\bx)=0\}=\Bigl(\frac{2^{r_2}}{\sqrt{D_K}}\Bigr)^{n-2}\fS_K\,\sigma_\infty(\mathcal B)\,P^{d(n-4)}+O\Bigl(\frac{P^{d(n-4)}}{\log P}\Bigr),
\]
where $\sigma_\infty(\mathcal B)$ is the singular integral of $\mathcal B$.
\end{corollary}

For $K=\QQ$, Theorem~\ref{thm:main} recovers Munshi's range $n\geqslant11$. The arithmetic input that brings the range down to $n=11$ is a subconvexity bound for Hecke $L$-functions of $K$. We use it in the following form, which is the case of $\mathrm{GL}_1$ in~\cite[Theorem~1.1]{MV10}.

\begin{theorem}\label{thm:MV}
There is an absolute constant $\delta_{\mathrm{MV}}>0$ such that $L(\frac12,\chi)\ll_K\mathfrak C(\chi)^{1/4-\delta_{\mathrm{MV}}}$ for every unitary Hecke character $\chi$ of $K$, where $\mathfrak C(\chi)$ is the analytic conductor of $\chi$.
\end{theorem}

The convexity bound is the case $\delta_{\mathrm{MV}}=0$. Michel and Venkatesh prove Theorem~\ref{thm:MV} for automorphic representations of $\mathrm{GL}_1$ and $\mathrm{GL}_2$ over $K$ with unitary central character, uniformly in all aspects; the constant $\delta_{\mathrm{MV}}$ does not depend on $K$. In our application the infinity types of the characters are bounded by $P^{\varepsilon}$, so that $\mathfrak C(\chi)$ is at most $P^{O(\varepsilon)}$ times the norm of the finite conductor. Put
\begin{equation}\label{eq:theta}
\theta_K=\begin{cases}\frac12-2\delta_{\mathrm{MV}}&\text{if $n$ is odd},\\ 0&\text{if $n$ is even},\end{cases}\qquad\Phi(n)=\frac{3n}4-\frac32+\frac{\theta_K}2.
\end{equation}
We shall show that the error term in Theorem~\ref{thm:main} is $O(P^{d\Phi(n)+\varepsilon}+P^{d(n-4)-d(n-8)/4+\varepsilon})$. The main term dominates when $n-4>\Phi(n)$, that is, when $n>10+2\theta_K$. For odd $n$ this needs $\theta_K<\frac12$, which is exactly the statement $\delta_{\mathrm{MV}}>0$. With only the convexity bound one would need $n\geqslant12$. For even $n$ no subconvexity is needed. For $K=\QQ$ and $\theta=\frac7{16}$, the value implicit in~\cite[Lemma~5]{Mun15}, our error term is $P^{3n/4-41/32+\varepsilon}$, which is Munshi's.

\subsection{The nested $\delta$-method over $K$}
We now explain the strategy. For $\bx\in\fo^n$ in the support of $W(\bx/P)$, both $Q_1(\bx)$ and $Q_2(\bx)$ have size $P^2$ at every archimedean place. A $\delta$-symbol detects the vanishing of an element of size $Q^2$ (at every place) with moduli of norm up to $Q^d$. Detecting the two equations independently would therefore require moduli $\fb_1,\fb_2$ with $\N\fb_1,\N\fb_2\ll P^d$, so that the total conductor $\N(\fb_1\fb_2)$ would reach $P^{2d}$. After Poisson summation over $\bx$, which has $n$ coordinates of size $P$, the dual sum then has as many terms as the original sum. Nothing is gained.

Munshi's idea is to use the first modulus to reduce the size of the second equation. Over $K$ this takes the following form. Suppose that $Q_1(\bx)=0$ has been detected by a $\delta$-symbol with modulus $\fb_1$. Then
\[
Q_2(\bx)=0\iff Q_2(\bx)\in\fb_1\quad\text{and}\quad(Q_2(\bx))\fb_1^{-1}=(0).
\]
The congruence $Q_2(\bx)\in\fb_1$ costs nothing: it is periodic modulo $\fb_1$, and so it can be absorbed into the exponential sum for the first equation. The ideal $(Q_2(\bx))\fb_1^{-1}$ is integral and has norm about $P^{2d}/\N\fb_1$, so it can be detected by the ideal $\delta$-symbol of Browning and Vishe with a modulus $\fb_2$ of norm about $P^{d/2}$, once a smooth cut-off has made its size uniform. The total conductor is then at most about $P^{3d/2}$. This is where the saving comes from. It is notable that the ideal formulation of~\cite{BV14}, which was introduced to compensate for the lack of a Farey dissection, is exactly what makes this nesting possible over $K$: the element $Q_2(\bx)/\beta$ for a generator $\beta$ of $\fb_1$ need not exist, but the ideal $(Q_2(\bx))\fb_1^{-1}$ always does.

After Poisson summation one obtains complete exponential sums $C_{\fb_1,\fb_2}(\bm)$ modulo $\fb_1\fb_2$, which are the analogues over $K$ of Munshi's sums, and oscillatory integrals $I_{\fb_1,\fb_2}(\bm)$. The sums admit a Kloosterman refinement in $\fb_1$ (square-root cancellation in the sum over the first modulus for generic $\bm$) and a ``double'' Kloosterman refinement in $\fb_2$: their dependence on $\fb_2$ is through a quadratic character, and averaging it over $\fb_2$ leads to Hecke $L$-functions, to which Theorem~\ref{thm:MV} applies. Section~\ref{sec:nested} gives the details.

\subsection{Archimedean sizes}
Over $\QQ$ the size of a number is a single quantity. Over $K$ it is a vector, one entry for each archimedean place, and norms control only the product of the entries. This creates a difficulty that has no counterpart over $\QQ$ or over imaginary quadratic fields.

Suppose that both equations were detected by the ideal $\delta$-symbol of Browning and Vishe. Its archimedean weight depends on $Q_i(\bx)$ only through the norm $|\Nm Q_i(\bx)|$. The dual variables $t_1,t_2\in V$ arising from the two weights are then controlled only through their heights $\cH(t_1)$ and $\cH(t_2)$, and in particular they may be concentrated at \emph{different} archimedean places; in other words, they may be highly skewed rather than balanced. The stationary phase analysis then produces dual frequencies $\bm$ far beyond the range predicted by the size of the moduli; for a real quadratic field one loses a power of $P$ which exceeds the main term. We explain this in detail in Section~\ref{sec:higher}.

We avoid the difficulty with two devices.
\begin{enumerate}[label=(\roman*)]
\item The first equation is detected with a $\delta$-symbol whose divisors are \emph{elements} $c\in\fo$, and whose archimedean weight $\bW(x)=\prod_lw_1(|x^{(l)}|)$ is a product over the places (Theorem~\ref{thm:elementdelta}). This is Heath-Brown's construction~\cite{HB96} applied to elements rather than ideals; the infinite unit group causes no trouble because $\bW$ is compactly supported in $V^\times$. The use of archimedean weights that are products over the places is in the spirit of Blomer and Brumley~\cite{BB11}. The dual variable of the first equation is then bounded at each place $l$ by $P/|c_1^{(l)}|$.
\item The cut-off for the second equation is $U_V(\bl Q_2(\bx)/P^2)$, where the vector $\bl=\bl(\fb_1,c_1)\in V$ has, up to a scalar, the same archimedean sizes as $c_1^{-1}$. This is reminiscent of the generators of principal ideals with prescribed archimedean sizes that Qi uses in~\cite{Qi24}. The norm of $\bl$ is still $P^d/\N\fb_1$, so the conductor is still at most about $P^{3d/2}$.
\end{enumerate}
With these choices the dual variable of the second equation dominates that of the first at every place separately, and the archimedean truncation of the dual sum follows place by place (Theorem~\ref{thm:truncation}). The second equation still uses the ideal $\delta$-symbol of Browning and Vishe. In particular, the exponential sums are exactly the analogues of Munshi's sums, and the arithmetic part of the argument is unaffected by the change.

\subsection{Organisation and conventions}
Section~\ref{sec:prelim} collects algebraic and geometric preliminaries. Section~\ref{sec:delta} treats the two $\delta$-symbols. Section~\ref{sec:nested} carries out the nesting and the Poisson summation. Section~\ref{sec:sums} studies the exponential sums and Section~\ref{sec:integrals} the exponential integrals. The error term is estimated in Section~\ref{sec:error} and the main term in Section~\ref{sec:main}, and Section~\ref{sec:proofs} completes the proofs. Section~\ref{sec:higher} explains why weights that depend only on norms fail when $K$ has more than one archimedean place.

To the best of our knowledge, this is the first time that the nested $\delta$-method has been carried out over a number field other than $\QQ$. For this reason we have tried to make the paper self-contained and to give complete proofs. Each result that we use from~\cite{BV14}, \cite{BM13} or~\cite{Mun15} is stated as a separate lemma and proved, with the changes needed over $K$. In several places we give a different proof from the original. The main-term property of the Browning--Vishe weight (Lemma~\ref{lem:h4}) is proved using local zeta integrals. The Fourier bounds for that weight (Lemma~\ref{lem:p2}) are proved by decomposing it into pieces supported on hyperbolic shells. The bounds for the exponential sums at the finitely many bad primes, including those above $2$, follow from a single uniform estimate (Lemma~\ref{lem:universal}). The stationary phase analysis is carried out on average over the pencil (Proposition~\ref{prop:statphase}), which avoids any assumption on the conditioning of the matrices $M(t_1,t_2)$. The new ingredients are the element $\delta$-symbol (Theorem~\ref{thm:elementdelta}), the place-wise nesting, adapted to the archimedean sizes of the divisor (Section~\ref{sec:nested}) and the resulting archimedean truncation (Theorem~\ref{thm:truncation}).

The following external results are used without proof: Dirichlet's unit theorem and basic facts about Dedekind zeta functions, Hecke characters and Hecke $L$-functions, including the prime ideal theorem in ideal classes and the zero-free region with Siegel's theorem~\cite{IK04,Neu99}; Deligne's bound for smooth complete intersections over finite fields~\cite{Del74}; Weil's bound for character sums; class field theory for quadratic extensions; Cohen's theorem on thin sets~\cite[Chap.~13]{Ser97}; the $\fp$-adic implicit function theorem and change of variables formula~\cite{Igu00}; Takagi's factorisation of complex symmetric matrices~\cite{HJ13}; Theorem~\ref{thm:MV}; and the classical description of smooth intersections of two quadrics~\cite{Rei72}.

\subsection{Use of AI tools}
In the preparation of this paper the author used the AI assistant Claude (Anthropic) as a research and writing assistant. The research problem, the overall strategy and the scope of the paper are the author's. Working from the author's material and instructions, Claude suggested key technical ingredients, notably the $\delta$-symbol over elements and the place-wise nesting (Theorems~\ref{thm:elementdelta} and~\ref{thm:nested}), and helped to write out proofs, adapt arguments from the literature to number fields, check computations and draft the text. All statements and proofs have been checked by the author, who takes full responsibility for the content of the paper.

\section{Preliminaries}\label{sec:prelim}

\subsection{Notation}\label{sec:notation}
Throughout, $K$ is a number field of degree $d=r_1+2r_2$, where $r_1$ and $r_2$ are the numbers of real and complex places. We write $\fo$ for its ring of integers, $\fd$ for the different, $D_K$ for the absolute value of the discriminant, $h_K$ for the class number, $U_K$ for the unit group and $w_K$ for the number of roots of unity in $K$. The residue of the Dedekind zeta function at $s=1$ is denoted by $\Delta_K$.

We follow~\cite[\S2]{BV14} for the archimedean notation. Let
\[
V=K\otimes_\QQ\RR=\bigoplus_{l=1}^{r_1+r_2}K_l,
\]
where $K_l=\RR$ for $1\leqslant l\leqslant r_1$ and $K_l=\CC$ for $r_1<l\leqslant r_1+r_2$. Put $c_l=1$ in the first case and $c_l=2$ in the second. For $v=(v^{(l)})\in V$ let
\[
\Nm(v)=\prod_{l=1}^{r_1}v^{(l)}\prod_{l=r_1+1}^{r_1+r_2}|v^{(l)}|^2,\qquad\Tr(v)=\sum_lc_l\Re v^{(l)},\qquad\ee(v)=e^{2\pi i\Tr(v)}.
\]
For $\alpha\in K$ these are the field norm and trace. We put $\ip v=\max_l|v^{(l)}|$ and $\|v\|^2=\sum_lc_l|v^{(l)}|^2$. On $V$ we use the measure $dv$ that is Lebesgue measure on each real factor and the measure $dx\,dy$ on each complex factor $x+iy$. With this normalisation a fractional ideal $\fa$, viewed as a lattice in $V$, has covolume $2^{-r_2}\sqrt{D_K}\N\fa$. Its dual with respect to the pairing $(v,w)\mapsto\Tr(vw)$ is $\hat\fa=\fa^{-1}\fd^{-1}$. For $\bv=(v_1,\dots,v_m)\in V^m$ we put $\bv\cdot\bw=\sum v_iw_i\in V$, $\ip\bv=\max_i\ip{v_i}$ and $|\bv^{(l)}|=\max_i|v_i^{(l)}|$, and we define
\begin{equation}\label{eq:height}
\cH(\bv)=\prod_l\max\bigl(1,|\bv^{(l)}|\bigr)^{c_l},\qquad\cN(\bv)=\prod_l|\bv^{(l)}|^{c_l}.
\end{equation}
Clearly $\cH(\bv)\geqslant\cN(\bv)$ and $\cH(\bv)\leqslant\max(1,\ip\bv)^d$. For $m=1$, $\cH$ is the function used in~\cite[(5.5)]{BV14}. We use the same letters for the Fourier transform on $V^m$,
\[
\hat f(\bt)=\int_{V^m}f(\bv)\ee(-\bt\cdot\bv)\,d\bv,
\]
for which $f(\bv)=\int\hat f(\bt)\ee(\bt\cdot\bv)\,d\bt$ when $f$ is a Schwartz function. Products and quotients of elements of $V$ are taken componentwise. In particular, for $\mu\in V$ and $v\in V$, $\mu v$ has components $\mu^{(l)}v^{(l)}$.

For a function $f\in C_c^\infty(V^m)$ we write $\|f\|_{C^k}$ for the maximum of the absolute values of the partial derivatives of $f$ of order at most $k$ with respect to the $dm$ real coordinates. For an ideal $\fa$ we write $\N\fa$ for its absolute norm and $\omega(\fa)$ for the number of its prime divisors. We write $\delta_K(\fa)=1$ if $\fa=(0)$ and $\delta_K(\fa)=0$ otherwise, and $\delta_K(\alpha)=\delta_K((\alpha))$ for $\alpha\in\fo$. Unless stated otherwise, implied constants may depend on $K$, $Q_1$, $Q_2$, $n$, $\varepsilon$ and the fixed weights below. The letter $\varepsilon$ denotes a small positive number whose value may change from one occurrence to the next. The large parameter is $P$, and we put
\begin{equation}\label{eq:L}
L=(\log P)^4.
\end{equation}

\subsection{The geometry of the pencil}\label{sec:geometry}
Write $Q_i(\bx)=\bx^tM_i\bx$, where $M_i$ is symmetric and $2M_i$ has entries in $\fo$, and put
\[
M(\bs)=s_1M_1+s_2M_2,\qquad \mathcal{P}(\bs)=\det M(\bs),\qquad F_{\bs}=s_1Q_1+s_2Q_2.
\]
For $\bm\in\bar K^n$ define the binary form
\begin{equation}\label{eq:Phi}
\Phi_{\bm}(\bs)=-\bm^t\operatorname{adj}(M(\bs))\,\bm=\det\begin{pmatrix}M(\bs)&\bm\\ \bm^t&0\end{pmatrix},
\end{equation}
which has degree $n-1$ in $\bs$ and degree $2$ in $\bm$. The second expression follows by expanding the determinant along the last row and column.

We use the following classical criterion; see Reid~\cite[Proposition~2.1]{Rei72}; compare also~\cite[\S1]{CSS87a}.

\begin{lemma}\label{lem:smoothcrit}
Let $F$ be a field of characteristic different from $2$, let $N\geqslant3$, and let $A,B$ be symmetric $N\times N$ matrices over $F$. Then $\{\bx^tA\bx=\bx^tB\bx=0\}\subset\PP^{N-1}_F$ is a non-singular complete intersection of codimension two if and only if the binary form $\det(s_1A+s_2B)$ is not identically zero and has $N$ distinct roots in $\PP^1(\bar F)$. In this case $\operatorname{rank}(s_1A+s_2B)\geqslant N-1$ for every $\bs\neq\mathbf0$ over $\bar F$.
\end{lemma}

We record the consequences that we need.

\begin{lemma}\label{lem:geometry}
Let $X$ be non-singular. Then the following hold.
\begin{enumerate}[label=\textup{(G\arabic*)}]
\item\label{G1} The binary form $\mathcal{P}(\bs)$ has non-zero discriminant, and $\operatorname{rank}M(\bs)\geqslant n-1$ for every $\bs\neq\mathbf0$ over $\bar K$.
\item\label{G2} After replacing $(Q_1,Q_2)$ by suitable $\fo$-linear combinations, which changes neither $X$ nor the counting problem, we have $\det M_1\det M_2\neq0$.
\item\label{G3} The set $X^\ast=\{\bx:\operatorname{rank}(M_1\bx,M_2\bx)\leqslant1\}$ is the union of the $n$ lines $\ker M(\bs)$, where $\bs$ runs over the roots of $\mathcal{P}$. In particular $X^\ast$ meets the affine cone over $X$ only in $\mathbf0$.
\item\label{G4} For $\bm\neq\mathbf0$, the hyperplane section $X\cap\{\bm\cdot\bx=0\}$ is a non-singular complete intersection in $\{\bm\cdot\bx=0\}\cong\PP^{n-2}$ if and only if $\Disc\Phi_{\bm}\neq0$. Moreover there is an absolutely irreducible form $G\in\fo[\bx]$ of degree $4(n-2)$, defining the dual variety of $X$, and a constant $c_G\in\fo\setminus\{0\}$ with $c_G\Disc\Phi_{\bm}=G(\bm)$ identically.
\end{enumerate}
\end{lemma}

\begin{proof}
Part~\ref{G1} is Lemma~\ref{lem:smoothcrit} over $\bar K$. For~\ref{G2}, note that $\mathcal{P}$ has only $n$ roots in $\PP^1$. Choose $\bs,\bs'\in\fo^2$, linearly independent, that are not roots, and replace $(Q_1,Q_2)$ by $(F_{\bs},F_{\bs'})$. For~\ref{G3}, a point $\bx\neq\mathbf0$ lies in $X^\ast$ if and only if $M(\bs)\bx=\mathbf0$ for some $\bs\neq\mathbf0$, which forces $\mathcal{P}(\bs)=0$; by~\ref{G1}, each such $\ker M(\bs)$ is a line. If $\bx\neq\mathbf0$ lies on the cone over $X$ and in $X^\ast$, then $\nabla Q_1(\bx)$ and $\nabla Q_2(\bx)$ are dependent, so $X$ is singular at $[\bx]$, a contradiction.

For~\ref{G4}, choose coordinates so that $\bm\cdot\bx=x_n$. The restriction of $M(\bs)$ to $\{x_n=0\}$ is its upper left $(n-1)\times(n-1)$ block $M'(\bs)$, and by~\eqref{eq:Phi} we have $\Phi_{\bm}(\bs)=-\det M'(\bs)$ in these coordinates. A linear change of coordinates multiplies $\Phi_{\bm}$ by a non-zero constant. Lemma~\ref{lem:smoothcrit} with $N=n-1$ now gives the first assertion. Since a general hyperplane section of $X$ is non-singular, $\Disc\Phi_{\bm}$ is a non-zero form in $\bm$. Its degree is $4(n-2)$, since the discriminant of a binary form of degree $n-1$ has degree $2(n-2)$ in the coefficients. By the first assertion it vanishes exactly on the dual variety $X^\vee$ of $X$. By Ein~\cite[Proposition~3.1]{Ein86} $X^\vee$ is a hypersurface, by Aznar~\cite[Theorem~3]{Azn80} its degree is $4(n-2)$, and it is irreducible as the image of the irreducible conormal variety; see also~\cite[\S2.2]{BM13}. Let $G_0$ be an irreducible form defining $X^\vee$. Then $\Disc\Phi_{\bm}$ is a constant multiple of a power of $G_0$, and comparing degrees shows that the power is the first. Clearing denominators, we obtain $G$ and $c_G$.
\end{proof}

From now on we assume that~\ref{G2} holds. Let $M_2^\ast=4^n\operatorname{adj}(M_2)=4^n(\det M_2)M_2^{-1}$, which has entries in $\fo$ since $4^n\operatorname{adj}(M_2)=2^{n+1}\operatorname{adj}(2M_2)$, and put $Q_2^\ast(\bx)=\bx^tM_2^\ast\bx$. For odd $n$ put $\epsilon_n=(-1)^{(n-1)/2}$. Fix a non-zero $\Delta_X\in\fo$ divisible by
\[
2\,D_K\det(2M_1)\det(2M_2)\,c_G\,\Disc\mathcal{P}.
\]
A prime ideal $\fp$ is called \emph{good} if $\fp\nmid\Delta_X$, and \emph{bad} otherwise. If $\fp$ is good, then the reductions modulo $\fp$ of $M_1$, $M_2$ and $\mathcal{P}$ satisfy the conclusions of Lemma~\ref{lem:geometry} over $\FF_{\fp}=\fo/\fp$: by Lemma~\ref{lem:smoothcrit}, $X$ has non-singular reduction modulo $\fp$; $\operatorname{rank}M(\bs)\geqslant n-1$ modulo $\fp$ for every $\bs\neq\mathbf0$ over $\bar\FF_{\fp}$; $\det M_1\det M_2\not\equiv0$; and, for $\bm\in\fo^n$, if $\fp\nmid G(\bm)$ then $\Phi_{\bm}$ has $n-1$ distinct roots in $\PP^1(\bar\FF_{\fp})$, so that $X\cap\{\bm\cdot\bx=0\}$ is non-singular modulo $\fp$.

\subsection{Units and lattices}
Let $\ell:V^\times\to\RR^{r_1+r_2}$ be the map $\ell(v)=(c_l\log|v^{(l)}|)_l$. By Dirichlet's unit theorem, $\ell(U_K)$ is a lattice of full rank in the hyperplane $H=\{\mathbf e:\sum_le_l=0\}$, and the kernel of $\ell|_{U_K}$ is the group of roots of unity.

\begin{lemma}\label{lem:units}
For every $v\in V$ with $\Nm(v)\neq0$ there is $u\in U_K$ such that $|(uv)^{(l)}|\asymp|\Nm v|^{1/d}$ for all $l$.
\end{lemma}

\begin{proof}
This is~\cite[Lemma~2.1]{BV14}; we recall the proof. Since $\ell(U_K)$ is a full lattice in $H$, there is a constant $R_K$ such that every point of $H$ lies within distance $R_K$ of $\ell(U_K)$. The vector $\mathbf e=\ell(v)-d^{-1}\log|\Nm v|\,(c_l)_l$ lies in $H$, since $\sum_lc_l=d$ and $\sum_l\ell(v)_l=\log|\Nm v|$. Choose $u\in U_K$ with $|\ell(u)+\mathbf e|\leqslant R_K$. Then $\ell(uv)=\ell(u)+\mathbf e+d^{-1}\log|\Nm v|\,(c_l)_l$, so that $|c_l\log|(uv)^{(l)}|-c_ld^{-1}\log|\Nm v||\leqslant R_K$ for every $l$.
\end{proof}

We call $\bv\in V^m$ \emph{balanced} if $|\bv^{(l)}|>0$ for all $l$ and $|\bv^{(l)}|\asymp\cN(\bv)^{1/d}$ for all $l$, with implied constants depending only on $K$ and fixed once and for all.

\begin{lemma}\label{lem:balance}
Let $\bv\in V^m$ with $|\bv^{(l)}|>0$ for all $l$. Then $u\bv$ is balanced for some $u\in U_K$. If $\bv\in\fo^m\setminus\{\mathbf0\}$ is balanced, then $\cN(\bv)\geqslant1$, $|\bv^{(l)}|\gg1$ for all $l$, and $\cH(\bv)\asymp\cN(\bv)$.
\end{lemma}

\begin{proof}
Apply Lemma~\ref{lem:units} to the element $w\in V$ with $w^{(l)}=|\bv^{(l)}|$, which has $\Nm w=\cN(\bv)$, and note that $|(u\bv)^{(l)}|=|u^{(l)}||\bv^{(l)}|$. If $\bv\in\fo^m\setminus\{\mathbf0\}$, some coordinate $v_i$ is non-zero, and then $\cN(\bv)\geqslant\prod_l|v_i^{(l)}|^{c_l}=|\N_{K/\QQ}v_i|\geqslant1$. If $\bv$ is also balanced, then $|\bv^{(l)}|\gg\cN(\bv)^{1/d}\gg1$, so that $\max(1,|\bv^{(l)}|)\asymp|\bv^{(l)}|$ for every $l$, and hence $\cH(\bv)\asymp\cN(\bv)$.
\end{proof}

\begin{lemma}\label{lem:boxcount}
Let $\fb$ be a non-zero fractional ideal and let $B\in V$ with $B^{(l)}>0$ for all $l$. Then
\[
\#\{c\in\fb\setminus\{0\}:|c^{(l)}|\leqslant B^{(l)}\ \forall l\}\ll\frac{\Nm(B)}{\N\fb},
\]
and the set is empty if $\Nm(B)<\N\fb$.
\end{lemma}

\begin{proof}
Every $c\in\fb\setminus\{0\}$ satisfies $|\N c|\geqslant\N\fb$, while every $c$ in the box satisfies $|\N c|\leqslant\Nm(B)$. This proves the second assertion, so we may assume that $\Nm(B)\geqslant\N\fb$. Choose $u\in U_K$ as in Lemma~\ref{lem:units} for $B$, and put $R=\Nm(B)^{1/d}$. Multiplication by $u$ preserves $\fb$ and maps the box into the cube $\{\ip c\leqslant C_KR\}$. Let $\lambda_1\leqslant\dots\leqslant\lambda_d$ be the successive minima of the lattice $\fb$ with respect to this cube. By the inequality of arithmetic and geometric means, every $c\in\fb\setminus\{0\}$ has $\ip c^d\geqslant|\N c|\geqslant\N\fb$, so $\lambda_1\geqslant(\N\fb)^{1/d}$. By Minkowski's second theorem $\lambda_1\cdots\lambda_d\asymp\operatorname{covol}(\fb)\asymp\N\fb$, and therefore $\lambda_i\asymp(\N\fb)^{1/d}$ for all $i$. The number of lattice points in a cube of side $R$ is $\ll\prod_i(1+R/\lambda_i)$, which is $\ll R^d/\N\fb$ because $R\geqslant(\N\fb)^{1/d}$.
\end{proof}

\subsection{Additive characters}
An additive character modulo a non-zero integral ideal $\fb$ is a character of $\fo/\fb$. It is \emph{primitive} if it is not a character modulo any proper divisor of $\fb$. We call $\gamma\in\fb^{-1}\fd^{-1}$ \emph{primitive} if $\gamma\notin\fc^{-1}\fd^{-1}$ for every divisor $\fc\mid\fb$ with $\fc\neq\fb$. A sum $\sumstar_{\gamma\bmod\fb}$ runs over the primitive classes in $\fb^{-1}\fd^{-1}/\fd^{-1}$.

\begin{lemma}\label{lem:primchar}
The map $\gamma\mapsto(\alpha\mapsto\ee(\gamma\alpha))$ is an isomorphism from $\fb^{-1}\fd^{-1}/\fd^{-1}$ onto the group of additive characters modulo $\fb$, and it maps the primitive classes onto the primitive characters. More generally, for integral ideals $\fc\mid\fb$, the classes of $\fb^{-1}\fd^{-1}/\fd^{-1}$ that give characters trivial on $\fc$ are those of $\fc^{-1}\fd^{-1}/\fd^{-1}$.
\end{lemma}

\begin{proof}
This is standard; compare~\cite[\S2.3]{BV14}. For $\gamma\in\fb^{-1}\fd^{-1}$ and $\alpha\in\fb$ we have $\gamma\alpha\in\fd^{-1}$, and $\Tr(\fd^{-1})\subseteq\ZZ$, so $\ee(\gamma\alpha)=1$. Thus the map is well defined. The character attached to $\gamma$ is trivial on an ideal $\fc\supseteq\fb$ if and only if $\Tr(\gamma\fc)\subseteq\ZZ$, that is, $\gamma\fc\subseteq\fd^{-1}$, that is, $\gamma\in\fc^{-1}\fd^{-1}$. Taking $\fc=\fo$ shows that the kernel is $\fd^{-1}$. Both groups have order $\N\fb$, so the map is an isomorphism. The statement about primitivity follows from the general statement.
\end{proof}

For an integral ideal $\fa$ and a character $\sigma$ modulo $\fb$ we put $\sigma(\fa)=\sigma(\alpha)$, where $\alpha\in\fo$ is any element with $\ord_{\fp}\alpha=\ord_{\fp}\fa$ for all $\fp\mid\fb$. Such $\alpha$ exist by the Chinese remainder theorem. For $\fa=(0)$ we put $\sigma(\fa)=1$.

\begin{lemma}\label{lem:orth}
Let $\fa,\fb,\fc$ be integral ideals with $\fb,\fc\neq(0)$.
\begin{enumerate}[label=\textup{(\roman*)}]
\item $\sum_{\sigma\bmod\fb}\sigma(\fa)=\N\fb$ if $\fb\mid\fa$, and $=0$ otherwise.
\item $\sum_{\sigma\bmod\fc}\sigma(\fa)=\sum_{\fb\mid\fc}\sumstar_{\sigma\bmod\fb}\sigma(\fa)$.
\item The sum $c_{\fb}(\fa):=\sumstar_{\sigma\bmod\fb}\sigma(\fa)$ does not depend on the choice of $\alpha$, and
\[
c_{\fb}(\fa)=\sum_{\fe\mid\fb,\ \fe\mid\fa}\mu(\fb\fe^{-1})\N\fe.
\]
\end{enumerate}
\end{lemma}

\begin{proof}
Parts~(i) and~(ii) are~\cite[(2.9), (2.10)]{BV14}, and the independence of~$\alpha$ in~(iii) is~\cite[Lemma~2.4]{BV14}; the formula in~(iii) is a consequence. (i) By orthogonality of characters on $\fo/\fb$, the sum equals $\N\fb$ if $\alpha\in\fb$ and $0$ otherwise. Since $\ord_{\fp}\alpha=\ord_{\fp}\fa$ for all $\fp\mid\fb$, we have $\alpha\in\fb$ if and only if $\fb\mid\fa$. For $\fa=(0)$ the sum is $\N\fb$ by definition. (ii) Each character modulo $\fc$ is a primitive character modulo exactly one divisor $\fb\mid\fc$, namely its conductor. (iii) By (ii) and M\"obius inversion over the divisors of $\fb$, $\sumstar_{\sigma\bmod\fb}\sigma(\fa)=\sum_{\fe\mid\fb}\mu(\fb\fe^{-1})\sum_{\sigma\bmod\fe}\sigma(\fa)$. By (i) the inner sum is $\N\fe\cdot\mathbf 1[\fe\mid\fa]$, which is independent of $\alpha$.
\end{proof}

\subsection{Poisson summation and sums over ideals}
\begin{lemma}\label{lem:poisson-lattice}
Let $\Lambda\subset V^m$ be a lattice, let $\Lambda^\ast=\{\bt\in V^m:\Tr(\bt\cdot\bv)\in\ZZ\ \forall\bv\in\Lambda\}$ be its dual, and let $F$ be a Schwartz function on $V^m$. Then, for $\ba\in V^m$,
\[
\sum_{\bv\in\Lambda}F(\bv+\ba)=\frac1{\operatorname{covol}(\Lambda)}\sum_{\bt\in\Lambda^\ast}\hat F(\bt)\ee(\bt\cdot\ba).
\]
In particular, for a non-zero integral ideal $\fb$ and $\ba\in\fo^n$,
\[
\sum_{\bx\equiv\ba\bmod\fb}F(\bx)=\frac{2^{r_2n}}{D_K^{n/2}(\N\fb)^n}\sum_{\bm\in\hat\fb^{\,n}}\ee(\bm\cdot\ba)\hat F(\bm).
\]
\end{lemma}

\begin{proof}
Identify $V^m$ with $\RR^{dm}$ through the real and imaginary parts of the coordinates, so that $dv$ becomes Lebesgue measure. Then $(\bt,\bv)\mapsto\Tr(\bt\cdot\bv)$ is a non-degenerate bilinear form $B(\bt,\bv)=\bt^tJ\bv$ with $J$ invertible. The usual Poisson summation formula on $\RR^{dm}$ gives $\sum_{\Lambda}F(\bv+\ba)=\operatorname{covol}(\Lambda)^{-1}\sum_{\boldsymbol\eta\in\Lambda'}\mathcal F F(\boldsymbol\eta)e^{2\pi i\boldsymbol\eta\cdot\ba}$, where $\Lambda'$ is the standard dual lattice and $\mathcal F$ is the standard Fourier transform. The substitution $\boldsymbol\eta=J^t\bt$ maps $\Lambda^\ast$ bijectively onto $\Lambda'$, and $\mathcal FF(J^t\bt)=\hat F(\bt)$. For the second statement take $\Lambda=\fb^n$, whose covolume is $(2^{-r_2}\sqrt{D_K}\N\fb)^n$ and whose dual is $\hat\fb^{\,n}$.
\end{proof}

We shall use the following facts about the Dedekind zeta function $\zeta_K(s)=\sum_{\fa}\N\fa^{-s}$. It extends to a meromorphic function on $\CC$ whose only pole is a simple pole at $s=1$ with residue $\Delta_K$, where
\begin{equation}\label{eq:classnumber}
\Delta_K=\frac{2^{r_1}(2\pi)^{r_2}h_KR_K}{w_K\sqrt{D_K}},
\end{equation}
and $R_K$ is the regulator. In every vertical strip $\sigma_1\leqslant\Re s\leqslant\sigma_2$ away from $s=1$, it satisfies $\zeta_K(s)\ll(1+|\Im s|)^{A}$ for some $A=A(\sigma_1,\sigma_2,d)$. Finally, it follows from the functional equation that $\zeta_K$ has a zero of order $r_1+r_2-1$ at $s=0$, with
\begin{equation}\label{eq:zeta0}
\zeta_K(s)=-\frac{h_KR_K}{w_K}s^{r_1+r_2-1}+O(|s|^{r_1+r_2}),
\end{equation}
a zero of order $r_1+r_2$ at each negative even integer, and a zero of order $r_2$ at each negative odd integer.

\begin{lemma}\label{lem:idealsum}
Let $f\in C_c^\infty((0,\infty))$ and $R>0$. For every $N\geqslant0$,
\[
\sum_{(0)\neq\fb\subseteq\fo}f(R\N\fb)=\frac{\Delta_K}R\int_0^\infty f(t)\,dt+O_{N,f}(R^N).
\]
\end{lemma}

\begin{proof}
This is~\cite[Lemma~2.6]{BV14}; we give a short proof. Let $\tilde f(s)=\int_0^\infty f(t)t^{s-1}\,dt$ be the Mellin transform. Since $f$ has compact support in $(0,\infty)$, $\tilde f$ is entire, and integration by parts $k$ times gives
\[
\tilde f(s)=\frac{(-1)^k}{s(s+1)\cdots(s+k-1)}\int_0^\infty f^{(k)}(t)t^{s+k-1}\,dt,
\]
so that $\tilde f$ decays faster than any power of $|\Im s|$, uniformly in vertical strips. By Mellin inversion and absolute convergence on $\Re s=2$,
\[
\sum_{\fb} f(R\N\fb)=\frac1{2\pi i}\int_{(2)}\tilde f(s)R^{-s}\zeta_K(s)\,ds.
\]
Move the line of integration to $\Re s=-N$. The only pole crossed is at $s=1$, with residue $\tilde f(1)\Delta_KR^{-1}$. On the new line the integrand is $\ll R^N|\tilde f(-N+it)|(1+|t|)^A$, which is integrable.
\end{proof}

\subsection{Counting lemmas}
\begin{lemma}\label{lem:count-trivial}
\begin{enumerate}[label=\textup{(\roman*)}]
\item For $H\geqslant1$, $\#\{\bx\in\fo^m:\ip\bx\leqslant H\}\ll H^{dm}$.
\item Let $Z\subseteq\mathbb A^m_K$ be a closed subvariety of dimension $k$ defined by forms. Then $\#\{\bx\in\fo^m\cap Z(K):\ip\bx\leqslant H\}\ll_ZH^{dk}$.
\item For $M\geqslant1$, the number of balanced $\bv\in\fo^m$ with $\cH(\bv)\leqslant M$ is $\ll M^m$. If $Z$ is as in~\textup{(ii)}, the number of such $\bv$ lying on $Z$ is $\ll M^k$.
\end{enumerate}
\end{lemma}

\begin{proof}
(i) Fix a bounded fundamental domain $F$ for $\fo\subset V$. The translates $\alpha+F$, for $\alpha\in\fo$ with $\ip\alpha\leqslant H$, are disjoint and lie in $\{\ip v\leqslant H+C\}$, whose volume is $\ll H^d$. This gives the case $m=1$, and the general case follows. (ii) Since $Z$ is a cone of dimension $k$, Noether normalisation provides a linear map $\pi:\mathbb A^m\to\mathbb A^k$, which we may take to have coefficients in $\fo$, whose restriction to $Z$ is finite. Its fibres on $Z$ have at most $\deg Z$ points, and $\ip{\pi(\bx)}\ll\ip\bx$. So (i) with $m=k$ applies. (iii) A balanced $\bv$ has $\ip\bv\ll\cN(\bv)^{1/d}\leqslant\cH(\bv)^{1/d}\leqslant M^{1/d}$, so (i) and (ii) apply with $H\asymp M^{1/d}$.
\end{proof}

A subset $\Omega\subseteq\mathbb A^N(K)$ is \emph{thin} if it is contained in a finite union of sets of the following two types: the $K$-points of a proper Zariski-closed subset, and the images $\pi(Y(K))$, where $Y$ is irreducible of dimension $N$ and $\pi:Y\to\mathbb A^N$ is a dominant morphism of degree at least $2$; see~\cite[Chap.~9]{Ser97}. The following counting result is due to S.~D. Cohen; see~\cite[Chap.~13]{Ser97}.

\begin{lemma}\label{lem:thin}
If $\Omega\subseteq\mathbb A^N(K)$ is thin, then $\#\{\bx\in\Omega\cap\fo^N:\ip\bx\leqslant H\}\ll_\Omega H^{d(N-1/2)}\log H$ for $H\geqslant2$.
\end{lemma}

\begin{lemma}\label{lem:count-BV}
Let $F\in\fo[x_1,\dots,x_N]$ be an absolutely irreducible form of degree $D\geqslant2$. Then
\[
\#\{\bx\in\fo^N:\ip\bx\leqslant H,\ F(\bx)=0\}\ll_F H^{d(N-3/2)}\log H\qquad(H\geqslant2).
\]
Let $F\in\fo[x_1,\dots,x_N]$ be a non-zero form which is not a constant multiple of the square of a form, and let $a\in K^\times$. Then
\[
\#\{\bx\in\fo^N:\ip\bx\leqslant H,\ F(\bx)\in aK^{\times2}\}\ll_{F,a}H^{d(N-1/2)}\log H\qquad(H\geqslant2).
\]
\end{lemma}

\begin{proof}
The first assertion is, with $\log H$ in place of $H^\varepsilon$, the case of trivial congruence conditions of~\cite[Lemma~2.5]{BV14}; the second is not stated there. We deduce both from Lemma~\ref{lem:thin}. For the first statement, choose $\mathbf e\in\fo^N$ with $F(\mathbf e)\neq0$, and a matrix $T\in M_N(\fo)$ with $\det T\neq0$ whose last column is $\mathbf e$. The form $F'(\by)=F(T\by)$ has coefficient $F(\mathbf e)\neq0$ at $y_N^D$. If $\bx\in\fo^N$ with $\ip\bx\leqslant H$ and $F(\bx)=0$, then $\by=(\det T)T^{-1}\bx$ lies in $\fo^N$, satisfies $\ip\by\ll H$ and $F'(\by)=0$, and $\bx\mapsto\by$ is injective. So it suffices to treat $F'$. Let $Y=\{F'=0\}\subset\mathbb A^N$, which is irreducible, and let $\pi(\by)=(y_1,\dots,y_{N-1})$. Since $F'$ is, up to the non-zero constant $F(\mathbf e)$, monic of degree $D$ in $y_N$, the map $\pi:Y\to\mathbb A^{N-1}$ is finite, surjective and of degree $D\geqslant2$. Hence $\pi(Y(K))$ is thin in $\mathbb A^{N-1}(K)$. Each of its points has at most $D$ preimages on $Y$, and Lemma~\ref{lem:thin} with $N-1$ in place of $N$ gives the bound.

For the second statement, let $Y_a=\{(\bx,y):ay^2=F(\bx)\}\subset\mathbb A^{N+1}$. Since $F/a$ is not a square in $K[\bx]$, the polynomial $ay^2-F(\bx)$ is irreducible, so $Y_a$ is irreducible, and the projection $Y_a\to\mathbb A^N$ has degree $2$. Its image on $K$-points is $\{\bx:F(\bx)\in aK^{\times2}\}\cup\{F=0\}$, which is therefore thin, and Lemma~\ref{lem:thin} applies.
\end{proof}

\begin{lemma}\label{lem:Hint}
Let $A\geqslant1$. Then $\operatorname{vol}\{v\in V:\cH(v)\leqslant A,\ \ip v\leqslant A\}\ll A(\log2A)^{r_1+r_2-1}$, and for every $\alpha<-1$ we have $\int_V\cH(v)^\alpha\,dv\ll_\alpha1$ and $\int_{\cH(v)\geqslant A}\cH(v)^\alpha\,dv\ll_\alpha A^{\alpha+1}(\log2A)^{r_1+r_2-1}$.
\end{lemma}

\begin{proof}
Split $V$ according to the set $S=\{l:|v^{(l)}|>1\}$. On the piece indexed by $S$, the coordinates $v^{(l)}$ with $l\notin S$ range over a set of bounded volume, and $\cH(v)=\prod_{l\in S}u_l$ with $u_l=|v^{(l)}|^{c_l}\geqslant1$. In the variables $u_l$ the measure is $\asymp\prod_{l\in S}du_l$ (up to an angular factor at complex places). Hence everything reduces to the elementary estimates $\int_{u_l\geqslant1,\ \prod u_l\leqslant A}du\ll A(\log2A)^{|S|-1}$ and $\int_{u_l\geqslant1,\ \prod u_l\geqslant A}\prod u_l^\alpha\,du\ll A^{\alpha+1}(\log2A)^{|S|-1}$, both of which follow by induction on $|S|$.
\end{proof}

\begin{lemma}\label{lem:smoothsupp}
Let $\fg$ be a non-zero integral ideal, $A\geqslant1$ and $x\geqslant1$. Then
\[
\#\{\fb:\fb\mid\fg^\infty,\ \N\fb\leqslant x\}\ll_\varepsilon(x\N\fg)^\varepsilon,\qquad\sum_{\N\fb\leqslant x}\frac{A^{\omega(\fb)}}{\N\rad\fb}\ll_{A,\varepsilon}x^\varepsilon.
\]
Moreover $A^{\omega(\fb)}\ll_{A,\varepsilon}\N\fb^\varepsilon$, where $\rad\fb=\prod_{\fp\mid\fb}\fp$.
\end{lemma}

\begin{proof}
For the first bound, the count is at most $x^\varepsilon\sum_{\fb\mid\fg^\infty}\N\fb^{-\varepsilon}=x^\varepsilon\prod_{\fp\mid\fg}(1-\N\fp^{-\varepsilon})^{-1}$. The factors with $\N\fp>2^{1/\varepsilon}$ are at most $2$, and there are $O_\varepsilon(1)$ others. So the product is $\ll_\varepsilon2^{\omega(\fg)}\ll_\varepsilon\N\fg^\varepsilon$, using the last assertion. For the second bound, the sum is at most $x^\varepsilon\prod_{\fp}(1+A\N\fp^{-1}\sum_{k\geqslant1}\N\fp^{-k\varepsilon})$, and this product converges because its factors are $1+O_\varepsilon(A\N\fp^{-1-\varepsilon})$ for $\N\fp$ large. Finally, $A^{\omega(\fb)}\leqslant\prod_{\fp\mid\fb}\max(A,\N\fp^{\varepsilon})\N\fp^{-\varepsilon}\cdot\N\fb^\varepsilon$, and the product is bounded by $\prod_{\N\fp\leqslant A^{1/\varepsilon}}A$.
\end{proof}

\subsection{The spectrum of the pencil}
Let $Z\subset\PP^1(\CC)$ be the set of roots of $\mathcal{P}$. By~\ref{G1} it has exactly $n$ elements, and by~\ref{G2} it contains neither $[1:0]$ nor $[0:1]$. We use the chordal metric on $\PP^1(\CC)$.

\begin{lemma}\label{lem:spectrum}
Let $K_l\in\{\RR,\CC\}$ and $\bs\in K_l^2\setminus\{\mathbf0\}$, and let $\sigma_1(\bs)\leqslant\dots\leqslant\sigma_n(\bs)$ be the singular values of $M(\bs)$. Then $\sigma_2(\bs)\asymp\sigma_n(\bs)\asymp|\bs|$ and $\sigma_1(\bs)\asymp|\bs|\operatorname{dist}([\bs],Z)$.
\end{lemma}

\begin{proof}
All quantities are homogeneous of degree one in $\bs$, so we may assume $|\bs|=1$. Since $M_1$ and $M_2$ are linearly independent, $\sigma_n(\bs)=\|M(\bs)\|\asymp1$. The function $\bs\mapsto\sigma_2(\bs)$ is continuous on the compact set $\{|\bs|=1\}$, and it is positive because $\operatorname{rank}M(\bs)\geqslant n-1$ by~\ref{G1}; hence $\sigma_2\gg1$. Since $|\mathcal{P}(\bs)|=|\det M(\bs)|=\prod_j\sigma_j$, it follows that $\sigma_1\asymp|\mathcal{P}(\bs)|$. Finally, $\mathcal{P}$ is a binary form of degree $n$ with $n$ distinct roots, none at infinity, so $|\mathcal{P}(\bs)|\asymp\prod_{z\in Z}\operatorname{dist}([\bs],z)$ for $|\bs|=1$. Near a root $z$ the other factors are $\asymp1$, and away from all roots every factor is $\asymp1$. Hence $|\mathcal{P}(\bs)|\asymp\operatorname{dist}([\bs],Z)$.
\end{proof}

\section{Two $\delta$-symbols over $K$}\label{sec:delta}

In this section we describe the two $\delta$-symbols that we shall use. The first is the ideal $\delta$-symbol of Browning and Vishe~\cite{BV14}, which will detect the second equation. The second is a $\delta$-symbol over elements, which will detect the first equation. Both are instances of Heath-Brown's identity~\cite[Theorem~1]{HB96}: if $\mathcal D(a)$ is the set of divisors of $a\neq0$, then $c\mapsto a/c$ is a bijection of $\mathcal D(a)$, and so any weight of the shape $w(c)-w(a/c)$ sums to zero over $\mathcal D(a)$.

\subsection{The ideal $\delta$-symbol}
Fix a non-negative function $w\in C^\infty(\RR)$ supported in $[\frac12,1]$ with $\int w=1$. For $x>0$ and $y\in\RR$ define
\begin{equation}\label{eq:hdef}
h(x,y)=\frac1{\Delta_K}\sum_{(0)\neq\fz\subseteq\fo}\frac1{x\N\fz}\Bigl\{w(x\N\fz)-w\Bigl(\frac{|y|}{x\N\fz}\Bigr)\Bigr\}.
\end{equation}
For fixed $(x,y)$ only finitely many $\fz$ contribute, since both terms vanish once $x\N\fz>\max(1,2|y|)$. The function $h$ depends on $y$ only through $|y|$, and we shall see in Lemma~\ref{lem:h1} that it is constant in $y$ for $|y|\leqslant x/2$. In particular $h$ is smooth on $(0,\infty)\times\RR$.

\begin{theorem}\label{thm:BV}
Let $Q\geqslant1$ and define $c_Q$ by
\[
c_Q^{-1}=\frac1{\Delta_KQ^d}\sum_{\fc\neq(0)}w\Bigl(\frac{\N\fc}{Q^d}\Bigr).
\]
Then, for every integral ideal $\fa$,
\[
\delta_K(\fa)=\frac{c_Q}{Q^{2d}}\sum_{(0)\neq\fb\subseteq\fo}\ \sumstar_{\sigma\bmod\fb}\sigma(\fa)\,h\Bigl(\frac{\N\fb}{Q^d},\frac{\N\fa}{Q^{2d}}\Bigr).
\]
Moreover $c_Q=1+O_N(Q^{-N})$ for every $N$. Finally $h(x,y)\ll x^{-1}$, and $h(x,y)=0$ unless $x\leqslant\max(1,2|y|)$.
\end{theorem}

\begin{proof}
This is~\cite[Theorem~1.2]{BV14}, the function $h$ being the one defined in~\cite[\S3]{BV14}; we include the short proof. Let $\fa\neq(0)$. The map $\fc\mapsto\fa\fc^{-1}$ is a bijection of the set of divisors of $\fa$, so
\[
\sum_{\fc\mid\fa}\Bigl\{w\Bigl(\frac{\N\fc}{Q^d}\Bigr)-w\Bigl(\frac{\N\fa}{Q^d\N\fc}\Bigr)\Bigr\}=0.
\]
If $\fa=(0)$, every non-zero $\fc$ divides $\fa$, and $w(\N\fa/(Q^d\N\fc))=w(0)=0$, so the left-hand side equals $\sum_{\fc} w(\N\fc/Q^d)=c_Q^{-1}\Delta_KQ^d$. In either case,
\[
\delta_K(\fa)=\frac{c_Q}{\Delta_KQ^d}\sum_{\fc\mid\fa}\Bigl\{w\Bigl(\frac{\N\fc}{Q^d}\Bigr)-w\Bigl(\frac{\N\fa}{Q^d\N\fc}\Bigr)\Bigr\}.
\]
By Lemma~\ref{lem:orth}(i) and~(ii),
\[
\mathbf 1[\fc\mid\fa]=\frac1{\N\fc}\sum_{\sigma\bmod\fc}\sigma(\fa)=\frac1{\N\fc}\sum_{\fb\mid\fc}\ \sumstar_{\sigma\bmod\fb}\sigma(\fa).
\]
Inserting this and writing $\fc=\fb\fz$, we obtain
\[
\delta_K(\fa)=\frac{c_Q}{\Delta_KQ^d}\sum_{\fb}\ \sumstar_{\sigma\bmod\fb}\sigma(\fa)\sum_{\fz}\frac1{\N\fb\N\fz}\Bigl\{w\Bigl(\frac{\N\fb\N\fz}{Q^d}\Bigr)-w\Bigl(\frac{\N\fa}{Q^d\N\fb\N\fz}\Bigr)\Bigr\}.
\]
With $x=\N\fb/Q^d$ and $y=\N\fa/Q^{2d}$ we have $(\N\fb\N\fz)^{-1}=Q^{-d}(x\N\fz)^{-1}$ and $\N\fa/(Q^d\N\fb\N\fz)=y/(x\N\fz)$. So the inner sum is $Q^{-d}\Delta_Kh(x,y)$, which gives the identity. All sums are finite for fixed $\fa$, so the rearrangements are justified.

For the asymptotic formula for $c_Q$, apply Lemma~\ref{lem:idealsum} with $f=w$ and $R=Q^{-d}$: it gives $\sum_{\fc} w(\N\fc/Q^d)=\Delta_KQ^d+O_N(Q^{-N})$. For the support of $h$, note that if $x>\max(1,2|y|)$, then $x\N\fz>1$ and $|y|/(x\N\fz)<\frac12$ for every $\fz$, so both terms in~\eqref{eq:hdef} vanish. Finally, the first term in~\eqref{eq:hdef} is supported on $\frac12\leqslant x\N\fz\leqslant1$, and the second on $|y|\leqslant x\N\fz\leqslant2|y|$. On each of these ranges the summand is $O(1)$, and the number of ideals with norm in a dyadic interval $[T,2T]$ is $O(1+T)$. Since $T=(2x)^{-1}$ in the first case, the first term contributes $O(1+x^{-1})$. In the second case $T=|y|/x$, and the bound $(x\N\fz)^{-1}\leqslant|y|^{-1}$ gives a contribution $O(|y|^{-1}(1+|y|/x))$. Since $h=0$ unless $x\leqslant\max(1,2|y|)$, both are $O(x^{-1})$ on the support.
\end{proof}

\subsection{Properties of $h$}
The following lemmas are~\cite[Lemmas~3.1, 3.2 and (3.7), (3.8)]{BV14}. Write $h=h_1-h_2$, where $h_1$ and $h_2$ are the contributions of the two terms in~\eqref{eq:hdef}.

\begin{lemma}\label{lem:h1}
If $|y|\leqslant x/2$, then $h_2(x,y)=0$, so that $h(x,y)=h_1(x)$ does not depend on $y$. If moreover $x\geqslant1$, then $h(x,y)=0$.
\end{lemma}

\begin{proof}
If $|y|\leqslant x/2$, then $|y|/(x\N\fz)\leqslant\frac12$ for every $\fz$, and $w$ vanishes on $(-\infty,\frac12]$ (at $\frac12$ by continuity). If also $x\geqslant1$, then $x\N\fz\geqslant1$, and $w$ vanishes on $[1,\infty)$.
\end{proof}

\begin{lemma}\label{lem:h2}
Let $i,j\geqslant0$. If $x\leqslant1$ and $|y|\leqslant x/2$, then $\partial_x^ih(x,y)\ll_ix^{-i-1}$. If $|y|\geqslant x/2$, then $\partial_x^i\partial_y^jh(x,y)\ll_{i,j}x^{-i-1}|y|^{-j}$ for $y\neq0$.
\end{lemma}

\begin{proof}
By induction on $i$ one finds constants $d_{i,\tau}$ and $c_{i,j,\tau}$ such that
\begin{align*}
\partial_x^i\Bigl(\frac{w(xt)}{xt}\Bigr)&=x^{-i}\sum_{0\leqslant\tau\leqslant i}d_{i,\tau}\frac{(xt)^\tau w^{(\tau)}(xt)}{xt},\\
\partial_x^i\partial_y^j\Bigl(\frac{w(|y|/(xt))}{xt}\Bigr)&=\frac{x^{-i}|y|^{-j}}{xt}\sum_{0\leqslant\tau\leqslant i+j}c_{i,j,\tau}\Bigl(\frac{|y|}{xt}\Bigr)^{\tau}w^{(\tau)}\Bigl(\frac{|y|}{xt}\Bigr).
\end{align*}
Take $t=\N\fz$ and sum over $\fz$. In the first line the summand is supported where $xt\asymp1$ and is $O(x^{-i})$ there; the number of $\fz$ with $\N\fz\asymp x^{-1}$ is $O(1+x^{-1})$, which is $O(x^{-1})$ when $x\leqslant1$. In the second line the summand is supported where $xt\asymp|y|$ and is $O(x^{-i}|y|^{-j-1})$ there; the number of such $\fz$ is $O(1+|y|/x)$, which is $O(|y|/x)$ when $|y|\geqslant x/2$. This proves the second assertion, and also the first, since in that case only $h_1$ is present by Lemma~\ref{lem:h1}.
\end{proof}

\begin{lemma}\label{lem:h3}
For all $i,j,N\geqslant0$ and $x>0$, $y\neq0$,
\[
\partial_x^i\partial_y^jh(x,y)\ll_{i,j,N}x^{-i-j-1}\Bigl(x^N+\min\Bigl\{1,\Bigl(\frac x{|y|}\Bigr)^N\Bigr\}\Bigr).
\]
The term $x^N$ may be omitted when $j\geqslant1$.
\end{lemma}

\begin{proof}
If $|y|\leqslant x/2$, this follows from Lemmas~\ref{lem:h1} and~\ref{lem:h2} (with $\partial_yh=0$). If $x/2\leqslant|y|\leqslant x$, it follows from Lemma~\ref{lem:h2}. It remains to treat the range $|y|>x$, and there we may assume $N\geqslant1$.

Consider first $h_1$, which does not depend on $y$; so we may take $j=0$. With $g_\tau(z)=z^{\tau-1}w^{(\tau)}(z)$, the first formula in the proof of Lemma~\ref{lem:h2} gives $\partial_x^ih_1=\Delta_K^{-1}x^{-i}\sum_\tau d_{i,\tau}\sum_{\fz} g_\tau(x\N\fz)$. By Lemma~\ref{lem:idealsum} with $R=x$,
\[
\sum_{\fz} g_\tau(x\N\fz)=\frac{\Delta_K}x\int_0^\infty g_\tau(z)\,dz+O_N(x^N).
\]
For $\tau\geqslant1$ we have $\int_0^\infty z^{\tau-1}w^{(\tau)}(z)\,dz=0$, by $\tau-1$ integrations by parts and $w^{(\tau-\tau')}$ having compact support in $(0,\infty)$. Since $d_{i,0}$ is the coefficient of $g_0$, we obtain
\[
\partial_x^ih_1(x)=x^{-i-1}\Bigl(d_{i,0}\int_0^\infty\frac{w(z)}z\,dz+O_N(x^{N+1})\Bigr).
\]
For $h_2$, put $u=|y|/(xt)$ in the second formula of the proof of Lemma~\ref{lem:h2}. Since $(xt)^{-1}=u/|y|$, the summand indexed by $\fz$, with $t=\N\fz$, equals
\[
x^{-i}|y|^{-j-1}\sum_\tau c_{i,j,\tau}u^{\tau+1}w^{(\tau)}(u)=x^{-i}|y|^{-j-1}\sum_\tau c_{i,j,\tau}F_\tau(R\N\fz),
\]
where $R=x/|y|$ and $F_\tau(z)=z^{-\tau-1}w^{(\tau)}(1/z)$, a smooth function supported in $[1,2]$. By Lemma~\ref{lem:idealsum}, $\sum_{\fz} F_\tau(R\N\fz)=\Delta_KR^{-1}\int_0^\infty F_\tau+O_N(R^N)$, and the substitution $z=1/u$ gives $\int_0^\infty F_\tau(z)\,dz=\int_0^\infty u^{\tau-1}w^{(\tau)}(u)\,du$. As before, this vanishes for $\tau\geqslant1$. Hence
\[
\partial_x^i\partial_y^jh_2=x^{-i-1}|y|^{-j}\Bigl(c_{i,j,0}\int_0^\infty\frac{w(u)}u\,du\Bigr)+O_N\bigl(x^{-i}|y|^{-j-1}(x/|y|)^N\bigr).
\]
The error term is $\ll x^{-i-j-1}(x/|y|)^N$, since $x<|y|$.

If $j\geqslant1$, then $h_1$ does not contribute, and $c_{i,j,0}=0$, because every derivative in $y$ falls on $w(|y|/(xt))$ and therefore produces a derivative of $w$. This gives the claim, without the term $x^N$. If $j=0$, the main terms of $\partial_x^ih_1$ and $\partial_x^ih_2$ are equal: both $d_{i,0}$ and $c_{i,0,0}$ are the coefficient $(-1)^ii!$ that arises when all the derivatives fall on the factor $x^{-1}$. So the main terms cancel, and what remains is $O(x^{-i-1}(x^{N+1}+(x/|y|)^N))$.
\end{proof}

The Fourier analysis of Section~\ref{sec:integrals} only uses the bounds of Lemma~\ref{lem:h3}, together with the fact that $h$ is constant in $y$ near $y=0$. Since we shall need to apply the same arguments to some derivatives of $h$, we isolate these properties.

\begin{definition}\label{def:admissible}
A smooth function $k:(0,\infty)\times\RR\to\CC$ is \emph{admissible} if $k(x,y)=k(x,|y|)$, if $k(x,y)$ does not depend on $y$ for $|y|\leqslant x/2$, and if $k$ satisfies the bounds of Lemma~\ref{lem:h3}.
\end{definition}

\begin{lemma}\label{lem:admissible}
The functions $h$, $x\partial_xh$ and $y\partial_yh$ are admissible.
\end{lemma}

\begin{proof}
For $h$ this is Lemmas~\ref{lem:h1} and~\ref{lem:h3}. Since $h$ is constant in $y$ for $|y|\leqslant x/2$, so are $x\partial_xh$ and $y\partial_yh=0$. Next, $\partial_x^i\partial_y^j(x\partial_xh)=x\partial_x^{i+1}\partial_y^jh+i\partial_x^i\partial_y^jh$, and each term satisfies the bound of Lemma~\ref{lem:h3}. Similarly $\partial_x^i\partial_y^j(y\partial_yh)=y\partial_x^i\partial_y^{j+1}h+j\partial_x^i\partial_y^jh$. The second term is fine. For the first, Lemma~\ref{lem:h3} with $N+1$ in place of $N$ gives $|y|x^{-i-j-2}\min\{1,(x/|y|)^{N+1}\}\leqslant x^{-i-j-1}\min\{1,(x/|y|)^N\}$, by considering $|y|\leqslant x$ and $|y|>x$ separately.
\end{proof}

\begin{lemma}\label{lem:hL1}
Let $k$ be admissible, let $R_0\geqslant1$ and $0<x\leqslant c$ for a constant $c$. Then
\[
\int_{\ip v\leqslant R_0}|k(x,\Nm v)|\,dv\ll_{R_0}(1+|\log x|)^{r_1+r_2}.
\]
\end{lemma}

\begin{proof}
This is the estimate in the proof of~\cite[Lemma~6.4]{BV14}. Let $\mathcal R_1$ be the set of $v$ with $\ip v\leqslant R_0$ and $|v^{(l)}|<x$ for some $l$. It has volume $O(x)$, and $|k|\ll x^{-1}$ by the admissibility bound with $i=j=N=0$, so $\mathcal R_1$ contributes $O(1)$. On the complement of $\mathcal R_1$, the bound with $N=2$ gives $|k(x,\Nm v)|\ll x+x^{-1}\min\{1,x^2/|\Nm v|^2\}$. The term $x$ contributes $O(1)$. The set of $v$ with $\ip v\leqslant R_0$ and $|\Nm v|\leqslant x$ has volume $\ll x(1+|\log x|)^{r_1+r_2-1}$, by the argument of Lemma~\ref{lem:Hint}; so it contributes $O((1+|\log x|)^{r_1+r_2-1})$. Finally, decomposing $\{|\Nm v|>x,\ \min_l|v^{(l)}|\geqslant x\}$ dyadically according to $|\Nm v|\asymp2^kx$, each piece has volume $\ll2^kx(1+|\log x|)^{r_1+r_2-1}$ and contributes $\ll2^{-k}(1+|\log x|)^{r_1+r_2-1}$. Summing over $k\geqslant0$ gives the claim.
\end{proof}

\subsection{Local zeta integrals and the main term}
The main term of the circle method comes from the following property of $h$, which says that $h(x,\Nm v)\,dv$ behaves like a multiple of the Dirac measure at $v=0$ as $x\to0$. It is~\cite[Lemma~4.1]{BV14}. We give a proof based on local zeta integrals, which makes the role of the zeros of $\zeta_K$ transparent.

For $f\in C_c^\infty(V)$ and $\Re s>-1$, put
\[
Z_f(s)=\int_Vf(v)|\Nm v|^s\,dv.
\]

\begin{lemma}\label{lem:zeta}
Let $f\in C_c^\infty(V)$ be supported in $\{\ip v\leqslant R_0\}$ and let $J\geqslant1$. Then $Z_f$ extends to a meromorphic function on $\Re s>-1-J$, whose poles are contained in $\{-1,-2,\dots,-J\}$. At $s=-1-j$ with $j$ odd the pole has order at most $r_2$, and at $s=-1-j$ with $j$ even it has order at most $r_1+r_2$. At $s=-1$,
\[
Z_f(s)=\frac{2^{r_1}\pi^{r_2}f(0)}{(s+1)^{r_1+r_2}}+O\bigl(|s+1|^{-r_1-r_2+1}\bigr).
\]
Moreover, for $\Re s>-1-J$ and $\operatorname{dist}(s,\ZZ)\geqslant\frac14$, we have $Z_f(s)\ll_{R_0,J}\|f\|_{C^{k}}(1+|\Im s|)^{-2}$ for some $k=k(J,d)$.
\end{lemma}

\begin{proof}
For each place $l$ fix a radial function $\chi_l\in C_c^\infty(K_l)$ with $\chi_l=1$ on $|v_l|\leqslant1$ and $\chi_l=0$ on $|v_l|\geqslant2$. For a smooth function $g$ of $v_l\in K_l$ and an integer $k\geqslant1$, define the Taylor operator
\[
T_l^kg(v_l)=\chi_l(v_l)\sum_{|\alpha|<k}\frac{\partial^\alpha g(0)}{\alpha!}v_l^\alpha,
\]
where at a real place $\alpha\in\ZZ_{\geqslant0}$ and $v_l^\alpha$ has its usual meaning, and at a complex place $\alpha=(a,b)$, $v_l^\alpha=v_l^a\bar v_l^b$, $\partial^\alpha=\partial_{v_l}^a\partial_{\bar v_l}^b$ and $|\alpha|=a+b$. Put $R_l^k=1-T_l^k$. By Taylor's theorem, $R_l^kg$ vanishes to order $k$ at $v_l=0$. Operators attached to different places act on different variables and commute, so
\[
f=\prod_l(T_l^k+R_l^k)f=\sum_{S}f_S,\qquad f_S=\prod_{l\in S}T_l^k\prod_{l\notin S}R_l^kf,
\]
where $S$ runs over the subsets of the set of places. Expanding the operators $T_l^k$ for $l\in S$, $f_S$ is a sum over multi-indices $(\alpha_l)_{l\in S}$ of terms
\[
\prod_{l\in S}\chi_l(v_l)\frac{v_l^{\alpha_l}}{\alpha_l!}\cdot g_{S,\alpha}(v_{S^c}),
\]
where $g_{S,\alpha}=\prod_{l\notin S}R_l^k\bigl[(\partial^{\alpha}f)(0_S,\cdot)\bigr]$ is a smooth function of the variables $v_{S^c}$ at the places outside $S$. It is supported in a fixed ball and satisfies $|g_{S,\alpha}(v_{S^c})|\ll\|f\|_{C^{k'}}\prod_{l\notin S}\min(1,|v_l|^k)$ for some $k'$ depending on $k$ and $d$. The same holds, with $k-1$ in place of $k$, for its first derivatives in the radial variables $|v_l|$.

Since $|\Nm v|^s=\prod_l|v_l|^{c_ls}$, the zeta integral of each term factorises as
\[
\prod_{l\in S}\zeta_{l,\alpha_l}(s)\cdot\int g_{S,\alpha}(v_{S^c})\prod_{l\notin S}|v_l|^{c_ls}\,dv_{S^c},\qquad\zeta_{l,\alpha}(s)=\int_{K_l}\chi_l(v)\frac{v^\alpha}{\alpha!}|v|^{c_ls}\,dv.
\]
The last integral converges absolutely and is holomorphic for $\Re s>-1-k/2$. One radial integration by parts at each place $l\notin S$ shows that it is $O(\|f\|_{C^{k'+1}}(1+|\Im s|)^{-1})$ in $\Re s>-1-(k-1)/2$; two give $(1+|\Im s|)^{-2}$.

Now we compute $\zeta_{l,\alpha}$. At a real place, by symmetry, $\zeta_{l,\alpha}(s)=0$ for odd $\alpha$, while for even $\alpha$
\[
\zeta_{l,\alpha}(s)=\frac2{\alpha!}\int_0^\infty\chi_l(r)r^{\alpha+s}\,dr=-\frac2{\alpha!}\cdot\frac1{s+\alpha+1}\int_0^\infty\chi_l'(r)r^{\alpha+s+1}\,dr.
\]
Since $\chi_l'$ is supported in $[1,2]$, the last integral is entire and decays rapidly in $|\Im s|$. So $\zeta_{l,\alpha}$ has a single simple pole, at $s=-1-\alpha$, and for $\alpha=0$ its residue is $-2\int_0^\infty\chi_l'(r)\,dr=2$. At a complex place, write $v=re^{i\phi}$ and $\alpha=(a,b)$. The integral over $\phi$ of $e^{i(a-b)\phi}$ vanishes unless $a=b$, and then
\[
\zeta_{l,\alpha}(s)=\frac{2\pi}{a!^2}\int_0^\infty\chi_l(r)r^{2a+2s+1}\,dr=-\frac{\pi}{a!^2}\cdot\frac1{s+a+1}\int_0^\infty\chi_l'(r)r^{2a+2s+2}\,dr,
\]
which has a single simple pole, at $s=-1-a$, with residue $\pi$ when $a=0$.

Choose $k>2J+4$. In $\Re s>-1-J$, each term is therefore meromorphic, with poles only at $s=-1-j$, $0\leqslant j<J$. The pole at $s=-1-j$ comes only from the factors $\zeta_{l,\alpha_l}$ with $l\in S$ that have a pole there: at real places this requires $j$ to be even, while at complex places any $j$ is possible. This gives the stated bounds for the orders. At $s=-1$ a pole of order $r_1+r_2$ can only come from $S$ equal to the set of all places and $\alpha=0$. That term is $f(0)\prod_l\chi_l(v_l)$, whose zeta integral is $f(0)\prod_l\zeta_{l,0}(s)$, with leading coefficient $2^{r_1}\pi^{r_2}f(0)$. The final bound follows from the decay of the factors $\zeta_{l,\alpha}$ and of the integrals over $v_{S^c}$.
\end{proof}

\begin{lemma}\label{lem:h4}
Let $f\in C_c^\infty(V)$ be supported in $\{\ip v\leqslant R_0\}$, and let $0<x\leqslant1$ and $N\geqslant1$. Then
\[
\int_Vf(v)h(x,\Nm v)\,dv=\frac{\sqrt{D_K}}{2^{r_2}}f(0)+O_{N,R_0}\bigl(x^N\|f\|_{C^k}\bigr)
\]
for some $k=k(N,d)$.
\end{lemma}

\begin{proof}
This is~\cite[Lemma~4.1]{BV14}. The proof given there is different; we use local zeta integrals. Write $G(T)=\int_Vf(v)w(|\Nm v|/T)\,dv$ for $T>0$, and
\[
\Psi(T)=\frac1T\Bigl\{w(T)\int_Vf-G(T)\Bigr\}.
\]
By~\eqref{eq:hdef},
\begin{equation}\label{eq:h4a}
\int_Vf(v)h(x,\Nm v)\,dv=\frac1{\Delta_K}\sum_{\fz}\Psi(x\N\fz).
\end{equation}
This sum is finite: $w(T)=0$ for $T\geqslant1$, and $G(T)=0$ once $T\geqslant2R_0^d$, since $|\Nm v|\leqslant R_0^d$ on the support of $f$. The function $\Psi$ is smooth on $(0,\infty)$. For $T<\frac12$ we have $\Psi(T)=-G(T)/T$, and $|G(T)|\ll\|f\|_\infty\operatorname{vol}\{\ip v\leqslant R_0:|\Nm v|\leqslant T\}\ll\|f\|_\infty T(1+|\log T|)^{r_1+r_2-1}$. The same bound holds for $(T\partial_T)^jG$, which has the same shape with $w$ replaced by $(-u\partial_u)^jw$. So the Mellin transform $\tilde\Psi(s)=\int_0^\infty\Psi(T)T^{s-1}\,dT$ converges absolutely for $\Re s>0$.

We compute $\tilde\Psi$ for $\Re s>0$. Let $M(s)=\int_0^\infty w(t)t^{s-1}\,dt$, which is entire, satisfies $M(1)=\int w=1$, and decays faster than any power of $|\Im s|$ in vertical strips. Then $\int_0^\infty T^{s-2}w(T)\,dT=M(s-1)$. By Fubini and the substitution $T=|\Nm v|/u$,
\[
\int_0^\infty T^{s-2}G(T)\,dT=\int_Vf(v)\int_0^\infty T^{s-2}w\Bigl(\frac{|\Nm v|}T\Bigr)\,dT\,dv=Z_f(s-1)\,M(1-s).
\]
Hence
\begin{equation}\label{eq:h4b}
\tilde\Psi(s)=M(s-1)\int_Vf-M(1-s)Z_f(s-1).
\end{equation}
On the line $\Re s=2$ this is integrable, so Mellin inversion and~\eqref{eq:h4a} give
\[
\int_Vf(v)h(x,\Nm v)\,dv=\frac1{\Delta_K}\cdot\frac1{2\pi i}\int_{(2)}\tilde\Psi(s)x^{-s}\zeta_K(s)\,ds,
\]
where the interchange of the sum over $\fz$ and the integral is justified by absolute convergence.

By Lemma~\ref{lem:zeta}, \eqref{eq:h4b} continues meromorphically to $\Re s>-N$. We move the line of integration to $\Re s=-N+\frac12$, and examine the possible poles of $\tilde\Psi(s)x^{-s}\zeta_K(s)$.

At $s=1$, $\zeta_K$ has a simple pole, but $\tilde\Psi(1)=M(0)\int f-M(0)Z_f(0)=0$ because $Z_f(0)=\int f$. So there is no pole.

At $s=0$, the first term of~\eqref{eq:h4b} is entire, and so is $\zeta_K$; so only the second term matters. By Lemma~\ref{lem:zeta}, $Z_f(s-1)$ has a pole of order at most $r_1+r_2$ at $s=0$, with leading coefficient $2^{r_1}\pi^{r_2}f(0)$, while $\zeta_K$ has a zero of order $r_1+r_2-1$, with leading coefficient $-h_KR_K/w_K$ by~\eqref{eq:zeta0}. So the product has at most a simple pole, and its residue only involves the two leading coefficients. Using $M(1)=1$, the residue of $\tilde\Psi(s)x^{-s}\zeta_K(s)$ at $s=0$ is
\[
-2^{r_1}\pi^{r_2}f(0)\cdot\Bigl(-\frac{h_KR_K}{w_K}\Bigr)=\frac{2^{r_1}\pi^{r_2}h_KR_K}{w_K}f(0).
\]

At $s=-j$ with $1\leqslant j<N$, the function $Z_f(s-1)$ has a pole of order at most $r_2$ if $j$ is odd and at most $r_1+r_2$ if $j$ is even, by Lemma~\ref{lem:zeta}. The zero of $\zeta_K$ at $s=-j$ has exactly this order. So there is no pole.

On the line $\Re s=-N+\frac12$, the decay of $M$ and the bound for $Z_f$ in Lemma~\ref{lem:zeta} (with $J=N+1$), together with the polynomial growth of $\zeta_K$, show that the integral is $O(x^{N-1/2}\|f\|_{C^k})$. Dividing by $\Delta_K$ and using~\eqref{eq:classnumber},
\[
\frac1{\Delta_K}\cdot\frac{2^{r_1}\pi^{r_2}h_KR_K}{w_K}=\frac{2^{r_1}\pi^{r_2}}{2^{r_1}(2\pi)^{r_2}}\sqrt{D_K}=\frac{\sqrt{D_K}}{2^{r_2}}.
\]
Replacing $N$ by $N+1$ gives the lemma.
\end{proof}

\subsection{The $\delta$-symbol over elements}\label{sec:elementdelta}
The weight $h(\N\fb/Q^d,\N\fa/Q^{2d})$ in Theorem~\ref{thm:BV} sees $\fa$ only through its norm. When $\fa=(a)$, it therefore forgets the sizes $|a^{(l)}|$ at the individual places, and keeps only their product. Section~\ref{sec:higher} explains why this loss is fatal for the first equation when $K$ has more than one archimedean place. We therefore also use a $\delta$-symbol whose divisors are elements, and whose archimedean weight is a product over the places.

Fix a non-negative $w_1\in C_c^\infty((\frac12,1))$ which is not identically zero, and put
\begin{equation}\label{eq:Wprod}
\bW(x)=\prod_{l=1}^{r_1+r_2}w_1\bigl(|x^{(l)}|\bigr)\qquad(x\in V).
\end{equation}
Then $\bW\in C_c^\infty(V)$, and $\bW(x)\neq0$ only if $\frac12<|x^{(l)}|<1$ at every place. In particular $\bW$ vanishes near every coordinate hyperplane, and $\bW(0)=0$. For $Q\geqslant1$ put
\[
C_Q=\sum_{c\in\fo\setminus\{0\}}\bW(c/Q).
\]
For a non-zero integral ideal $\fb$ and $a\in\fo$ define
\begin{equation}\label{eq:Hb}
\mathcal H_{\fb}(a)=\sum_{c\in\fb\setminus\{0\}}\frac1{|\N c|}\Bigl\{\bW\Bigl(\frac cQ\Bigr)-\bW\Bigl(\frac a{cQ}\Bigr)\Bigr\}.
\end{equation}
For fixed $a$ only finitely many $c$ contribute: the first term requires $|c^{(l)}|<Q$ for every $l$, and the second requires $|c^{(l)}|<2|a^{(l)}|/Q$ for every $l$ (and is absent if $a=0$).

\begin{theorem}[Element $\delta$-symbol]\label{thm:elementdelta}
Let $Q\geqslant1$. For every $a\in\fo$,
\[
\delta_K(a)=\frac1{C_Q}\sum_{(0)\neq\fb\subseteq\fo}\ \sumstar_{\gamma\bmod\fb}\ee(\gamma a)\,\mathcal H_{\fb}(a).
\]
Moreover
\[
C_Q=\frac{2^{r_2}}{\sqrt{D_K}}Q^d\int_V\bW(x)\,dx+O_N(Q^{-N}).
\]
\end{theorem}

\begin{proof}
First let $a\neq0$, and let $\mathcal D(a)=\{c\in\fo\setminus\{0\}:a/c\in\fo\}$ be the set of divisors of $a$ in $\fo$. It is a finite union of $U_K$-orbits, hence infinite when $U_K$ is infinite. However, only finitely many $c\in\mathcal D(a)$ have $\bW(c/Q)\neq0$, and only finitely many have $\bW(a/(cQ))\neq0$, because $\bW$ has compact support in $V^\times$. Since $c\mapsto a/c$ is a bijection of $\mathcal D(a)$,
\[
\sum_{c\in\mathcal D(a)}\Bigl\{\bW\Bigl(\frac cQ\Bigr)-\bW\Bigl(\frac a{cQ}\Bigr)\Bigr\}=\sum_{c\in\mathcal D(a)}\bW\Bigl(\frac cQ\Bigr)-\sum_{c'\in\mathcal D(a)}\bW\Bigl(\frac{c'}Q\Bigr)=0,
\]
both sums being finite. If $a=0$, then $\mathcal D(0)=\fo\setminus\{0\}$ and $\bW(0)=0$, so the same sum is $C_Q$. Therefore, for all $a\in\fo$,
\[
\delta_K(a)=\frac1{C_Q}\sum_{c\in\fo\setminus\{0\}}\mathbf 1[c\mid a]\Bigl\{\bW\Bigl(\frac cQ\Bigr)-\bW\Bigl(\frac a{cQ}\Bigr)\Bigr\}.
\]
By Lemma~\ref{lem:primchar} applied to the ideal $(c)$, and orthogonality on $\fo/(c)$,
\[
\mathbf 1[c\mid a]=\frac1{|\N c|}\sum_{\gamma\in(c)^{-1}\fd^{-1}/\fd^{-1}}\ee(\gamma a).
\]
By Lemma~\ref{lem:primchar} again, each class $\gamma\in(c)^{-1}\fd^{-1}/\fd^{-1}$ is primitive modulo exactly one divisor $\fb\mid(c)$, and the classes primitive modulo $\fb$ are exactly the primitive classes of $\fb^{-1}\fd^{-1}/\fd^{-1}$. Moreover $\fb\mid(c)$ if and only if $c\in\fb$. Inserting this and interchanging the finite sums over $c$ and $\fb$ gives the formula.

For the asymptotic formula, note that $C_Q=\sum_{c\in\fo}\bW(c/Q)$, since $\bW(0)=0$. Apply Lemma~\ref{lem:poisson-lattice} to the lattice $\fo\subset V$, whose covolume is $2^{-r_2}\sqrt{D_K}$ and whose dual is $\fd^{-1}$:
\[
C_Q=\frac{2^{r_2}}{\sqrt{D_K}}Q^d\sum_{\xi\in\fd^{-1}}\hat\bW(Q\xi).
\]
The term $\xi=0$ gives the main term. Every $\xi\in\fd^{-1}\setminus\{0\}$ satisfies $\ip\xi\gg1$, and $\hat\bW$ decays faster than any power. The sum over $\xi\neq0$ is therefore $O_N(Q^{-N})$.
\end{proof}

The ideal $\fb$ plays the role of the modulus, and the element $c\in\fb$ records the archimedean size of the divisor at each place. Unlike $h$, the weight $\mathcal H_{\fb}(a)$ depends on $a$ place by place, through $a/(cQ)$. Note that $c$ need not be balanced: when $a$ is small at some places, the relevant $c$ are skewed. The price is an extra summation over $c\in\fb$, which is harmless by the next lemma.

\begin{lemma}\label{lem:latticesum}
Let $\fb$ be a non-zero integral ideal and $X\geqslant1$. Then
\[
\sum_{\substack{c\in\fb\setminus\{0\}\\ \ip c\leqslant X}}\frac1{|\N c|}\ll\frac{(\log2X)^{r_1+r_2}}{\N\fb}.
\]
\end{lemma}

\begin{proof}
Decompose the sum according to the dyadic ranges of the sizes $|c^{(l)}|$, that is, according to the vector $(k_l)_l\in\ZZ^{r_1+r_2}$ with $2^{k_l}<|c^{(l)}|\leqslant2^{k_l+1}$. For given dyadic ranges, $|\N c|\asymp2^{\sum_lc_lk_l}$, and by Lemma~\ref{lem:boxcount} the number of $c\in\fb$ in these dyadic ranges is $\ll2^{\sum_lc_lk_l}/\N\fb$. Hence each choice of $(k_l)_l$ contributes $O(1/\N\fb)$. It remains to count the vectors $(k_l)_l$ that occur. We have $k_l\leqslant\log_2X$ for each $l$. Also $|\N c|\geqslant\N\fb\geqslant1$, so that $c_lk_l\geqslant-\sum_{l'\neq l}c_{l'}k_{l'}-O(1)\geqslant-d\log_2X-O(1)$. Thus each $k_l$ lies in an interval of length $O(\log2X)$, and there are $O((\log2X)^{r_1+r_2})$ such vectors.
\end{proof}

\begin{lemma}\label{lem:Hmain}
Let $Q\geqslant1$, let $\fb$ be a non-zero integral ideal with $\rho=\N\fb/Q^d\leqslant1$, and let $f\in C_c^\infty(V)$ be supported in $\{\ip v\leqslant R_0\}$. Then:
\begin{enumerate}[label=\textup{(\roman*)}]
\item for every $N\geqslant1$,
\[
\int_Vf(v)\,\mathcal H_{\fb}(Q^2v)\,dv=Q^{-d}\Bigl(f(0)\int_V\bW+O_{N,R_0}\bigl(\rho^N\|f\|_{C^k}\bigr)\Bigr)
\]
for some $k=k(N,d)$;
\item $\displaystyle\int_V\bigl|f(v)\,\mathcal H_{\fb}(Q^2v)\bigr|\,dv\ll(1+R_0)^d\frac{\|f\|_\infty}{\N\fb}$.
\end{enumerate}
\end{lemma}

\begin{proof}
Put $\Lambda=Q^{-1}\fb\subset V$, a lattice of covolume $2^{-r_2}\sqrt{D_K}\rho$. Writing $x=c/Q$, we have $|\N c|^{-1}=Q^{-d}|\N x|^{-1}$ and $Q^2v/(cQ)=v/x$, so that
\[
\mathcal H_{\fb}(Q^2v)=Q^{-d}\sum_{x\in\Lambda\setminus\{0\}}\frac1{|\N x|}\Bigl\{\bW(x)-\bW\Bigl(\frac vx\Bigr)\Bigr\}.
\]
For $x\in V^\times$ the substitution $v=xu$, for which $dv=|\N x|\,du$, gives
\[
\frac1{|\N x|}\int_Vf(v)\bW\Bigl(\frac vx\Bigr)\,dv=g(x),\qquad g(x)=\int_Vf(xu)\bW(u)\,du.
\]
The function $g$ is defined and smooth on all of $V$. Since $\bW(u)\neq0$ forces $|u^{(l)}|>\frac12$, $g$ is supported in $\{\ip x\leqslant2R_0\}$, and $\partial^\alpha g(x)=\int(\partial^\alpha f)(xu)u^\alpha\bW(u)\,du$, so that $\|g\|_{C^k}\ll\|f\|_{C^k}$. Only finitely many $x\in\Lambda$ meet the supports involved, so we may integrate term by term:
\[
\int_Vf(v)\mathcal H_{\fb}(Q^2v)\,dv=Q^{-d}\sum_{x\in\Lambda\setminus\{0\}}F(x),\qquad F(x)=\frac{\bW(x)}{|\N x|}\int_Vf-g(x).
\]
The function $|\N x|^{-1}\bW(x)$ is smooth with compact support, because $\bW$ vanishes near the coordinate hyperplanes. Hence $F\in C_c^\infty(V)$, with $F(0)=-g(0)=-f(0)\int\bW$. Moreover, the substitution $x=v/u$ shows that $\int_Vg(x)\,dx=\int_Vf\cdot\int_V\bW(u)|\N u|^{-1}\,du$, so that $\int_VF=0$.

By Lemma~\ref{lem:poisson-lattice},
\[
\sum_{x\in\Lambda}F(x)=\frac1{\operatorname{covol}(\Lambda)}\sum_{\xi\in\Lambda^\ast}\hat F(\xi)=\frac{2^{r_2}}{\sqrt{D_K}\rho}\sum_{\xi\in\Lambda^\ast\setminus\{0\}}\hat F(\xi),
\]
the term $\xi=0$ vanishing because $\int F=0$. Here $\Lambda^\ast=Q\hat\fb=Q\fb^{-1}\fd^{-1}$. A non-zero element $\xi=Q\eta$ of $\Lambda^\ast$ has $\ip\eta^d\geqslant|\N\eta|\geqslant\N(\fb^{-1}\fd^{-1})=(D_K\N\fb)^{-1}$, so that $\ip\xi\gg\rho^{-1/d}$. By Lemma~\ref{lem:boxcount}, the number of $\xi\in\Lambda^\ast$ with $\ip\xi\leqslant T$ is $\ll T^d\rho$. Integration by parts gives $\hat F(\xi)\ll_k(1+R_0)^d\|F\|_{C^k}\ip\xi^{-k}$. Taking $k=d(N+2)$ and summing over dyadic ranges of $\ip\xi$,
\[
\sum_{\xi\in\Lambda^\ast\setminus\{0\}}|\hat F(\xi)|\ll(1+R_0)^d\|f\|_{C^k}\sum_{T\gg\rho^{-1/d}}T^{d-k}\rho\ll(1+R_0)^d\|f\|_{C^k}\rho^{N+1}.
\]
Hence $\sum_{x\in\Lambda\setminus\{0\}}F(x)=-F(0)+O(\rho^N\|f\|_{C^k})=f(0)\int\bW+O(\rho^N\|f\|_{C^k})$, which proves~(i).

For~(ii), we bound the two parts of $\mathcal H_{\fb}$ separately. The terms with $\bW(x)$ require $\frac12<|x^{(l)}|<1$ for all $l$, so $|\N x|^{-1}\leqslant2^d$, and by Lemma~\ref{lem:boxcount} there are $\ll Q^d/\N\fb$ such $x\in\Lambda$. Their contribution is $\ll Q^{-d}\cdot(Q^d/\N\fb)\int|f|$. For the other terms,
\[
\frac1{|\N x|}\int_V|f(v)|\bW\Bigl(\frac vx\Bigr)\,dv=\int_V|f(xu)|\bW(u)\,du\leqslant\|f\|_\infty\int\bW,
\]
and this vanishes unless $\ip x\leqslant2R_0$. By Lemma~\ref{lem:boxcount} there are $\ll(R_0Q)^d/\N\fb$ such $x$. Multiplying by $Q^{-d}$ gives~(ii).
\end{proof}

With $Q=P$ and Theorem~\ref{thm:elementdelta}, Lemma~\ref{lem:Hmain}(i) reads
\[
\frac1{C_P}\int_Vf(v)\mathcal H_{\fb}(P^2v)\,dv=\frac{\sqrt{D_K}}{2^{r_2}}P^{-2d}\bigl(f(0)+O(\rho^N)\bigr).
\]
This has the same normalisation as Theorem~\ref{thm:BV} combined with Lemma~\ref{lem:h4}. The proof uses only Poisson summation on a lattice; the zeros of $\zeta_K$ play no role.

\section{The nested $\delta$-method over $K$}\label{sec:nested}

In this section we set up the nested $\delta$-method. We begin with a heuristic account of the sizes involved, then construct the nested $\delta$-symbol, and finally apply Poisson summation to $N(P)$.

\subsection{The counting function}\label{sec:counting}
Recall that $\bxi\in V^n\setminus X^\ast(V)$ is fixed, which means that $M_1\bxi^{(l)}$ and $M_2\bxi^{(l)}$ are linearly independent at every place $l$. Fix $\delta_0>0$ so small that the closed $\delta_0$-neighbourhood of $\bxi$ (with respect to $\ip\cdot$) still avoids $X^\ast(V)$, and fix $w_0\in C_c^\infty(V^n)$ with $0\leqslant w_0\leqslant1$, equal to $1$ on $\ip{\by-\bxi}\leqslant\delta_0/2$ and to $0$ outside $\ip{\by-\bxi}\leqslant\delta_0$. With $L=(\log P)^4$ put
\begin{equation}\label{eq:weight}
W(\by)=w_0(\by)\exp\bigl(-L\|\by-\bxi\|^2\bigr),\qquad N(P)=\sum_{\substack{\bx\in\fo^n\\ Q_1(\bx)=Q_2(\bx)=0}}W(\bx/P).
\end{equation}
On the support of $W$ the map $\by\mapsto(Q_1(\by),Q_2(\by))\in V^2$ is a submersion. Let $\cW(v_1,v_2)$ be the density of the push-forward of $W(\by)\,d\by$ under this map, so that $\int_{V^n}g(Q_1(\by),Q_2(\by))W(\by)\,d\by=\int_{V^2}g(v_1,v_2)\cW(v_1,v_2)\,dv_1\,dv_2$ for continuous $g$. By the coarea formula, $\cW(\bv)$ is the integral of $W$ against a smooth measure on the fibre over $\bv$. Hence $\cW$ is smooth with compact support, $0\leqslant\cW\ll1$, and its derivatives of order $j$ are $O(L^{j/2})$, because those of $W$ are. We put
\begin{equation}\label{eq:sigmainf}
\sigma_\infty=\cW(0,0).
\end{equation}
If $\bxi\in X(V)$, the fibre over $(0,0)$ is a smooth submanifold of real dimension $d(n-2)$ near $\bxi$, and $W$ restricted to it is a Gaussian of width $L^{-1/2}$. Therefore $\sigma_\infty\asymp L^{-d(n-2)/2}=(\log P)^{-2d(n-2)}$.

\subsection{Heuristics}\label{sec:heuristics}
Let $\bx\in\fo^n$ with $\bx/P$ in the support of $W$. Then $|Q_i(\bx)^{(l)}|\ll P^2$ at every place. Theorem~\ref{thm:BV} or Theorem~\ref{thm:elementdelta} with $Q=P$ detects $Q_1(\bx)=0$ with moduli $\fb_1$ of norm at most $P^d$. If we did the same for $Q_2(\bx)$, the total modulus $\fb_1\fb_2$ would have norm up to $P^{2d}$. Poisson summation in $\bx$ modulo $\fb_1\fb_2$ replaces the sum over $\bx$, which has about $P^{dn}$ terms, by a dual sum over $\bm\in(\fb_1\fb_2)^{\wedge n}$ of length about $(\N(\fb_1\fb_2)/P^d)^n$. For the largest moduli this is again about $P^{dn}$, and there is no saving.

The nested method reduces the size of the second equation using the first modulus. Once $\fb_1$ has been fixed, we write
\[
\delta_K(Q_2(\bx))=\mathbf 1\bigl[Q_2(\bx)\in\fb_1\bigr]\cdot\delta_K\bigl((Q_2(\bx))\fb_1^{-1}\bigr).
\]
The first factor is a congruence modulo $\fb_1$. The second factor concerns an integral ideal of norm about $P^{2d}/\N\fb_1$. If we insert a smooth cut-off forcing $|\Nm Q_2(\bx)|\ll P^{2d}/\N\fb_1$, the ideal $\delta$-symbol with parameter $Q$ satisfying $Q^{2d}=P^{2d}/\N\fb_1$ detects $(Q_2(\bx))\fb_1^{-1}=(0)$, with moduli $\fb_2$ of norm up to $Q^d=P^d/(\N\fb_1)^{1/2}$. In order to apply Theorem~\ref{thm:BV} with a parameter independent of $\fb_1$, we instead take $Q=P^{1/2}$ and a cut-off of the shape $|\Nm Q_2(\bx)|\ll P^{d}\N\fb_1$, which gives $\N\fb_2\ll P^{d/2}$. In either form, the total modulus satisfies
\[
\N(\fb_1\fb_2)\ll P^{3d/2}.
\]
After Poisson summation modulo $\fb_1\fb_2$, the dual sum has length about $(P^{3d/2}/P^d)^n=P^{dn/2}$. In return the cut-off is only correct when it does not destroy the solutions; it does not, since $Q_2(\bx)=0$ lies in the region where the cut-off equals $1$.

The congruence $Q_2(\bx)\in\fb_1$ is periodic modulo $\fb_1$, so it is absorbed into the complete exponential sum modulo $\fb_1\fb_2$. The resulting sums $C_{\fb_1,\fb_2}(\bm)$ are the analogues of Munshi's sums~\cite[(7)]{Mun15}. For generic $\bm$, the sum over $\fb_1$ exhibits square-root cancellation, which is the Kloosterman refinement in the first modulus. The sum over $\fb_2$ factors through a quadratic character attached to $Q_2^\ast(\bm)$; averaging over $\fb_2$ leads to Hecke $L$-functions, and this gives a second saving. These savings are worked out in Sections~\ref{sec:sums} and~\ref{sec:error}.

\textbf{Archimedean sizes.} Over $K$ the size of $Q_2(\bx)$ is a vector of sizes. The cut-off for the second equation must be compatible, at each place, with the weight from the first equation, and this is where we use the element $\delta$-symbol. The terms of Theorem~\ref{thm:elementdelta} applied to $a_1=Q_1(\bx)$ carry an element $c_1\in\fb_1$, and the weight $\bW(Q_1(\bx)/(c_1P))$ varies, at place $l$, on the scale $|c_1^{(l)}|P$ in $Q_1(\bx)^{(l)}$. We choose a place-wise cut-off for $Q_2(\bx)$, with the same archimedean sizes; see~\eqref{eq:lambdavec} below. Section~\ref{sec:higher} explains why a scalar cut-off fails.

\subsection{The place-wise cut-off}
Fix $U\in C^\infty([0,\infty))$ with $0\leqslant U\leqslant1$, $U=1$ on $[0,1]$ and $U=0$ on $[2,\infty)$, and put
\[
U_V(v)=\prod_lU\bigl(|v^{(l)}|\bigr)\qquad(v\in V).
\]
Then $U_V\in C_c^\infty(V)$ and $U_V$ is identically $1$ near $v=0$. For a non-zero integral ideal $\fb_1$ and $c_1\in\fb_1\setminus\{0\}$, define $\bl=\bl(\fb_1,c_1)\in V$ by
\begin{equation}\label{eq:lambdavec}
\bl^{(l)}=\frac{P}{\kappa\,|c_1^{(l)}|},\qquad\kappa=\kappa(\fb_1,c_1)=\Bigl(\frac{\N\fb_1}{|\N c_1|}\Bigr)^{1/d}.
\end{equation}

\begin{lemma}\label{lem:lambda}
Let $\fb_1$ and $c_1$ be as above.
\begin{enumerate}[label=\textup{(\roman*)}]
\item $\Nm(\bl)=P^d/\N\fb_1$; in particular, $\Nm\bl$ depends only on $\fb_1$.
\item $\kappa\leqslant1$, and $\bl^{(l)}\geqslant P/|c_1^{(l)}|$ for every $l$.
\item If $\ip{c_1}\leqslant AP$, then $\bl^{(l)}\geqslant A^{-1}$ for every $l$.
\item If $\ip{c_1}\leqslant AP$, then $\bl^{(l)}\ll_AP^{d+2}$ for every $l$.
\end{enumerate}
\end{lemma}

\begin{proof}
(i) We have $\prod_l(\bl^{(l)})^{c_l}=P^d\kappa^{-d}|\N c_1|^{-1}=P^d/\N\fb_1$. (ii) Since $c_1\in\fb_1\setminus\{0\}$, $|\N c_1|\geqslant\N\fb_1$, so $\kappa\leqslant1$. (iii) follows from (ii). (iv) Since $|\N c_1|\geqslant1$ and $|c_1^{(l')}|\leqslant AP$ for all $l'$, we have $|c_1^{(l)}|^{c_l}\geqslant(AP)^{-(d-c_l)}$, hence $|c_1^{(l)}|\geqslant(AP)^{-d}$. Also $\kappa^d=\N\fb_1/|\N c_1|\geqslant(AP)^{-d}$. Hence $\bl^{(l)}\leqslant P(AP)(AP)^d=A^{d+1}P^{d+2}$.
\end{proof}

\begin{lemma}[Nesting identity]\label{lem:nest-identity}
Let $\fb_1$ be a non-zero integral ideal, $\mu\in V$ and $a_2\in\fo$. Then
\[
\delta_K(a_2)=\mathbf 1[a_2\in\fb_1]\cdot\delta_K\bigl((a_2)\fb_1^{-1}\bigr)\cdot U_V(\mu a_2).
\]
\end{lemma}

\begin{proof}
If $a_2=0$, all three factors on the right equal $1$. If $a_2\neq0$, the left side is $0$; on the right either $a_2\notin\fb_1$, or $(a_2)\fb_1^{-1}$ is a non-zero integral ideal and $\delta_K$ of it vanishes.
\end{proof}

The identity holds for every $\mu$. In particular $\mu$ may depend on $c_1$: we apply it separately to each term of the element $\delta$-symbol for $a_1$.

\subsection{The nested $\delta$-symbol}
\begin{theorem}[Nested $\delta$-symbol]\label{thm:nested}
Let $P\geqslant1$ and $A\geqslant1$, and let $a_1,a_2\in\fo$ with $\ip{a_i}\leqslant AP^2$. Put $c_{\ast}=2^{d+1}$. Then
\begin{multline*}
\delta_K(a_1)\delta_K(a_2)=\frac{c_{P^{1/2}}}{C_PP^{d}}\sum_{\fb_1}\ \sumstar_{\gamma_1\bmod\fb_1}\ \sum_{c_1\in\fb_1\setminus\{0\}}\frac{\mathbf 1[a_2\in\fb_1]\,\ee(\gamma_1a_1)}{|\N c_1|}\Bigl\{\bW\Bigl(\frac{c_1}P\Bigr)-\bW\Bigl(\frac{a_1}{c_1P}\Bigr)\Bigr\}\\
\times\sum_{\N\fb_2\leqslant c_{\ast}P^{d/2}}c_{\fb_2}\bigl((a_2)\fb_1^{-1}\bigr)\,h\Bigl(\frac{\N\fb_2}{P^{d/2}},\frac{|\Nm a_2|}{P^d\N\fb_1}\Bigr)\,U_V\Bigl(\frac{\bl(\fb_1,c_1)a_2}{P^2}\Bigr),
\end{multline*}
where $c_{\fb_2}$ is the Ramanujan sum of Lemma~\ref{lem:orth}. Only terms with $\ip{c_1}\leqslant2AP$ and $\N\fb_1\leqslant(2AP)^d$ occur, and all sums are finite.
\end{theorem}

\begin{proof}
Apply Theorem~\ref{thm:elementdelta} with $Q=P$ to $a_1$. The term indexed by $(\fb_1,\gamma_1,c_1)$ is non-zero only if $\bW(c_1/P)\neq0$, which forces $|c_1^{(l)}|<P$ for all $l$, or $\bW(a_1/(c_1P))\neq0$, which forces $|c_1^{(l)}|<2|a_1^{(l)}|/P\leqslant2AP$ for all $l$. In either case $\ip{c_1}\leqslant2AP$, and then $\N\fb_1\leqslant|\N c_1|\leqslant(2AP)^d$. Multiply each term by $\delta_K(a_2)$, written by Lemma~\ref{lem:nest-identity} with $\mu=\bl(\fb_1,c_1)/P^2$. In the terms with $a_2\in\fb_1$, apply Theorem~\ref{thm:BV} with $Q=P^{1/2}$ to the integral ideal $\fa_2=(a_2)\fb_1^{-1}$. Here $Q^d=P^{d/2}$, $Q^{2d}=P^d$ and $\N\fa_2=|\Nm a_2|/\N\fb_1$, which gives the stated arguments of $h$. On the support of $U_V(\bl a_2/P^2)$ we have $|\bl^{(l)}a_2^{(l)}|\leqslant2P^2$ for every $l$, and hence, by Lemma~\ref{lem:lambda}(i),
\[
\frac{\N\fa_2}{P^d}=\frac{|\Nm a_2|}{P^d\N\fb_1}=\frac{|\Nm(\bl a_2)|}{P^{2d}}\leqslant2^d.
\]
Since $h(x,y)=0$ unless $x\leqslant\max(1,2|y|)$, only $\fb_2$ with $\N\fb_2/P^{d/2}\leqslant2^{d+1}$ occur.
\end{proof}

\begin{remark}\label{rem:cutoff}
We record why the archimedean sizes of $\bl$ matter; the precise statement is Theorem~\ref{thm:truncation}. For $a_1=Q_1(\bx)$ with $\bx=P\by$, the first $\delta$-symbol produces the weight $\bW(PQ_1(\by)/c_1)$. At place $l$ it varies on the scale $|c_1^{(l)}|/P$ in $Q_1(\by)^{(l)}$, so its Fourier transform lives on frequencies $|t_1^{(l)}|\lesssim P/|c_1^{(l)}|\leqslant\bl^{(l)}$. The cut-off $U_V(\bl Q_2(\by))$, together with the weight $h$ for the second equation, produces frequencies of size $\bl^{(l)}\max(1,|t_2^{(l)}|)$ at place $l$, where $\cH(t_2)\lesssim P^{d/2}/\N\fb_2$. So the frequency coming from the second equation dominates the one from the first equation at \emph{every} place, and the product over the places is controlled by norms alone. With a weight $h(\rho_1,|\Nm Q_1(\by)|)$ in place of $\bW$, the frequency $t_1$ would only be controlled through $\cH(t_1)$, and this domination can fail.
\end{remark}

\subsection{Additive characters for the second modulus}
To apply Poisson summation we express the Ramanujan sum $c_{\fb_2}((a_2)\fb_1^{-1})$ through additive characters in $a_2$.

\begin{lemma}\label{lem:additive}
Let $\fb_1,\fb_2$ be non-zero integral ideals and $a_2\in\fb_1$. Then
\[
c_{\fb_2}\bigl((a_2)\fb_1^{-1}\bigr)=\sumstar_{\gamma_2}\ee(\gamma_2a_2),
\]
where $\gamma_2$ runs over the classes of $(\fb_1\fb_2)^{-1}\fd^{-1}/\fb_1^{-1}\fd^{-1}$ that do not lie in $(\fb_1\fe)^{-1}\fd^{-1}$ for any divisor $\fe\mid\fb_2$ with $\fe\neq\fb_2$.
\end{lemma}

\begin{proof}
Let $\fe\mid\fb_2$. For $\gamma\in(\fb_1\fe)^{-1}\fd^{-1}$, the value $\ee(\gamma a_2)$ depends only on the class of $\gamma$ modulo $\fb_1^{-1}\fd^{-1}$, because $a_2\in\fb_1$. The group $(\fb_1\fe)^{-1}\fd^{-1}/\fb_1^{-1}\fd^{-1}$ has order $\N\fe$. The character $\gamma\mapsto\ee(\gamma a_2)$ on it is trivial if and only if $a_2(\fb_1\fe)^{-1}\fd^{-1}\subseteq\fd^{-1}$, that is, $a_2\in\fb_1\fe$, that is, $\fe\mid(a_2)\fb_1^{-1}$. Hence
\[
\sum_{\gamma\in(\fb_1\fe)^{-1}\fd^{-1}/\fb_1^{-1}\fd^{-1}}\ee(\gamma a_2)=\N\fe\cdot\mathbf 1\bigl[\fe\mid(a_2)\fb_1^{-1}\bigr].
\]
The group for $\fe$ is the union of the sets of classes that are ``primitive at level $\fe'$'' for $\fe'\mid\fe$. M\"obius inversion over the divisors of $\fb_2$ therefore gives
\[
\sumstar_{\gamma_2}\ee(\gamma_2a_2)=\sum_{\fe\mid\fb_2}\mu(\fb_2\fe^{-1})\N\fe\cdot\mathbf 1\bigl[\fe\mid(a_2)\fb_1^{-1}\bigr],
\]
which is $c_{\fb_2}((a_2)\fb_1^{-1})$ by Lemma~\ref{lem:orth}(iii). For $a_2=0$ both sides equal the number of primitive classes, as $\sigma((0))=1$.
\end{proof}

\begin{lemma}\label{lem:periodic}
Let $\gamma_1$ be a class of $\fb_1^{-1}\fd^{-1}/\fd^{-1}$ and $\gamma_2$ a class of $(\fb_1\fb_2)^{-1}\fd^{-1}/\fb_1^{-1}\fd^{-1}$. Then, for $\ba\in\fo^n$ with $Q_2(\ba)\in\fb_1$, the value $\ee(\gamma_1Q_1(\ba)+\gamma_2Q_2(\ba))$ does not depend on the representatives $\gamma_1,\gamma_2$, and it depends only on $\ba$ modulo $\fb_1\fb_2$. The condition $Q_2(\ba)\in\fb_1$ depends only on $\ba$ modulo $\fb_1$.
\end{lemma}

\begin{proof}
Changing $\gamma_1$ by an element of $\fd^{-1}$ changes $\gamma_1Q_1(\ba)$ by an element of $\fd^{-1}$, and changing $\gamma_2$ by an element of $\fb_1^{-1}\fd^{-1}$ changes $\gamma_2Q_2(\ba)$ by an element of $\fd^{-1}$, because $Q_2(\ba)\in\fb_1$. For $\by\in(\fb_1\fb_2)^n$ we have
\[
Q_i(\ba+\by)-Q_i(\ba)=\by^t(2M_i)\ba+Q_i(\by)\in\fb_1\fb_2,
\]
since $2M_i$ has entries in $\fo$. As $\gamma_1$ and $\gamma_2$ lie in $(\fb_1\fb_2)^{-1}\fd^{-1}$, the phase changes by an element of $\fd^{-1}$. The last assertion is proved in the same way.
\end{proof}

\subsection{Poisson summation}
For non-zero integral ideals $\fb_1,\fb_2$ and $\bm\in(\fb_1\fb_2)^{\wedge n}$ put
\begin{equation}\label{eq:Csum}
C_{\fb_1,\fb_2}(\bm)=\sumstar_{\gamma_1}\ \sumstar_{\gamma_2}\sum_{\substack{\ba\in(\fo/\fb_1\fb_2)^n\\ Q_2(\ba)\in\fb_1}}\ee\bigl(\gamma_1Q_1(\ba)+\gamma_2Q_2(\ba)+\bm\cdot\ba\bigr),
\end{equation}
where $\gamma_1$ runs over the primitive classes of $\fb_1^{-1}\fd^{-1}/\fd^{-1}$ and $\gamma_2$ as in Lemma~\ref{lem:additive}. By Lemma~\ref{lem:periodic} the summand is well defined; note that $\ee(\bm\cdot\ba)$ is well defined modulo $\fb_1\fb_2$ because $\bm\in(\fb_1^{-1}\fb_2^{-1}\fd^{-1})^n$. When $K=\QQ$ these are the sums of~\cite[(7)]{Mun15}. Next, for $c_1\in\fb_1\setminus\{0\}$, put $\rho_2=\N\fb_2/P^{d/2}$ and $\bl=\bl(\fb_1,c_1)$, and define
\begin{multline}\label{eq:Iint}
I_{c_1,\fb_2}(\bm)=\int_{V^n}W(\by)\Bigl\{\bW\Bigl(\frac{c_1}P\Bigr)-\bW\Bigl(\frac{PQ_1(\by)}{c_1}\Bigr)\Bigr\}\,h\bigl(\rho_2,\Nm(\bl Q_2(\by))\bigr)\\
\times U_V\bigl(\bl Q_2(\by)\bigr)\,\ee(-P\bm\cdot\by)\,d\by.
\end{multline}
Finally put
\begin{equation}\label{eq:Ibig}
I_{\fb_1,\fb_2}(\bm)=\frac{P^{2d}}{C_P}\sum_{\substack{c_1\in\fb_1\setminus\{0\}\\ \ip{c_1}\leqslant A_0P}}\frac{I_{c_1,\fb_2}(\bm)}{|\N c_1|},
\end{equation}
where $A_0$ is a constant such that $\ip{Q_i(\by)}\leqslant A_0/2$ on the support of $W$.

\begin{proposition}\label{prop:poisson}
For every $A>0$,
\begin{multline*}
N(P)=\Bigl(\frac{2^{r_2}}{\sqrt{D_K}}\Bigr)^nP^{d(n-3)}\sum_{\N\fb_1\leqslant(A_0P)^d}\ \sum_{\N\fb_2\leqslant c_{\ast}P^{d/2}}\frac1{\N(\fb_1\fb_2)^n}\\
\times\sum_{\bm\in(\fb_1\fb_2)^{\wedge n}}C_{\fb_1,\fb_2}(\bm)\,I_{\fb_1,\fb_2}(\bm)+O_A(P^{-A}).
\end{multline*}
The series over $\bm$ converges absolutely.
\end{proposition}

\begin{proof}
For $\bx\in\fo^n$ with $W(\bx/P)\neq0$ we have $\ip{Q_i(\bx)}\leqslant(A_0/2)P^2$. Apply Theorem~\ref{thm:nested} with $A=A_0/2$ to $a_i=Q_i(\bx)$, multiply by $W(\bx/P)$ and sum over $\bx$. For each $\bx$ all the sums in Theorem~\ref{thm:nested} are finite, and only finitely many $\bx$ occur, so we may interchange the order of summation freely. Replace $c_{\fb_2}((Q_2(\bx))\fb_1^{-1})$ by the sum over $\gamma_2$ in Lemma~\ref{lem:additive}. We obtain
\[
N(P)=\frac{c_{P^{1/2}}}{C_PP^d}\sum_{\fb_1,c_1,\fb_2}\frac1{|\N c_1|}\sumstar_{\gamma_1}\sumstar_{\gamma_2}\sum_{\substack{\bx\in\fo^n\\ Q_2(\bx)\in\fb_1}}\ee\bigl(\gamma_1Q_1(\bx)+\gamma_2Q_2(\bx)\bigr)\,\Psi(\bx),
\]
where $\Psi(\bx)=\Psi_{c_1,\fb_2}(\bx)$ is the product of $W(\bx/P)$ with the archimedean weights of Theorem~\ref{thm:nested}. Split $\bx$ into residue classes $\ba$ modulo $\fb_1\fb_2$. By Lemma~\ref{lem:periodic}, the condition and the phase depend only on $\ba$. Lemma~\ref{lem:poisson-lattice} gives
\[
\sum_{\bx\equiv\ba\bmod\fb_1\fb_2}\Psi(\bx)=\frac{2^{r_2n}}{D_K^{n/2}\N(\fb_1\fb_2)^n}\sum_{\bm\in(\fb_1\fb_2)^{\wedge n}}\ee(\bm\cdot\ba)\hat\Psi(\bm).
\]
Summing over $\ba$, $\gamma_1$ and $\gamma_2$ produces $C_{\fb_1,\fb_2}(\bm)$. Substituting $\bx=P\by$, we have $d\bx=P^{dn}d\by$, and, with $\bl=\bl(\fb_1,c_1)$,
\begin{gather*}
\frac{Q_1(\bx)}{c_1P}=\frac{PQ_1(\by)}{c_1},\qquad\frac{\bl Q_2(\bx)}{P^2}=\bl Q_2(\by),\\
\frac{|\Nm Q_2(\bx)|}{P^d\N\fb_1}=\frac{P^d|\Nm Q_2(\by)|}{\N\fb_1}=|\Nm(\bl Q_2(\by))|,
\end{gather*}
where we used Lemma~\ref{lem:lambda}(i). Therefore $\hat\Psi(\bm)=P^{dn}I_{c_1,\fb_2}(\bm)$. The terms with $\ip{c_1}>A_0P$ vanish by Theorem~\ref{thm:nested}. Collecting the factors,
\[
\frac{c_{P^{1/2}}}{C_PP^d}\cdot\frac{2^{r_2n}P^{dn}}{D_K^{n/2}}=\Bigl(\frac{2^{r_2}}{\sqrt{D_K}}\Bigr)^nP^{d(n-3)}\cdot\frac{P^{2d}}{C_P}\cdot c_{P^{1/2}}.
\]
Since $c_{P^{1/2}}=1+O(P^{-A})$ and trivially $0\leqslant N(P)\ll P^{dn}$, replacing $c_{P^{1/2}}$ by $1$ costs $O(P^{-A})$. The absolute convergence of the sum over $\bm$ follows from Lemma~\ref{lem:ibp} below.
\end{proof}

By Lemma~\ref{lem:latticesum} and $C_P\asymp P^d$, the trivial bound $I_{c_1,\fb_2}(\bm)\ll\rho_2^{-1}$, which follows from $h\ll\rho_2^{-1}$, gives
\begin{equation}\label{eq:Itrivial}
I_{\fb_1,\fb_2}(\bm)\ll(\log P)^{r_1+r_2}\frac{P^d}{\N\fb_1}\cdot\frac{P^{d/2}}{\N\fb_2}=(\log P)^{r_1+r_2}\frac{P^{3d/2}}{\N(\fb_1\fb_2)}.
\end{equation}
This is the analogue of the trivial bound $B^{3/2}/(q_1q_2)$ in~\cite{Mun15}.

\subsection{The shape of the argument}
We write
\[
N(P)=M(P)+E(P),
\]
where $M(P)$ is the contribution of $\bm=\mathbf0$ in Proposition~\ref{prop:poisson} and $E(P)$ is the contribution of $\bm\neq\mathbf0$. For the main term $M(P)$ we need an asymptotic formula for $I_{\fb_1,\fb_2}(\mathbf0)$; this is where Lemmas~\ref{lem:h4} and~\ref{lem:Hmain} enter (Section~\ref{sec:main}). For the error term $E(P)$ we need three things:
\begin{enumerate}[label=(\alph*)]
\item an archimedean truncation, showing that $I_{\fb_1,\fb_2}(\bm)$ is negligible unless $\cH(P\bm)\leqslant P^{3d/2+\varepsilon}/\N(\fb_1\fb_2)$ (Theorem~\ref{thm:truncation});
\item a bound of the shape $I_{\fb_1,\fb_2}(\bm)\ll P^\varepsilon P^{3d/2}\N(\fb_1\fb_2)^{-1}\cH(P\bm)^{-n/2}$ (Corollary~\ref{cor:Ibound}), and some control on its dependence on $\fb_2$ (Lemma~\ref{lem:smooth});
\item bounds for $C_{\fb_1,\fb_2}(\bm)$ and its averages over $\fb_1$ and $\fb_2$ (Section~\ref{sec:sums}).
\end{enumerate}
The element $\delta$-symbol and the vector $\bl$ are designed to make (a) work for every number field.

\begin{center}
\renewcommand{\arraystretch}{1.25}
\begin{tabular}{c|c}
Munshi, over $\QQ$ & this paper, over $K$\\\hline
$q_1\leqslant B$, $q_2\leqslant B^{1/2}$ & $\N\fb_1\ll P^d$, $\N\fb_2\ll P^{d/2}$\\
weight $h(q_1/B,Q_1/B^2)$ & $\sum_{c_1\in\fb_1}|\N c_1|^{-1}\{\bW(c_1/P)-\bW(Q_1/(c_1P))\}$\\
$Q_2(\bx)/q_1$ & the ideal $(Q_2(\bx))\fb_1^{-1}$\\
cut-off $U(Q_2/(q_1B))$ & $U_V(\bl(\fb_1,c_1)Q_2/P^2)$, matched to $c_1$\\
truncation $|\bm|\ll B^{1/2+\varepsilon}$ & $\cH(\tv)\ll P^{d/2+\varepsilon}$ (Section~\ref{sec:error})\\
Dirichlet $L$-functions & Hecke $L$-functions with infinity type
\end{tabular}
\end{center}

\section{The exponential sums}\label{sec:sums}

In this section we study the sums $C_{\fb_1,\fb_2}(\bm)$ defined in~\eqref{eq:Csum}. They factor into local sums (Lemma~\ref{lem:mult}), and the local sums are the analogues of the sums studied by Browning and Munshi~\cite{BM13} and Munshi~\cite{Mun15}, with $\ZZ/p^k\ZZ$ replaced by $\fo/\fp^k$. The main results are a uniform bound valid at every prime (Lemma~\ref{lem:universal}), the Kloosterman refinement in the first modulus (Lemmas~\ref{lem:first-generic} and~\ref{lem:first-degenerate}), the vanishing of the mixed sums (Lemma~\ref{lem:mixed}), and the evaluation of the sums for the second modulus in terms of a quadratic character (Lemma~\ref{lem:gauss}). The average of this character over the second modulus is treated in~\S\ref{sec:hecke}.

\subsection{Local sums}
Let $\fp$ be a prime ideal. We write $K_{\fp}$ and $\cO_{\fp}$ for the completions, $\varpi$ for a uniformiser, $q=\N\fp$, and $\delta_{\fp}$ for a generator of the local different $\fd\cO_{\fp}$. Let $p$ be the rational prime below $\fp$, let $\{\cdot\}_p$ denote the $p$-adic fractional part, and put
\[
\psi_{\fp}(x)=\exp\bigl(2\pi i\{\Tr_{K_{\fp}/\QQ_p}x\}_p\bigr)\qquad(x\in K_{\fp}).
\]
This character is trivial on $\delta_{\fp}^{-1}\cO_{\fp}$ but not on $\varpi^{-1}\delta_{\fp}^{-1}\cO_{\fp}$. For $\gamma\in K$ we have $\Tr_{K/\QQ}\gamma=\sum_{\fp\mid p}\Tr_{K_{\fp}/\QQ_p}\gamma$, and for $r\in\QQ$ the difference $r-\sum_p\{r\}_p$ is an integer. Hence
\[
\ee(\gamma)=\prod_{\fp}\psi_{\fp}(\gamma)\qquad(\gamma\in K),
\]
where only finitely many factors differ from $1$. For $j\geqslant0$ and $x\in\cO_{\fp}$ put
\begin{equation}\label{eq:ep}
\ep j(x)=\psi_{\fp}\bigl(x/(\varpi^j\delta_{\fp})\bigr).
\end{equation}
This is a primitive additive character of $\cO_{\fp}/\fp^j$, and $\ep j(\varpi x)=\ep{j-1}(x)$. For $\bm\in(\fp^{-j}\fd^{-1}\cO_{\fp})^n$ we write $\bm=\tm/(\varpi^j\delta_{\fp})$ with $\tm=\tm^{(j)}\in\cO_{\fp}^n$, and call $\tm$ the \emph{normalised frequency} at level $j$. It depends on the choice of $\varpi$ and $\delta_{\fp}$ only up to multiplication by a unit, and Lemma~\ref{lem:unitinv} below shows that this does not matter.

For $r,\ell\geqslant0$ and $\bm\in(\fp^{-r-\ell}\fd^{-1}\cO_{\fp})^n$ let $C_{\fp}(r,\ell;\bm)$ be defined by~\eqref{eq:Csum}, with $\fo$, $\fb_1$, $\fb_2$ and $\ee$ replaced by $\cO_{\fp}$, $\fp^r$, $\fp^\ell$ and $\psi_{\fp}$. Following~\cite{BM13,Mun15} we also put
\begin{align*}
S_{\fp}(r,\ell;\bm)&=\sumstar_{a\bmod\fp^\ell}\ \sum_{\substack{\bb\bmod\fp^{r+\ell}\\ Q_1(\bb)\equiv Q_2(\bb)\equiv0\ (\fp^r)}}\ep{r+\ell}\bigl(aQ_2(\bb)+\tm\cdot\bb\bigr),\qquad\mathcal Q_{\fp}(\ell;\bm)=S_{\fp}(0,\ell;\bm),
\end{align*}
where $\tm=\tm^{(r+\ell)}$ and $\sumstar$ denotes a sum over units. These are the analogues of the sums $\mathsf Q_q$ and $\mathsf M_{d,q}$ in~\cite[\S\S4--6]{BM13}.

\begin{lemma}[Dictionary]\label{lem:dictionary}
With $\tm=\tm^{(r+\ell)}$,
\[
C_{\fp}(r,\ell;\bm)=\sumstar_{a_1\bmod\fp^r}\ \sumstar_{a_2\bmod\fp^\ell}\ \sum_{\substack{\bb\bmod\fp^{r+\ell}\\ Q_2(\bb)\equiv0\ (\fp^r)}}\ep{r+\ell}\bigl(a_1\varpi^\ell Q_1(\bb)+a_2Q_2(\bb)+\tm\cdot\bb\bigr).
\]
In particular $C_{\fp}(0,\ell;\bm)=\mathcal Q_{\fp}(\ell;\bm)$ and $C_{\fp}(0,0;\bm)=1$.
\end{lemma}

\begin{proof}
The primitive classes of the group $\fp^{-r}\delta_{\fp}^{-1}\cO_{\fp}/\delta_{\fp}^{-1}\cO_{\fp}$ are the classes of $\gamma_1=a_1/(\varpi^r\delta_{\fp})=a_1\varpi^\ell/(\varpi^{r+\ell}\delta_{\fp})$ with $a_1\in(\cO_{\fp}/\fp^r)^\times$. For $\ell\geqslant1$, the admissible classes of $\fp^{-r-\ell}\delta_{\fp}^{-1}\cO_{\fp}/\fp^{-r}\delta_{\fp}^{-1}\cO_{\fp}$ are those of $\gamma_2=a_2/(\varpi^{r+\ell}\delta_{\fp})$ with $a_2\in(\cO_{\fp}/\fp^\ell)^\times$; for $\ell=0$ there is the single class $\gamma_2=0$, which corresponds to $a_2$ running over the one-element set $(\cO_{\fp}/\cO_{\fp})^\times$. Finally $\bm\cdot\bb=\tm\cdot\bb/(\varpi^{r+\ell}\delta_{\fp})$.
\end{proof}

Thus $C_{\fp}(r,\ell;\bm)$ is~\cite[(7)]{Mun15} with $q_1\mapsto\varpi^r$, $q_2\mapsto\varpi^\ell$ and $e_{q_1q_2}\mapsto\ep{r+\ell}$.

\begin{lemma}[Multiplicativity]\label{lem:mult}
Let $\fb_1=\prod_{\fp}\fp^{r_{\fp}}$ and $\fb_2=\prod_{\fp}\fp^{\ell_{\fp}}$. Then, for $\bm\in(\fb_1\fb_2)^{\wedge n}$,
\[
C_{\fb_1,\fb_2}(\bm)=\prod_{\fp\mid\fb_1\fb_2}C_{\fp}(r_{\fp},\ell_{\fp};\bm).
\]
\end{lemma}

\begin{proof}
By the Chinese remainder theorem, $\fo/\fb_1\fb_2\cong\prod_{\fp}\cO_{\fp}/\fp^{r_{\fp}+\ell_{\fp}}$ and
\[
\begin{gathered}
\fb_1^{-1}\fd^{-1}/\fd^{-1}\cong\bigoplus_{\fp}\fp^{-r_{\fp}}\delta_{\fp}^{-1}\cO_{\fp}/\delta_{\fp}^{-1}\cO_{\fp},\\
(\fb_1\fb_2)^{-1}\fd^{-1}/\fb_1^{-1}\fd^{-1}\cong\bigoplus_{\fp}\fp^{-r_{\fp}-\ell_{\fp}}\delta_{\fp}^{-1}\cO_{\fp}/\fp^{-r_{\fp}}\delta_{\fp}^{-1}\cO_{\fp}.
\end{gathered}
\]
A class is primitive (in the sense of~\eqref{eq:Csum}) if and only if each of its local components is. The condition $Q_2(\ba)\in\fb_1$ holds if and only if $Q_2(\ba)\in\fp^{r_{\fp}}\cO_{\fp}$ for every $\fp$. Finally, the argument $z=\gamma_1Q_1(\ba)+\gamma_2Q_2(\ba)+\bm\cdot\ba$ lies in $K$, and $\ee(z)=\prod_{\fp}\psi_{\fp}(z)$. For $\fp\nmid\fb_1\fb_2$ we have $z\in\delta_{\fp}^{-1}\cO_{\fp}$ and $\psi_{\fp}(z)=1$, while for $\fp\mid\fb_1\fb_2$ the value $\psi_{\fp}(z)$ depends only on the $\fp$-components. So the sum factorises.
\end{proof}

\begin{lemma}\label{lem:unitinv}
For $u\in\cO_{\fp}^\times$ we have $C_{\fp}(r,\ell;u\bm)=C_{\fp}(r,\ell;\bm)$, and the same holds for $S_{\fp}$.
\end{lemma}

\begin{proof}
In Lemma~\ref{lem:dictionary}, substitute $\bb=u^{-1}\bb'$. Then $Q_i(\bb)=u^{-2}Q_i(\bb')$, the congruence condition is unchanged, $u\tm\cdot\bb=\tm\cdot\bb'$, and $a_i\mapsto u^{-2}a_i$ permutes the units.
\end{proof}

Over $\QQ$, the multiplicativity relation in~\cite[\S3]{Mun15} involves twisting $\bm$ by inverses of the coprime parts of the moduli. Over $K$ the local definition avoids this, and Lemma~\ref{lem:unitinv} shows that such twists would be harmless anyway.

\subsection{Local tools}
For $A\in M_n(\cO_{\fp})$ and $\ba\in\cO_{\fp}^n$ put $K_{\fp^k}(A;\ba)=\#\{\bx\bmod\fp^k:A\bx\equiv\ba\bmod\fp^k\}$.

\begin{lemma}\label{lem:smith}
Let $A\in M_n(\cO_{\fp})$ with $\det A\neq0$, and let $\varpi^{e_1},\dots,\varpi^{e_n}$ be its elementary divisors, so that $e_1+\dots+e_n=\ord_{\fp}\det A$. Then $K_{\fp^k}(A;\ba)\leqslant K_{\fp^k}(A;\mathbf0)=q^{\sum_i\min(e_i,k)}$. Moreover, $e_1+\dots+e_{n-1}$ is the minimum of the orders of the $(n-1)\times(n-1)$ minors of $A$.
\end{lemma}

\begin{proof}
This is the analogue of~\cite[Lemma~1]{BM13}. Since $\cO_{\fp}$ is a principal ideal domain, there are $U,U'\in GL_n(\cO_{\fp})$ with $UAU'=\operatorname{diag}(\varpi^{e_1},\dots,\varpi^{e_n})$. The system $A\bx\equiv\ba$ becomes $\varpi^{e_i}x_i'\equiv a_i'\bmod\fp^k$, $1\leqslant i\leqslant n$. Each equation has either no solution or exactly $q^{\min(e_i,k)}$ solutions modulo $\fp^k$, and all of them are solvable when $\ba=\mathbf0$. The last statement is the standard description of the elementary divisors through determinantal ideals.
\end{proof}

\begin{lemma}\label{lem:weyl}
Let $M$ be a symmetric $n\times n$ matrix with $2M\in M_n(\cO_{\fp})$, let $k\geqslant1$ and $\tm\in\cO_{\fp}^n$. Then
\[
\Bigl|\sum_{\bb\bmod\fp^k}\ep k\bigl(\bb^tM\bb+\tm\cdot\bb\bigr)\Bigr|\leqslant q^{nk/2}K_{\fp^k}(2M;\mathbf0)^{1/2}.
\]
\end{lemma}

\begin{proof}
This is the analogue of~\cite[Lemma~13]{BM13}. The square of the absolute value of the sum is $\sum_{\bx,\by}\ep k(\bx^tM\bx-\by^tM\by+\tm\cdot(\bx-\by))$. Substitute $\bx=\by+\bz$. Since $(\by+\bz)^tM(\by+\bz)-\by^tM\by=\by^t(2M)\bz+\bz^tM\bz$, the summand becomes $\ep k(\bz^tM\bz+\tm\cdot\bz)\ep k(\by^t(2M)\bz)$. As $\ep k$ is primitive, the sum over $\by$ is $q^{nk}$ if $2M\bz\equiv\mathbf0\bmod\fp^k$ and $0$ otherwise. Bounding the sum over $\bz$ trivially gives the claim.
\end{proof}

For the next two lemmas assume that $\fp\nmid2$. Let $\chi_{\fp}$ be the quadratic residue symbol modulo $\fp$, extended to $\cO_{\fp}$ by $\chi_{\fp}(x)=\chi_{\fp}(x\bmod\fp)$, with $\chi_{\fp}(0)=0$, and put $g_{\fp}=\sum_{x\bmod\fp}\ep1(x^2)$. Then $g_{\fp}^2=\chi_{\fp}(-1)q$.

\begin{lemma}\label{lem:gausssum}
Let $\fp\nmid2$, let $k\geqslant1$, and let $A$ be a symmetric matrix over $\cO_{\fp}$ with $\det A\in\cO_{\fp}^\times$. Put $\gamma_k(A)=\sum_{\bb\bmod\fp^k}\ep k(\bb^tA\bb)$. Then
\[
\gamma_k(A)=\begin{cases}q^{nk/2}&\text{if $k$ is even},\\ q^{n(k-1)/2}\chi_{\fp}(\det A)g_{\fp}^n&\text{if $k$ is odd},\end{cases}
\]
and, for $\tm\in\cO_{\fp}^n$,
\[
\sum_{\bb\bmod\fp^k}\ep k\bigl(\bb^tA\bb+\tm\cdot\bb\bigr)=\ep k\bigl(-\overline{4}\,\tm^tA^{-1}\tm\bigr)\gamma_k(A),
\]
where $\overline4$ is the inverse of $4$ in $\cO_{\fp}$. In particular, if $\alpha\in\cO_{\fp}^\times$ then $\gamma_k(\alpha A)=\chi_{\fp}(\alpha)^{nk}\gamma_k(A)$.
\end{lemma}

\begin{proof}
For the second formula, complete the square: since $2$ and $A$ are invertible, $\bb^tA\bb+\tm\cdot\bb=(\bb+\overline2A^{-1}\tm)^tA(\bb+\overline2A^{-1}\tm)-\overline4\,\tm^tA^{-1}\tm$, and $\bb\mapsto\bb+\overline2A^{-1}\tm$ permutes the residue classes. For the first formula, since $\fp\nmid2$ there is $T\in GL_n(\cO_{\fp})$ with $T^tAT=\operatorname{diag}(\alpha_1,\dots,\alpha_n)$ and $\alpha_i\in\cO_{\fp}^\times$. Then $\gamma_k(A)=\prod_ig_{(k)}(\alpha_i)$ with $g_{(k)}(\alpha)=\sum_{x\bmod\fp^k}\ep k(\alpha x^2)$. For $k\geqslant2$, write $x=x_0+\varpi^{k-1}t$ with $x_0\bmod\fp^{k-1}$ and $t\bmod\fp$. Then $\alpha x^2\equiv\alpha x_0^2+2\alpha\varpi^{k-1}x_0t\bmod\fp^k$, and the sum over $t$ vanishes unless $\fp\mid x_0$. Writing $x_0=\varpi x_1$ gives $g_{(k)}(\alpha)=q\,g_{(k-2)}(\alpha)$, with $g_{(0)}=1$. Also $g_{(1)}(\alpha)=\chi_{\fp}(\alpha)g_{\fp}$. Hence $g_{(k)}(\alpha)=q^{k/2}$ for even $k$ and $q^{(k-1)/2}\chi_{\fp}(\alpha)g_{\fp}$ for odd $k$. Finally $\prod_i\chi_{\fp}(\alpha_i)=\chi_{\fp}(\det A)\chi_{\fp}(\det T)^{-2}=\chi_{\fp}(\det A)$.
\end{proof}

\begin{lemma}\label{lem:ramanujan}
Let $\fp\nmid2$, $k\geqslant1$, $R\in\cO_{\fp}$ and $\iota\in\{0,1\}$, and put $\mathcal R=\sumstar_{a\bmod\fp^k}\chi_{\fp}(a)^\iota\ep k(aR)$. Then $|\mathcal R|\leqslant q^k$, and if $k\geqslant2$ then $\mathcal R=0$ unless $\fp^{k-1}\mid R$. If $\iota=0$, then $\mathcal R=q^k\mathbf 1[\fp^k\mid R]-q^{k-1}\mathbf 1[\fp^{k-1}\mid R]$. If $\iota=1$ and $k=1$, then $|\mathcal R|=q^{1/2}$ when $\fp\nmid R$ and $\mathcal R=0$ when $\fp\mid R$.
\end{lemma}

\begin{proof}
The case $\iota=0$ is the evaluation of a Ramanujan sum. Let $\iota=1$. For $k=1$, $\mathcal R=\chi_{\fp}(R)g_{\fp}$, by the substitution $a\mapsto aR^{-1}$ when $\fp\nmid R$, and $\mathcal R=\sum_a\chi_{\fp}(a)=0$ when $\fp\mid R$. For $k\geqslant2$, write $a=a_0+\varpi^{k-1}t$ with $a_0$ running over units modulo $\fp^{k-1}$ and $t\bmod\fp$. Then $\chi_{\fp}(a)=\chi_{\fp}(a_0)$, and the sum over $t$ of $\ep1(tR)$ is $q\mathbf 1[\fp\mid R]$. If $\fp\mid R$ we get $\mathcal R=q\sumstar_{a_0\bmod\fp^{k-1}}\chi_{\fp}(a_0)\ep{k-1}(a_0R/\varpi)$, and the claim follows by induction on $k$.
\end{proof}

\subsection{A bound valid at every prime}
The following bound, valid at all primes including those dividing $2\Delta_X$, replaces several of the case distinctions in~\cite[\S\S4--6]{BM13} and~\cite[\S3]{Mun15}.

\begin{lemma}\label{lem:universal}
For every prime $\fp$ there is a constant $c_{\fp}\geqslant1$, depending only on $\fp$, $Q_1$ and $Q_2$, such that
\[
|C_{\fp}(r,\ell;\bm)|\leqslant c_{\fp}\,q^{(r+\ell)(n/2+1)}
\]
for all $r,\ell\geqslant0$ and all $\bm$. If $\fp$ is good, one can take $c_{\fp}=4n$.
\end{lemma}

\begin{proof}
We may assume $r+\ell\geqslant1$. Detect the condition $\fp^r\mid Q_2(\bb)$ in Lemma~\ref{lem:dictionary} by
\[
\mathbf 1[\fp^r\mid Q_2(\bb)]=q^{-r}\sum_{c\bmod\fp^r}\ep r(cQ_2(\bb))=q^{-r}\sum_{c\bmod\fp^r}\ep{r+\ell}(c\varpi^\ell Q_2(\bb)).
\]
With $\bs=(a_1\varpi^\ell,a_2+c\varpi^\ell)$ this gives
\[
C_{\fp}(r,\ell;\bm)=q^{-r}\sumstar_{a_1\bmod\fp^r}\ \sumstar_{a_2\bmod\fp^\ell}\ \sum_{c\bmod\fp^r}\ \sum_{\bb\bmod\fp^{r+\ell}}\ep{r+\ell}\bigl(\bb^tM(\bs)\bb+\tm\cdot\bb\bigr),
\]
and Lemma~\ref{lem:weyl} bounds the inner sum by $q^{n(r+\ell)/2}K_{\fp^{r+\ell}}(2M(\bs);\mathbf0)^{1/2}$.

\emph{The case $\ell\geqslant1$.} Here $u=a_2+c\varpi^\ell$ is a unit, and as $(a_2,c)$ runs over $(\cO_{\fp}/\fp^\ell)^\times\times\cO_{\fp}/\fp^r$, $u$ runs over $(\cO_{\fp}/\fp^{r+\ell})^\times$ exactly once. Moreover $M(\bs)=uM(x,1)$ with $x=a_1\varpi^\ell/u$, and for fixed $u$, as $a_1$ runs over the units modulo $\fp^r$, so does $x/\varpi^\ell$. The matrix $2M(x,1)$ modulo $\fp^{r+\ell}$ only depends on $x/\varpi^\ell$ modulo $\fp^r$. Hence
\begin{equation}\label{eq:univ1}
|C_{\fp}(r,\ell;\bm)|\leqslant q^{-r}q^{n(r+\ell)/2}q^{r+\ell}\sumstar_{w\bmod\fp^r}K_{\fp^{r+\ell}}\bigl(2M(\varpi^\ell w,1);\mathbf0\bigr)^{1/2}.
\end{equation}
If $\fp$ is good, then $\det2M(\varpi^\ell w,1)\equiv2^n\det M_2\not\equiv0\bmod\fp$, so the $K$ factor is $1$ and we obtain $|C_{\fp}|\leqslant q^{(r+\ell)(n/2+1)}$.

\emph{The case $\ell=0$.} Now $r\geqslant1$, $M(\bs)=a_1M(1,x)$ with $x=c/a_1$, and for fixed $a_1$ the map $c\mapsto x$ is a bijection of $\cO_{\fp}/\fp^r$. Hence
\begin{equation}\label{eq:univ2}
|C_{\fp}(r,0;\bm)|\leqslant q^{-r}q^{nr/2}q^r\sum_{x\bmod\fp^r}K_{\fp^r}\bigl(2M(1,x);\mathbf0\bigr)^{1/2}.
\end{equation}
Suppose that $\fp$ is good. By~\ref{G1} modulo $\fp$, $M(1,x)$ has rank at least $n-1$ modulo $\fp$, so its first $n-1$ elementary divisors are units, and Lemma~\ref{lem:smith} gives $K_{\fp^r}(2M(1,x);\mathbf0)=q^{\min(e,r)}$ with $e=\ord_{\fp} \mathcal{P}(1,x)$. The polynomial $\mathcal{P}(1,x)$ has degree $n$ and leading coefficient $\det M_2$, a unit, and its reduction modulo $\fp$ has at most $n$ roots in $\fo/\fp$, all simple. By Hensel's lemma, each such root lifts to a unique root $\theta\in\cO_{\fp}$, and $\ord_{\fp} \mathcal{P}(1,x)=\ord_{\fp}(x-\theta)$ for $x$ congruent to it. Hence the number of $x\bmod\fp^r$ with $\ord_{\fp} \mathcal{P}(1,x)\geqslant j$ is at most $nq^{r-j}$ for $1\leqslant j\leqslant r$. Therefore
\[
\sum_{x\bmod\fp^r}K_{\fp^r}(2M(1,x);\mathbf0)^{1/2}\leqslant q^r+\sum_{j=1}^rnq^{r-j}q^{j/2}\leqslant\Bigl(1+\frac{n}{q^{1/2}-1}\Bigr)q^r\leqslant4nq^r,
\]
and~\eqref{eq:univ2} gives $|C_{\fp}(r,0;\bm)|\leqslant4nq^{r(n/2+1)}$.

\emph{Bad primes.} Let $f(x)=\det2M(x,1)=2^n\mathcal{P}(x,1)$ and $f^\ast(x)=\det2M(1,x)=2^n\mathcal{P}(1,x)$. By~\ref{G1} and~\ref{G2}, these polynomials have degree $n$, non-zero constant terms, and distinct roots in $\bar K_{\fp}$. Extend $\ord_{\fp}$ to $\bar K_{\fp}$. We claim that there are constants $C_{\fp}$ and $\eta_{\fp}>0$ such that, for every $x\in\cO_{\fp}$, either $\ord_{\fp} f(x)\leqslant C_{\fp}$, or there is a root $\theta$ of $f$ with $\ord_{\fp}(x-\theta)\geqslant\eta_{\fp}^{-1}$ and
\[
e_1+\dots+e_{n-1}\leqslant C_{\fp},\qquad e_n\leqslant\ord_{\fp}(x-\theta)+C_{\fp},
\]
where $e_1\leqslant\dots\leqslant e_n$ are the exponents of the elementary divisors of $2M(x,1)$; and similarly for $f^\ast$. Indeed, choose $\eta_{\fp}$ so small that the discs $\{\ord_{\fp}(z-\theta)\geqslant\eta_{\fp}^{-1}\}$ around the roots are disjoint. Outside these discs, $|f(x)|_{\fp}$ is bounded below on the compact set $\cO_{\fp}$. Inside the disc around $\theta$, the factors of $f$ belonging to the other roots have bounded order, so $\ord_{\fp} f(x)=\ord_{\fp}(x-\theta)+O(1)$. Moreover $M(\theta,1)$ has rank $n-1$ by~\ref{G1}, so some $(n-1)\times(n-1)$ minor $\mu(x)$ of $2M(x,1)$ has $\mu(\theta)\neq0$; shrinking $\eta_{\fp}$, $\ord_{\fp}\mu(x)$ is bounded on the disc, and Lemma~\ref{lem:smith} bounds $e_1+\dots+e_{n-1}$.

Consider first~\eqref{eq:univ1}, with $x=\varpi^\ell w$. Since $f(0)\neq0$, the roots of $f$ have bounded order; let $\ell_{\fp}$ exceed all of these orders by more than $\eta_{\fp}^{-1}$. If $\ell\geqslant\ell_{\fp}$, then $x=\varpi^\ell w$ lies outside all the discs, so $K_{\fp^{r+\ell}}\leqslant q^{C_{\fp}}$, and~\eqref{eq:univ1} gives the claim. If $\ell<\ell_{\fp}$, observe that the set of $x\in\cO_{\fp}$ with $\ord_{\fp}(x-\theta)\geqslant j$ is either empty or contained in a ball $\{\ord_{\fp}(x-x_0)\geqslant j\}$ with $x_0\in\cO_{\fp}$. Hence the number of $w\bmod\fp^r$ with $\ord_{\fp}(\varpi^\ell w-\theta)\geqslant j$ is at most $q^{r+\ell-j}$ for $0\leqslant j\leqslant r+\ell$. By Lemma~\ref{lem:smith}, $K_{\fp^{r+\ell}}\leqslant q^{C_{\fp}+\min(j+C_{\fp}+1,\,r+\ell)}$ when $j\leqslant\ord_{\fp}(x-\theta)<j+1$. Summing over $j$ and over the at most $n$ roots,
\[
\sumstar_{w\bmod\fp^r}K_{\fp^{r+\ell}}^{1/2}\ll_{\fp} q^r+\sum_{j\geqslant0}q^{r+\ell-j}q^{j/2}\ll_{\fp} q^{r+\ell}\leqslant q^{\ell_{\fp}}q^r.
\]
Inserting this into~\eqref{eq:univ1} gives the claim. The bound~\eqref{eq:univ2} is treated in the same way, using $f^\ast$, with $\ell=0$.
\end{proof}

\begin{corollary}\label{cor:pointwise}
For all $\fb_1,\fb_2$ and $\bm$ we have $|C_{\fb_1,\fb_2}(\bm)|\ll(4n)^{\omega(\fb_1\fb_2)}\N(\fb_1\fb_2)^{n/2+1}\ll_\varepsilon\N(\fb_1\fb_2)^{n/2+1+\varepsilon}$.
\end{corollary}

\begin{proof}
Combine Lemmas~\ref{lem:mult} and~\ref{lem:universal}. The finitely many bad primes contribute a constant factor, and $(4n)^{\omega(\fb)}\ll\N\fb^\varepsilon$ by Lemma~\ref{lem:smoothsupp}.
\end{proof}

\subsection{The sums for the second modulus}
We now evaluate $\mathcal Q_{\fp}(\ell;\bm)=C_{\fp}(0,\ell;\bm)$. Recall that $Q_2^\ast(\bx)=\bx^tM_2^\ast\bx$ with $M_2^\ast=4^n\operatorname{adj}(M_2)$, and that $\epsilon_n=(-1)^{(n-1)/2}$ for odd $n$.

\begin{lemma}\label{lem:gauss}
Let $\ell\geqslant1$ and $\tm=\tm^{(\ell)}$.
\begin{enumerate}[label=\textup{(\roman*)}]
\item If $\fp$ is good and $n$ is odd, then $\mathcal Q_{\fp}(1;\bm)=q^{(n+1)/2}\chi_{\fp}\bigl(\epsilon_nQ_2^\ast(\tm)\bigr)$.
\item If $\fp$ is good and $n$ is even, then $\mathcal Q_{\fp}(1;\bm)=\chi_{\fp}\bigl((-1)^{n/2}\det M_2\bigr)q^{n/2}\bigl(q\mathbf 1[\fp\mid Q_2^\ast(\tm)]-1\bigr)$.
\item If $\fp$ is good and $\ell\geqslant2$, then $\mathcal Q_{\fp}(\ell;\bm)=0$ unless $\fp^{\ell-1}\mid Q_2^\ast(\tm)$. In particular $\mathcal Q_{\fp}(\ell;\bm)=0$ if $\ell\geqslant2$ and $\fp\nmid Q_2^\ast(\tm)$.
\item If $\fp$ is good, then $|\mathcal Q_{\fp}(\ell;\bm)|\leqslant q^{\ell(n/2+1)}$. For every $\fp$, $|\mathcal Q_{\fp}(\ell;\bm)|\leqslant c_{\fp} q^{\ell(n/2+1)}$.
\end{enumerate}
\end{lemma}

\begin{proof}
These are the analogues of~\cite[Lemmas~14, 15]{BM13} and~\cite[Lemma~4]{Mun15}. Part~(iv) is contained in Lemma~\ref{lem:universal}. Let $\fp$ be good. By Lemma~\ref{lem:gausssum} with $A=aM_2$,
\[
\begin{split}
\mathcal Q_{\fp}(\ell;\bm)&=\sumstar_{a\bmod\fp^\ell}\ep\ell\bigl(-\overline{4a}\,\tm^tM_2^{-1}\tm\bigr)\gamma_\ell(aM_2)\\
&=\gamma_\ell(M_2)\sumstar_{a\bmod\fp^\ell}\chi_{\fp}(a)^{n\ell}\ep\ell(-aR),
\end{split}
\]
where $R=\overline4\,\tm^tM_2^{-1}\tm$ and we have replaced $a$ by $a^{-1}$. Since
\[
Q_2^\ast(\tm)=4^n\det(M_2)\,\tm^tM_2^{-1}\tm=4^{n+1}\det(M_2)R,
\]
and $4\det M_2$ is a unit, $R$ and $Q_2^\ast(\tm)$ have the same order. Part~(iii) now follows from Lemma~\ref{lem:ramanujan}. For~(i), $n\ell=n$ is odd, so the sum over $a$ is $\chi_{\fp}(-R)g_{\fp}$ if $\fp\nmid R$ and $0$ otherwise. With $\gamma_1(M_2)=\chi_{\fp}(\det M_2)g_{\fp}^n$ and $g_{\fp}^{n+1}=(\chi_{\fp}(-1)q)^{(n+1)/2}$,
\[
\mathcal Q_{\fp}(1;\bm)=\chi_{\fp}(-R\det M_2)\chi_{\fp}(-1)^{(n+1)/2}q^{(n+1)/2}=q^{(n+1)/2}\chi_{\fp}\bigl((-1)^{(n+3)/2}Q_2^\ast(\tm)\bigr),
\]
since $R\det M_2$ and $Q_2^\ast(\tm)$ differ by the square $4^{n+1}$. Finally $(-1)^{(n+3)/2}=\epsilon_n$. For~(ii), $n$ is even, so the sum over $a$ is the Ramanujan sum $q\mathbf 1[\fp\mid R]-1$, and $\gamma_1(M_2)=\chi_{\fp}(\det M_2)(\chi_{\fp}(-1)q)^{n/2}$.
\end{proof}

\begin{remark}[The quadratic character]\label{rem:hecke}
Let $n$ be odd and $\tv\in\fo^n$ with $Q_2^\ast(\tv)\neq0$. The function $\fp\mapsto\chi_{\fp}(\epsilon_nQ_2^\ast(\tv))$ on primes $\fp\nmid2Q_2^\ast(\tv)$ is the Artin symbol of $\fp$ in the extension $K(\sqrt{\epsilon_nQ_2^\ast(\tv)})/K$. By class field theory it is the restriction of a Hecke character $\chi_{\tv}$ of $K$ of order dividing $2$, whose finite conductor $\ff_{\tv}$ divides $4Q_2^\ast(\tv)$. The character $\chi_{\tv}$ is trivial if and only if $\epsilon_nQ_2^\ast(\tv)$ is a square in $K$. For $K=\QQ$ it is the Jacobi symbol of~\cite[Lemma~18]{BM13}.
\end{remark}

\subsection{The first modulus}
The next two lemmas give the Kloosterman refinement in the first modulus. Recall from~\eqref{eq:Phi} the binary form $\Phi_{\bm}(\bs)=-\bm^t\operatorname{adj}(M(\bs))\bm$ and the form $G$ of~\ref{G4}.

\begin{lemma}\label{lem:first-generic}
Let $\fp$ be good, let $r\geqslant1$ and $\tm=\tm^{(r)}$, and suppose that $\fp\nmid G(\tm)$. Then
\[
|C_{\fp}(1,0;\bm)|\ll q^{n/2},\qquad|C_{\fp}(r,0;\bm)|\leqslant2nq^{nr/2+r/2}\leqslant2nq^{r(n/2+1)-1}\quad(r\geqslant2),
\]
where the implied constant depends only on $n$.
\end{lemma}

\begin{proof}
This is the analogue of~\cite[Lemma~6]{Mun15}. \emph{The case $r=1$.} In Lemma~\ref{lem:dictionary} the sum over $a_1$ is $q\mathbf 1[\fp\mid Q_1(\bb)]-1$, so
\[
C_{\fp}(1,0;\bm)=q\mathcal D-\mathcal E,\qquad\mathcal D=\sum_{\substack{\bb\bmod\fp\\ Q_1\equiv Q_2\equiv0}}\ep1(\tm\cdot\bb),\qquad\mathcal E=\sum_{\substack{\bb\bmod\fp\\ Q_2\equiv0}}\ep1(\tm\cdot\bb).
\]
Since $G$ is a form of positive degree, $\fp\nmid G(\tm)$ implies $\tm\not\equiv\mathbf0\bmod\fp$. To bound $\mathcal E$, detect $Q_2(\bb)\equiv0$ by $q^{-1}\sum_{c\bmod\fp}\ep1(cQ_2(\bb))$. The term $c=0$ vanishes because $\tm\not\equiv\mathbf0$, and for $c\neq0$ the sum over $\bb$ has absolute value $q^{n/2}$ by Lemma~\ref{lem:gausssum}. Hence $|\mathcal E|\leqslant q^{n/2}$.

For $\mathcal D$, let $\hat X\subset\FF_q^n$ be the affine cone over the reduction of $X$ modulo $\fp$, and $\hat X_{\tm}$ its intersection with the hyperplane $\tm\cdot\bx=0$. Since $\hat X$ is stable under multiplication by $\FF_q^\times$,
\[
(q-1)\mathcal D=\sum_{t\in\FF_q^\times}\sum_{\bb\in\hat X}\ep1(t\,\tm\cdot\bb)=q\#\hat X_{\tm}-\#\hat X.
\]
Since $\fp$ is good, $X$ is a non-singular complete intersection of dimension $n-3$ modulo $\fp$, and, since $\fp\nmid G(\tm)$, so is $X\cap\{\tm\cdot\bx=0\}$, of dimension $n-4$ (see~\S\ref{sec:geometry}). For a non-singular complete intersection $Y$ of dimension $e$ over $\FF_q$, Deligne's theorem~\cite[Th\'eor\`eme~8.1]{Del74} gives $\#Y(\FF_q)=\#\PP^e(\FF_q)+O(q^{e/2})$, with an implied constant depending only on the dimension and the degrees. For the affine cone this gives $1+(q-1)\#Y(\FF_q)=q^{e+1}+O(q^{e/2+1})$. Hence $\#\hat X=q^{n-2}+O(q^{(n-1)/2})$ and $\#\hat X_{\tm}=q^{n-3}+O(q^{(n-2)/2})$, so that $(q-1)\mathcal D=O(q^{n/2})$. Therefore $C_{\fp}(1,0;\bm)\ll q^{n/2}$.

\emph{The case $r\geqslant2$.} As in the proof of Lemma~\ref{lem:universal}, with $x=c/a_1$ and $M_x=M(1,x)$,
\[
C_{\fp}(r,0;\bm)=q^{-r}\sum_{x\bmod\fp^r}\ \sumstar_{a_1\bmod\fp^r}T(a_1,x),\qquad T(a_1,x)=\sum_{\bb\bmod\fp^r}\ep r\bigl(a_1\bb^tM_x\bb+\tm\cdot\bb\bigr).
\]
Put $\Psi(x)=\Phi_{\tm}(1,x)$ and $\mathcal{P}_x=\mathcal{P}(1,x)=\det M_x$. We first show that, for each $x$,
\begin{equation}\label{eq:caseAB}
\Bigl|\sumstar_{a_1}T(a_1,x)\Bigr|\leqslant q^{nr/2+r+j/2}\quad\text{with }j=\min(r,\ord_{\fp} \mathcal{P}_x),
\end{equation}
and that the left-hand side vanishes unless $\ord_{\fp}\Psi(x)\geqslant r-1$, and unless $\ord_{\fp}\Psi(x)\geqslant r$ in the case $j\geqslant1$.

Suppose first that $\mathcal{P}_x$ is a unit. By Lemma~\ref{lem:gausssum} with $A=a_1M_x$,
\[
T(a_1,x)=\ep r\bigl(-\overline{4a_1}\,\tm^tM_x^{-1}\tm\bigr)\chi_{\fp}(a_1)^{nr}\gamma_r(M_x).
\]
Replacing $a_1$ by $a_1^{-1}$, the sum over $a_1$ is $\gamma_r(M_x)\mathcal R$ with $\mathcal R$ as in Lemma~\ref{lem:ramanujan}, for $R_x=-\overline4\,\tm^tM_x^{-1}\tm$. Since $|\gamma_r(M_x)|=q^{nr/2}$, \eqref{eq:caseAB} follows, and the sum vanishes unless $\fp^{r-1}\mid R_x$. Since $\Psi(x)=-\mathcal{P}_x\,\tm^t M_x^{-1}\tm=4\mathcal{P}_xR_x$ and $\mathcal{P}_x$ is a unit, this means $\ord_{\fp}\Psi(x)\geqslant r-1$.

Suppose now that $\fp\mid \mathcal{P}_x$. Since $M_x$ has rank at least $n-1$ modulo $\fp$ and $\fp\nmid2$, there is $T\in GL_n(\cO_{\fp})$ with $T^tM_xT=\operatorname{diag}(\alpha_1,\dots,\alpha_{n-1},\beta)$, where the $\alpha_i$ are units and $\ord_{\fp}\beta=\ord_{\fp} \mathcal{P}_x$. Put $\tm'=T^t\tm$. Then $T(a_1,x)$ is the product of the sums $\sum_{b}\ep r(a_1\alpha_ib^2+m_i'b)$ for $i<n$, each equal to $\ep r(-\overline{4a_1\alpha_i}m_i'^2)$ times a one-variable Gauss sum of absolute value $q^{r/2}$, and of the sum $\Sigma_n=\sum_{b\bmod\fp^r}\ep r(a_1\beta b^2+m_n'b)$. Write $b=b_0+\varpi^{r-j}t$ with $t\bmod\fp^j$. Since $\ord_{\fp}\beta\geqslant j$, we have $a_1\beta b^2\equiv a_1\beta b_0^2\bmod\fp^r$, and the sum over $t$ gives $q^j\mathbf 1[\fp^j\mid m_n']$. If $j=r$, then $\Sigma_n=q^r\mathbf 1[\fp^r\mid m_n']$. If $j<r$, write $\beta=\varpi^j\beta_0$ and $m_n'=\varpi^j\mu$; then $\Sigma_n=q^j\ep{r-j}(-\overline{4a_1\beta_0}\mu^2)$ times a Gauss sum of absolute value $q^{(r-j)/2}$. In all cases $|T(a_1,x)|\leqslant q^{nr/2+j/2}$, and
\[
T(a_1,x)=\eta\,\chi_{\fp}(a_1)^{\iota}\,\ep r\bigl(-\overline{a_1}R'_x\bigr)\,|T(a_1,x)|,\qquad R_x'=\overline4\Bigl(\sum_{i<n}\frac{m_i'^2}{\alpha_i}+\mathbf 1[j<r]\frac{\varpi^j\mu^2}{\beta_0}\Bigr),
\]
where $|\eta|=1$ and $\iota\in\{0,1\}$ do not depend on $a_1$. Lemma~\ref{lem:ramanujan} gives~\eqref{eq:caseAB}, and shows that the sum over $a_1$ vanishes unless $\fp^{j}\mid m_n'$ and $\fp^{r-1}\mid R_x'$. To relate this to $\Psi$, note that $\operatorname{adj}(M_x)=(\det T)^{-2}T\operatorname{adj}(D)T^t$ with $D=\operatorname{diag}(\alpha_1,\dots,\alpha_{n-1},\beta)$. Hence
\[
\Psi(x)=-(\det T)^{-2}\prod_{i<n}\alpha_i\cdot\Bigl(\beta\sum_{i<n}\frac{m_i'^2}{\alpha_i}+m_n'^2\Bigr).
\]
If $j<r$, the bracket equals $4\varpi^j\beta_0R_x'$, so $\ord_{\fp}\Psi(x)\geqslant j+r-1\geqslant r$. If $j=r$, the bracket has order at least $\min(r+(r-1),2r)\geqslant r$. This proves the claim.

Since $\fp$ is good and $\fp\nmid G(\tm)$, the binary form $\Phi_{\tm}$ has distinct roots modulo $\fp$, by~\ref{G4} and the choice of $\Delta_X$. So $\Psi$ is not identically zero modulo $\fp$, and its roots modulo $\fp$ are simple. By Hensel's lemma, the number of $x\bmod\fp^r$ with $\ord_{\fp}\Psi(x)\geqslant k$ is at most $(n-1)q^{\max(r-k,0)}$. Using~\eqref{eq:caseAB} with $j=0$ for the at most $(n-1)q$ values of $x$ with $\mathcal{P}_x$ a unit and $\ord_{\fp}\Psi(x)\geqslant r-1$, and with $j\leqslant r$ for the at most $n-1$ values of $x$ with $\ord_{\fp}\Psi(x)\geqslant r$, we obtain
\[
|C_{\fp}(r,0;\bm)|\leqslant q^{-r}(n-1)\bigl(q\cdot q^{nr/2+r}+q^{nr/2+r+r/2}\bigr)\leqslant2nq^{nr/2+r/2},
\]
since $q\leqslant q^{r/2}$ for $r\geqslant2$. Finally $nr/2+r/2\leqslant r(n/2+1)-1$ for $r\geqslant2$.
\end{proof}

\begin{lemma}\label{lem:first-degenerate}
Let $\fp$ be good. Then $|C_{\fp}(1,0;\bm)|\leqslant3nq^{(n+1)/2}$ for every $\bm$.
\end{lemma}

\begin{proof}
We use the notation of the previous proof with $r=1$, so that $C_{\fp}(1,0;\bm)=q^{-1}\sum_{x\bmod\fp}\sumstar_{a_1}T(a_1,x)$. For the at most $n$ values of $x$ with $\fp\mid \mathcal{P}_x$, the bound~\eqref{eq:caseAB} with $j=1$ gives a total contribution of at most $q^{-1}nq^{n/2+3/2}=nq^{(n+1)/2}$. For the other $x$, $\sumstar_{a_1}T(a_1,x)=\gamma_1(M_x)\mathcal R(x)$, where $\mathcal R(x)=\sumstar_a\chi_{\fp}(a)^n\ep1(aR_x)$.

If $n$ is odd, $|\mathcal R(x)|\leqslant q^{1/2}$ by Lemma~\ref{lem:ramanujan}, and the contribution is at most $q^{-1}\cdot q\cdot q^{n/2}q^{1/2}$.

If $n$ is even, $\gamma_1(M_x)=\chi_{\fp}(\mathcal{P}_x)\gamma$ with $\gamma=(\chi_{\fp}(-1)q)^{n/2}$ independent of $x$, and $\mathcal R(x)=q\mathbf 1[\fp\mid R_x]-1$. The contribution is therefore
\[
q^{-1}\gamma\Bigl(q\sum_{\substack{x:\ \fp\nmid \mathcal{P}_x\\ \fp\mid R_x}}\chi_{\fp}(\mathcal{P}_x)-\sum_{x:\ \fp\nmid \mathcal{P}_x}\chi_{\fp}(\mathcal{P}_x)\Bigr).
\]
The polynomial $\mathcal{P}(1,x)$ has degree $n$ and distinct roots modulo $\fp$, so it is not a constant multiple of a square, and Weil's bound gives $|\sum_{x\bmod\fp}\chi_{\fp}(\mathcal{P}(1,x))|\leqslant(n-1)q^{1/2}$. If $\Psi=\Phi_{\tm}(1,\cdot)$ is not identically zero modulo $\fp$, then $\fp\mid R_x$ for at most $n-1$ values of $x$, since $\Psi(x)=4\mathcal{P}_xR_x$. If $\Psi\equiv0\bmod\fp$, then $\fp\mid R_x$ for every $x$ with $\fp\nmid \mathcal{P}_x$, and the first sum is again bounded by Weil's bound. In either case the contribution is at most $q^{-1}q^{n/2}(q(n-1)q^{1/2}+(n-1)q^{1/2})\leqslant2(n-1)q^{(n+1)/2}$.
\end{proof}

\subsection{Mixed moduli}
\begin{lemma}\label{lem:CtoS}
Let $r,\ell\geqslant1$ and $\tm=\tm^{(r+\ell)}$. Then
\[
C_{\fp}(r,\ell;\bm)=q^rS_{\fp}(r,\ell;\bm)-q^{r-2}S_{\fp}(r-1,\ell+1;\bm).
\]
\end{lemma}

\begin{proof}
This corresponds to~\cite[Lemma~3]{Mun15}, whose statement contains a third term. In the range $r,\ell\geqslant1$, which is the only one we need, the direct computation below shows that only two terms occur. In Lemma~\ref{lem:dictionary}, $\ep{r+\ell}(a_1\varpi^\ell Q_1(\bb))=\ep r(a_1Q_1(\bb))$, and the sum over $a_1$ is the Ramanujan sum $q^r\mathbf 1[\fp^r\mid Q_1(\bb)]-q^{r-1}\mathbf 1[\fp^{r-1}\mid Q_1(\bb)]$. The first part gives $q^rS_{\fp}(r,\ell;\bm)$. For the second, let
\[
A_k=\sum_{\substack{\bb\bmod\fp^{r+\ell}\\ \fp^{r-1}\mid Q_1(\bb),\ \fp^{k}\mid Q_2(\bb)}}\ep{r+\ell}(\tm\cdot\bb)\qquad(k\geqslant r).
\]
Since $\fp^r\mid Q_2(\bb)$ there, the sum over $a_2$ of $\ep{r+\ell}(a_2Q_2(\bb))=\ep\ell(a_2Q_2(\bb)/\varpi^r)$ is the Ramanujan sum $q^\ell\mathbf 1[\fp^{r+\ell}\mid Q_2]-q^{\ell-1}\mathbf 1[\fp^{r+\ell-1}\mid Q_2]$. So the second part equals $-q^{r-1}(q^\ell A_{r+\ell}-q^{\ell-1}A_{r+\ell-1})$. On the other hand, in $S_{\fp}(r-1,\ell+1;\bm)$ we have $\fp^{r-1}\mid Q_2(\bb)$, and the sum over $a\bmod\fp^{\ell+1}$ is $q^{\ell+1}\mathbf 1[\fp^{r+\ell}\mid Q_2]-q^\ell\mathbf 1[\fp^{r+\ell-1}\mid Q_2]$. Thus $S_{\fp}(r-1,\ell+1;\bm)=q^{\ell+1}A_{r+\ell}-q^\ell A_{r+\ell-1}$, and the second part equals $-q^{r-2}S_{\fp}(r-1,\ell+1;\bm)$.
\end{proof}

\begin{lemma}\label{lem:Svanish}
Let $\fp$ be good, let $r,\ell\geqslant1$ and $\tm=\tm^{(r+\ell)}$. Then $S_{\fp}(r,\ell;\bm)=0$ unless $\fp^{\min(r,\ell)}\mid Q_2^\ast(\tm)$.
\end{lemma}

\begin{proof}
This is the analogue, at good primes, of~\cite[Lemmas~26, 27]{BM13}. Write $\bb=\bk+\varpi^r\bx$ with $\bk\bmod\fp^r$ and $\bx\bmod\fp^\ell$. The conditions $\fp^r\mid Q_i(\bb)$ depend only on $\bk$, since $Q_i(\bk+\varpi^r\bx)=Q_i(\bk)+\varpi^r\bx^t(2M_i)\bk+\varpi^{2r}Q_i(\bx)$. The phase is
\[
\ep{r+\ell}\bigl(aQ_2(\bk)+\tm\cdot\bk\bigr)\,\ep\ell\bigl(\varpi^raQ_2(\bx)+(a\nabla Q_2(\bk)+\tm)\cdot\bx\bigr),
\]
where $\nabla Q_2(\bk)=2M_2\bk$. Put $k=\min(r,\ell)$ and write $\bx=\by+\varpi^{\ell-k}\bz$ with $\bz\bmod\fp^k$. Then $\varpi^rQ_2(\bx)\equiv\varpi^rQ_2(\by)\bmod\fp^\ell$, because $r+\ell-k\geqslant\ell$, and the sum over $\bz$ is $q^{nk}\mathbf 1[a\nabla Q_2(\bk)+\tm\equiv\mathbf0\bmod\fp^k]$. Hence $S_{\fp}(r,\ell;\bm)=0$ unless there are a unit $a$ and $\bk$ with $\fp^r\mid Q_2(\bk)$ and $2aM_2\bk\equiv-\tm\bmod\fp^k$. In that case $\bk\equiv-\overline{2a}M_2^{-1}\tm\bmod\fp^k$, and therefore
\[
Q_2^\ast(\tm)=4^n\det(M_2)\,\tm^tM_2^{-1}\tm\equiv4^{n+1}a^2\det(M_2)\,Q_2(\bk)\equiv0\bmod\fp^k.\qedhere
\]
\end{proof}

\begin{lemma}\label{lem:mixed}
Let $\fp$ be good, let $r,\ell\geqslant1$ and $\tm=\tm^{(r+\ell)}$. Then $C_{\fp}(r,\ell;\bm)=0$ unless $\fp\mid Q_2^\ast(\tm)$, and in any case $|C_{\fp}(r,\ell;\bm)|\leqslant q^{(r+\ell)(n/2+1)}$.
\end{lemma}

\begin{proof}
Compare~\cite[Lemma~9]{Mun15}. The bound is Lemma~\ref{lem:universal}. By Lemma~\ref{lem:CtoS}, $C_{\fp}(r,\ell;\bm)$ is a combination of $S_{\fp}(r,\ell;\bm)$ and $S_{\fp}(r-1,\ell+1;\bm)$. The first vanishes unless $\fp\mid Q_2^\ast(\tm)$ by Lemma~\ref{lem:Svanish}. So does the second: for $r\geqslant2$ by Lemma~\ref{lem:Svanish}, and for $r=1$ because $S_{\fp}(0,\ell+1;\bm)=\mathcal Q_{\fp}(\ell+1;\bm)$ with $\ell+1\geqslant2$, by Lemma~\ref{lem:gauss}(iii). Note that $\tm^{(r+\ell)}$ is the same normalised frequency for both terms.
\end{proof}

\subsection{Averages over the first modulus}
For $\tv\in\fo^n\setminus\{\mathbf0\}$ let $\fg(\tv)$ be the principal ideal generated by $\Delta_XG(\tv)$ if $G(\tv)\neq0$, and by $\Delta_X$ if $G(\tv)=0$.

\begin{lemma}\label{lem:Mun8}
Let $\tv\in\fo^n\setminus\{\mathbf0\}$. Suppose that for every prime $\fp$ and every $\fb_1$ in the sums below, the normalised frequency $\tm^{(r_{\fp})}$ of $\bm$ at $\fp$ is a unit multiple of $\tv$. Then, for $x\geqslant1$,
\[
\sum_{\substack{\N\fb_1\leqslant x\\ (\fb_1,\fg(\tv))=1}}\frac{|C_{\fb_1,\fo}(\bm)|}{\N\fb_1^{n/2+1}}\ll_\varepsilon\begin{cases}x^\varepsilon&\text{if }G(\tv)\neq0,\\ x^{1/2+\varepsilon}&\text{if }G(\tv)=0.\end{cases}
\]
\end{lemma}

\begin{proof}
This is the analogue of~\cite[Lemma~8]{Mun15}. All primes dividing $\fb_1$ are good. If $G(\tv)\neq0$, then $\fp\nmid G(\tm)$ for every $\fp\mid\fb_1$, by homogeneity. Lemmas~\ref{lem:mult} and~\ref{lem:first-generic} give $|C_{\fb_1,\fo}(\bm)|\leqslant(Cn)^{\omega(\fb_1)}\N\fb_1^{n/2+1}/\N\rad\fb_1$, and Lemma~\ref{lem:smoothsupp} gives the claim. If $G(\tv)=0$, Lemma~\ref{lem:first-degenerate} for the primes exactly dividing $\fb_1$ and Lemma~\ref{lem:universal} for the others give
\[
\frac{|C_{\fb_1,\fo}(\bm)|}{\N\fb_1^{n/2+1}}\leqslant(4n)^{\omega(\fb_1)}\N(\fs)^{-1/2},
\]
where $\fs$ is the square-free part of $\fb_1$, that is, the product of the primes exactly dividing $\fb_1$. Write $\fb_1=\fs\ff$ with $\ff$ square-full and coprime to $\fs$. The number of square-full $\ff$ with $\N\ff\leqslant y$ is $O(y^{1/2})$, so
\[
\sum_{\N\fb_1\leqslant x}(4n)^{\omega(\fb_1)}\N(\fs)^{-1/2}\ll x^\varepsilon\sum_{\N\ff\leqslant x}\ \sum_{\N\fs\leqslant x/\N\ff}\N\fs^{-1/2}\ll x^{1/2+\varepsilon}\sum_{\ff\ \text{square-full}}\N\ff^{-1/2-\varepsilon}\ll x^{1/2+2\varepsilon}.\qedhere
\]
\end{proof}

\subsection{The average over the second modulus}\label{sec:hecke}
This subsection replaces~\cite[Lemma~5]{Mun15} and~\cite[Lemma~18]{BM13}. Over $\QQ$ one sums a Jacobi symbol over the integers $q_2$ in a dyadic range. Over $K$ the dual lattice $(\fb_1\fb_2)^\wedge$ moves with $\fb_2$, and in Section~\ref{sec:error} we shall fix it by restricting $\fb_2$ to an ideal class and summing over all generators $\beta$ of $\fb_2\fc$ for a fixed ideal $\fc$. This leads to sums over elements $\beta$ weighted by a function of $\beta\in V^\times$, which we treat by Fourier analysis on $V^\times/U_K$ and Hecke $L$-functions.

For each place $l$, let $\mathcal D_l$ be the set of differential operators $\beta^{(l)}\partial_{\beta^{(l)}}$ (if $l$ is real) or $\beta^{(l)}\partial_{\beta^{(l)}}$ and $\bar\beta^{(l)}\partial_{\bar\beta^{(l)}}$ (if $l$ is complex), acting on functions of $\beta\in V^\times$. These operators commute with multiplication by any fixed element of $V^\times$, in particular by units. For a function $F$ on $V^\times$ we put $F^U(\beta)=\sum_{u\in U_K}F(u\beta)$ whenever the sum converges absolutely.

\begin{proposition}\label{prop:hecke}
Let $n$ be odd, let $\fc'$ be a non-zero integral ideal, and let $\cS$ be a finite set of primes containing those dividing $2\fc'$. Let $\chi$ be a Hecke character of $K$ of order at most $2$ whose finite conductor $\ff$ is divisible only by primes in $\cS$. Let $X\geqslant1$, let $B>0$, and let $F$ be a smooth function on $V^\times$, supported in $\{X/2\leqslant|\Nm\beta|/\N\fc'\leqslant2X\}$, such that $F^U$ converges absolutely, $\sup|F^U|\leqslant B$, and
\[
\sup_\beta|D_1\cdots D_jF^U(\beta)|\leqslant P^{j\varepsilon}B\qquad(D_i\in\textstyle\bigcup_l\mathcal D_l,\ j\leqslant J)
\]
for a sufficiently large $J$ depending on $d$ and $\varepsilon$. Suppose that $X$, $\N\fc'$ and $\N\ff$ are at most $P^{O(1)}$, and that every prime in $\cS$ divides a fixed integral ideal $\fa_{\cS}$ with $\N\fa_{\cS}\leqslant P^{O(1)}$. For $\beta\in\fc'\setminus\{0\}$ write $\fb=(\beta)\fc'^{-1}$. Then
\[
\Sigma:=\sum_{\substack{\beta\in\fc'\setminus\{0\}\\ (\fb,\cS)=1}}\mu^2(\fb)\chi(\fb)(\N\fb)^{\frac{n+1}2}F(\beta)\ll P^{\varepsilon}B\,\mathcal E,
\]
where $\mathcal E=X^{(n+2)/2}(\N\ff)^{1/4-\delta_{\mathrm{MV}}}$ if $\chi\psi\neq1$ for every character $\psi$ of the class group, and $\mathcal E=X^{(n+3)/2}$ in general.
\end{proposition}

\begin{proof}
\emph{Step 1: grouping by ideals.} The summand depends on $\beta$ only through $\fb$, apart from $F(\beta)$. Grouping the elements $\beta$ according to the principal ideal $(\beta)$, we get
\[
\Sigma=\sum_{\substack{\fb\ \text{integral},\ (\fb,\cS)=1\\ \fb\fc'\ \text{principal}}}\mu^2(\fb)\chi(\fb)(\N\fb)^{\frac{n+1}2}F^U(\beta_{\fb}),
\]
where $\beta_{\fb}$ is any generator of $\fb\fc'$; this is well defined because $F^U$ is invariant under $U_K$.

\emph{Step 2: Fourier and Mellin expansion.} Every $\beta\in V^\times$ can be written uniquely as $\beta=|\Nm\beta|^{1/d}\beta_1$ with $\beta_1\in V^1=\{|\Nm\beta|=1\}$, where $|\Nm\beta|^{1/d}$ is embedded diagonally. This decomposition commutes with multiplication by units. The group $\mathbb T=V^1/U_K$ is a compact abelian Lie group of dimension $d-1$. Let $\lambda$ run over its characters, and let $\|\lambda\|$ be the size of the (integral and real) parameters describing $\lambda$; the number of $\lambda$ with $\|\lambda\|\leqslant T$ is $O((1+T)^{d-1})$. Fourier expansion on $\mathbb T$ and Mellin inversion in $|\Nm\beta|$ give, for $\sigma\in\RR$,
\[
F^U(\beta)=\sum_\lambda\lambda(\beta_1)\frac1{2\pi i}\int_{(\sigma)}\tilde F(\lambda,s)|\Nm\beta|^{-s}\,ds,
\]
where
\[
\tilde F(\lambda,s)=\int_{V^\times/U_K}F^U(\beta)\overline{\lambda(\beta_1)}|\Nm\beta|^s\,d^\times\beta,
\]
for a suitably normalised Haar measure $d^\times\beta$. The functions $\beta\mapsto\lambda(\beta_1)|\Nm\beta|^{-s}$ are joint eigenfunctions of the operators in $\bigcup_l\mathcal D_l$, with eigenvalues that are linear forms in $s$ and the parameters of $\lambda$, and some eigenvalue has absolute value $\gg1+|s|+\|\lambda\|$. Integrating by parts $J'$ times with the corresponding operator, and using the support of $F$, we obtain
\[
\tilde F(\lambda,s)\ll(X\N\fc')^{\Re s}B\Bigl(\frac{P^\varepsilon}{1+|s|+\|\lambda\|}\Bigr)^{J'}\qquad(J'\leqslant J).
\]
Hence the terms with $\|\lambda\|>P^{2\varepsilon}$, and the parts of the integrals with $|\Im s|>P^{2\varepsilon}$, contribute $O(P^{-10}B)$ to $\Sigma$ in total, once $J$ is large enough.

\emph{Step 3: Hecke $L$-functions.} For a character $\lambda$ of $\mathbb T$, the map $(\beta)\mapsto\lambda(\beta_1)$ is a character of the group of principal fractional ideals, since $\lambda$ is trivial on units. Since this group has finite index $h_K$ in the group of all fractional ideals, it extends to a character $\tilde\lambda$ of the latter, which is a Hecke character unramified at every finite place, with infinity type determined by $\lambda$. Two extensions differ by a character of the class group. Detect the condition that $\fb\fc'$ is principal by $h_K^{-1}\sum_\psi\psi(\fb\fc')$, where $\psi$ runs over the characters of the class group. Since $\lambda((\beta_{\fb})_1)=\tilde\lambda(\fb\fc')$ and $|\Nm\beta_{\fb}|=\N\fb\N\fc'$, we obtain, for $\Re s=\frac{n+3}2+\varepsilon$,
\[
\Sigma=\sum_\lambda\frac1{2\pi i}\int_{(\sigma)}\frac{\tilde F(\lambda,s)}{\N\fc'^{\,s}}\cdot\frac1{h_K}\sum_\psi(\tilde\lambda\psi)(\fc')\,\frac{L^{\cS}(s-\frac{n+1}2,\xi)}{L^{\cS}(2s-n-1,\xi^2)}\,ds,\qquad\xi=\chi\psi\tilde\lambda.
\]
Here $L^{\cS}$ denotes an $L$-function with the Euler factors at $\cS$ removed, and we have used the identity $\sum_{(\fb,\cS)=1}\mu^2(\fb)\xi(\fb)\N\fb^{-w}=\prod_{\fp\notin\cS}(1+\xi(\fp)\N\fp^{-w})=L^{\cS}(w,\xi)/L^{\cS}(2w,\xi^2)$.

\emph{Step 4: contour shift.} Move the line of integration to $\Re s=\frac{n+2}2$, that is, $\Re(s-\frac{n+1}2)=\frac12$. The function $L(w,\xi)$ has a pole, at $w=1$, only if $\xi$ is trivial, which requires $\tilde\lambda$ to have trivial infinity type, so that $\lambda=1$ and $\tilde\lambda$ is a class group character, and $\chi\psi\tilde\lambda=1$. If $\chi\psi\neq1$ for every class group character $\psi$, this never happens. On the new line we use the following bounds for $|\Im s|\leqslant P^{2\varepsilon}$ and $\|\lambda\|\leqslant P^{2\varepsilon}$:
\begin{itemize}
\item the removed Euler factors satisfy
\[
\bigl|L^{\cS}(w,\xi)/L(w,\xi)\bigr|\leqslant2^{|\cS|},\qquad\bigl|L(2w,\xi^2)/L^{\cS}(2w,\xi^2)\bigr|\leqslant2^{|\cS|};
\]
moreover $2^{|\cS|}\leqslant2^{\omega(\fa_{\cS})}$, which is $O(P^\varepsilon)$ by Lemma~\ref{lem:smoothsupp};
\item $|L(1+2it,\xi^2)|^{-1}\ll P^\varepsilon$, by the standard zero-free region and Siegel's theorem for the possible exceptional real character (and trivially if $\xi^2$ is trivial, since then $L=\zeta_K$ has a pole at $1$);
\item the analytic conductor of $\xi$ is $\ll\N\ff P^{O(\varepsilon)}$, because $\psi$ and $\tilde\lambda$ are unramified at the finite places and the archimedean parameters are $O(P^{2\varepsilon})$; so Theorem~\ref{thm:MV} gives $L(\frac12+it,\xi)\ll(\N\ff)^{1/4-\delta_{\mathrm{MV}}}P^{O(\varepsilon)}$.
\end{itemize}
Since $(X\N\fc')^{(n+2)/2}\N\fc'^{-(n+2)/2}=X^{(n+2)/2}$, and the number of $\lambda$ involved is $P^{O(\varepsilon)}$, the integrals on the new line contribute $\ll P^\varepsilon BX^{(n+2)/2}(\N\ff)^{1/4-\delta_{\mathrm{MV}}}$. This proves the first bound. The second bound follows by keeping the line $\Re s=\frac{n+3}2+\varepsilon$ and bounding everything trivially.
\end{proof}

\begin{remark}\label{rem:unramified}
For $K=\QQ$, the condition $\chi\psi=1$ for some $\psi$ means that $\epsilon_nQ_2^\ast(\bv)$ is a square, which is Munshi's exceptional case. Over $K$ with $h_K>1$, it also happens when $K(\sqrt{\epsilon_nQ_2^\ast(\bv)})/K$ is a non-trivial \emph{everywhere unramified} quadratic extension, that is, a subextension of the Hilbert class field. Let $\mathcal A_K\subset K^\times/K^{\times2}$ be the finite set of classes $a$ such that $K(\sqrt a)$ is contained in the Hilbert class field; it contains the trivial class. The exceptional set is $\{\bv:\epsilon_nQ_2^\ast(\bv)\in\mathcal A_K\}$. In Section~\ref{sec:error} it is treated together with the case of squares.
\end{remark}

For even $n$, and for the case $Q_2^\ast(\tv)=0$, we only need trivial bounds.

\begin{lemma}\label{lem:betatrivial}
Let $\fc'$, $X$, $F$ and $B$ be as in Proposition~\ref{prop:hecke}, except that we only assume $\sup_\beta\sum_{u\in U_K}|F(u\beta)|\leqslant B$. Let $\kappa\geqslant0$. Then
\[
\sum_{\beta\in\fc'\setminus\{0\}}(\N\fb)^\kappa|F(\beta)|\ll BX^{\kappa+1},\qquad\sum_{\substack{\beta\in\fc'\setminus\{0\}\\ \fb\ \text{square-full}}}(\N\fb)^\kappa|F(\beta)|\ll BX^{\kappa+1/2}.
\]
\end{lemma}

\begin{proof}
As in Step~1 above, the sums equal sums over integral ideals $\fb$ with $\N\fb\asymp X$ of $(\N\fb)^\kappa\sum_u|F(u\beta_{\fb})|$. There are $O(X)$ such ideals, and $O(X^{1/2})$ square-full ones.
\end{proof}

\section{The exponential integrals}\label{sec:integrals}

In this section we study the integrals $I_{c_1,\fb_2}(\bm)$ defined in~\eqref{eq:Iint}. Throughout, $\fb_1$ and $\fb_2$ are non-zero integral ideals with $\N\fb_1\leqslant(A_0P)^d$ and $\N\fb_2\leqslant c_{\ast}P^{d/2}$, and $c_1\in\fb_1\setminus\{0\}$ satisfies $\ip{c_1}\leqslant A_0P$. We put
\[
\rho_1=\frac{\N\fb_1}{P^d},\qquad\rho_2=\frac{\N\fb_2}{P^{d/2}},\qquad\bl=\bl(\fb_1,c_1).
\]
Then $P^{-d}\leqslant\rho_1\leqslant A_0^d$ and $P^{-d/2}\leqslant\rho_2\leqslant c_{\ast}$. By Lemma~\ref{lem:lambda},
\begin{equation}\label{eq:lambdabounds}
\Nm(\bl)=\rho_1^{-1},\qquad\frac{P}{|c_1^{(l)}|}\leqslant\bl^{(l)}\ll P^{d+2},\qquad\bl^{(l)}\geqslant A_0^{-1}\qquad(1\leqslant l\leqslant r_1+r_2).
\end{equation}
We always assume that $P$ is larger than a suitable constant. A quantity is called \emph{negligible} if it is $O_A(P^{-A})$ for every $A>0$.

\subsection{A class of integrals}
In Section~\ref{sec:error} we shall need to differentiate $I_{c_1,\fb_2}(\bm)$ with respect to the second modulus. The derivatives are integrals of the same shape with slightly different weights. To handle them simultaneously we work with the following class.

\begin{definition}\label{def:class}
\begin{enumerate}[label=(\alph*)]
\item A \emph{shell weight} is a function of the form
\[
k(x,y)=\frac1{\Delta_K}\sum_{(0)\neq\fz\subseteq\fo}\frac1{x\N\fz}\Bigl\{a(x\N\fz)-b\Bigl(\frac{|y|}{x\N\fz}\Bigr)\Bigr\}\qquad(x>0,\ y\in\RR),
\]
where $a,b\in C^\infty(\RR)$ are supported in $[\frac12,1]$.
\item A \emph{product cut-off} is a function $u_V(v)=\prod_lu_l(|v^{(l)}|)$ on $V$, where each $u_l\in C^\infty([0,\infty))$ is supported in $[0,2]$ and is either constant near $0$ or vanishes near $0$.
\item An \emph{element weight} (attached to $c_1$) is a function $A(v)=\kappa_0-\prod_la_l(|Pv^{(l)}/c_1^{(l)}|)$ on $V$, where $|\kappa_0|\leqslant1$ and each $a_l\in C^\infty(\RR)$ is supported in $[\frac12,1]$.
\item A \emph{Gaussian weight of degree $D$} is a function on $V^n$ of the form
\[
W^\natural(\by)=\Pi(\by-\bxi)\,w_0(\by)\,e^{-L\|\by-\bxi\|^2}+R(\by),
\]
where $\Pi$ is a polynomial in the $dn$ real coordinates of degree at most $D$ whose coefficients are $O(L^D)$, and $R\in C^\infty(V^n)$ is supported in $\{\ip{\by-\bxi}\leqslant\delta_0\}$ and satisfies $\|R\|_{C^k}\ll_kL^{D+k}e^{-c_0L}$ for every $k$, with $c_0=\delta_0^2/16$.
\end{enumerate}
An \emph{integral of type $\mathfrak J_D$} is an integral
\begin{equation}\label{eq:typeJ}
I^\natural(\bm)=\int_{V^n}W^\natural(\by)\,A\bigl(Q_1(\by)\bigr)\,k\bigl(\rho_2,\Nm(\bl Q_2(\by))\bigr)\,u_V\bigl(\bl Q_2(\by)\bigr)\,\ee(-P\bm\cdot\by)\,d\by,
\end{equation}
where $W^\natural$ is a Gaussian weight of degree $D$, $A$ is an element weight, $k$ is a shell weight and $u_V$ is a product cut-off.
\end{definition}

By~\eqref{eq:hdef}, $h$ is the shell weight with $a=b=w$, and $U_V$ is a product cut-off. With $\kappa_0=\bW(c_1/P)$ and $a_l=w_1$, the weight $\bW(c_1/P)-\bW(Q_1/(c_1P))$ in~\eqref{eq:Iint} is an element weight. Finally $W$ is a Gaussian weight of degree $0$, with $\Pi=1$ and $R=0$. So $I_{c_1,\fb_2}(\bm)$ is an integral of type $\mathfrak J_0$.

In the rest of this section, the implied constants depend on $D$, on the functions $a,b,u_l,a_l$ that enter the definition (through finitely many of their derivatives), and on the data fixed in Section~\ref{sec:nested}. They do not depend on $P$, $\fb_1$, $c_1$, $\fb_2$, $\bm$, $\kappa_0$ or on the coefficients of $\Pi$ beyond the bound $O(L^D)$.

We first record some elementary properties.

\begin{lemma}\label{lem:classbasic}
Let $k$ be a shell weight.
\begin{enumerate}[label=\textup{(\roman*)}]
\item $k$ is smooth on $(0,\infty)\times\RR$, and $k(x,y)=0$ unless $x\leqslant\max(1,2|y|)$.
\item Write $k=k_1-k_2$, where $k_1(x)$ is the contribution of $a$. Then $k_1(x)\ll x^{-1}$, and $k_1(x)=0$ for $x>1$. Moreover $k_2(x,y)=0$ unless $|y|\geqslant x/2$, and for $j\geqslant0$ we have $\partial_y^jk_2(x,y)\ll_jx^{-1}|y|^{-j}$.
\item In particular $|k(x,y)|\ll x^{-1}$ and $|\partial_y^jk(x,y)|\ll_jx^{-1-j}$ for $j\geqslant1$.
\item The functions $x\partial_xk$ and $y\partial_yk$ are shell weights.
\end{enumerate}
\end{lemma}

\begin{proof}
For (i) and (ii), note that the term of $k_1$ indexed by $\fz$ is non-zero only if $\frac12\leqslant x\N\fz\leqslant1$, and that of $k_2$ only if $x\N\fz\leqslant|y|\leqslant2x\N\fz$. Since $\N\fz\geqslant1$, the first requires $x\leqslant1$ and the second requires $|y|\geqslant x$. Both conditions define locally finite families, and the functions $y\mapsto b(|y|/(x\N\fz))$ are smooth because $b$ vanishes near $0$. The number of ideals $\fz$ with $\N\fz$ in a dyadic interval $[T,2T]$ is $O(1+T)$. For $k_1$ we have $T\asymp x^{-1}$ and each term is $O(1)$, so $k_1\ll1+x^{-1}\ll x^{-1}$. For $k_2$ we have $T\asymp|y|/x\gg1$ and each term is $O((x\N\fz)^{-1})=O(|y|^{-1})$, while $\partial_y^j$ of the term is $O(|y|^{-1}(x\N\fz)^{-j})=O(|y|^{-1-j})$. Hence $\partial_y^jk_2\ll(|y|/x)|y|^{-1-j}=x^{-1}|y|^{-j}$. Part~(iii) follows, since $|y|\geqslant x/2$ on the support of $k_2$ and $k_1$ does not depend on $y$. For~(iv), a direct computation with $t=\N\fz$ gives
\begin{gather*}
x\partial_x\Bigl(\frac{a(xt)}{xt}\Bigr)=\frac{(xt)a'(xt)-a(xt)}{xt},\\
x\partial_x\Bigl(\frac{b(u)}{xt}\Bigr)=\frac{-b(u)-ub'(u)}{xt},\qquad y\partial_y\Bigl(\frac{b(u)}{xt}\Bigr)=\frac{ub'(u)}{xt},
\end{gather*}
where $u=|y|/(xt)$. The new functions $za'(z)-a(z)$, $-b(z)-zb'(z)$ and $zb'(z)$ are again smooth and supported in $[\frac12,1]$.
\end{proof}

\subsection{Polynomial truncation}
\begin{lemma}\label{lem:ibp}
Let $I^\natural$ be an integral of type $\mathfrak J_D$. Then $|I^\natural(\bm)|\ll L^D\rho_2^{-1}$ for all $\bm$, and for every $N\geqslant0$,
\[
I^\natural(\bm)\ll_NP^{C_0N}\bigl(P\ip\bm\bigr)^{-N},
\]
where $C_0$ depends only on $d$. In particular $I^\natural(\bm)$ is negligible when $\ip\bm\geqslant P^{2C_0}$, and the same holds for $I_{\fb_1,\fb_2}(\bm)$.
\end{lemma}

\begin{proof}
On the support of the integrand we have $\ip{\by-\bxi}\leqslant\delta_0$, so the integrand is supported in a fixed compact set. There $|W^\natural|\ll L^D$, $|A|\leqslant2$, $|u_V|\ll1$ and, by Lemma~\ref{lem:classbasic}(iii), $|k|\ll\rho_2^{-1}$. This gives the first bound.

For the second, we estimate the partial derivatives of the integrand with respect to the $dn$ real coordinates of $\by$. A derivative of $W^\natural$ costs a factor $O(L)$. A derivative of $a_l(|PQ_1^{(l)}(\by)/c_1^{(l)}|)$ costs $O(P/|c_1^{(l)}|)$, which is $O(P^{d+2})$ by~\eqref{eq:lambdabounds}. On the support of $u_V(\bl Q_2(\by))$ we have $|\bl^{(l)}Q_2^{(l)}(\by)|\leqslant2$ for every $l$, so every partial derivative of $\Nm(\bl Q_2(\by))$ is $O(\max_l\bl^{(l)})=O(P^{d+2})$. By Lemma~\ref{lem:classbasic}(iii), every derivative of $k(\rho_2,\cdot)$ costs $O(\rho_2^{-1})=O(P^{d/2})$. Finally a derivative of $u_V(\bl Q_2(\by))$ costs $O(P^{d+2})$. By the product and chain rules, any partial derivative of order $N$ of the integrand, apart from the phase, is $O_N(P^{(2d+2)N}\rho_2^{-1})$, since $\frac d2+d+2\leqslant2d+2$.

Choose a real coordinate $y_\ast$ of $\by$ such that the corresponding coefficient of the linear form $\by\mapsto\Tr(P\bm\cdot\by)$ has absolute value $\gg P\ip\bm$. This is possible because, at a complex place, a complex number $z$ satisfies $|z|\leqslant\sqrt2\max(|\Re z|,|\Im z|)$, and the trace doubles the real part. Integrating by parts $N$ times with respect to $y_\ast$, which is legitimate since the integrand has compact support, we obtain the second bound with $C_0=2d+2+d/2$. For $I_{\fb_1,\fb_2}(\bm)$ we use~\eqref{eq:Ibig}: the sum over $c_1$ has $O(P^{d})$ terms, each with $|\N c_1|^{-1}\leqslant1$, and $P^{2d}/C_P\ll P^d$.
\end{proof}

\subsection{The Fourier transforms of the weights}
The next lemma is the basic estimate for the second equation. It replaces the bounds of~\cite[Lemmas~6.1--6.4]{BV14}, from which it differs only in the method of proof: we decompose the weight into pieces supported on hyperbolic shells $\{|\Nm v|\asymp T\}$.

For $T>0$ and a product cut-off $u_V$, let $b\in C^\infty(\RR)$ be supported in $[\frac12,1]$ and put
\[
g_T(v)=u_V(v)\,b\bigl(|\Nm v|/T\bigr)\qquad(v\in V).
\]
Since $b$ vanishes near $0$, $g_T$ is smooth. It is supported in the shell
\[
\mathcal S_T=\{v\in V:\ \ip v\leqslant2,\ T/2\leqslant|\Nm v|\leqslant T\},
\]
which is empty if $T>2^{d+1}$. Note that every $v\in\mathcal S_T$ has $v^{(l)}\neq0$ for all $l$. Write $r=r_1+r_2$.

\begin{lemma}[Shells]\label{lem:shell}
Let $0<T\leqslant2^{d+1}$ and $N\geqslant0$. Then
\[
\hat g_T(t)\ll_NT\bigl(1+|\log T|\bigr)^{r-1}\min\Bigl\{1,\bigl(T\cH(t)\bigr)^{-N}\Bigr\}.
\]
\end{lemma}

\begin{proof}
We use the coordinates $r_l=|v^{(l)}|^{c_l}$ and the Euler operators $\theta_l=r_l\partial_{r_l}$. A function of $v$ that depends only on $(r_l)_l$ will be written as a function of these variables. The function $g_T$ is of this form, and so is each of the functions obtained from it by applying any monomial in the $\theta_l$. All of these are bounded, with bounds depending on the monomial, and supported in $\mathcal S_T$. Indeed $\theta_l$ maps $u_l(r_l^{1/c_l})$ to $c_l^{-1}su_l'(s)|_{s=r_l^{1/c_l}}$, and it maps $b(\prod_{l'}r_{l'}/T)$ to $zb'(z)|_{z=\prod r_{l'}/T}$.

For each place $l$ let $\mathsf D_l$ denote $\partial/\partial v^{(l)}$ if $l$ is real, and the Laplacian $\partial_x^2+\partial_y^2$ in $v^{(l)}=x+iy$ if $l$ is complex. We claim that, for any smooth function $F$ of $r_l$ alone,
\begin{equation}\label{eq:euler}
\mathsf D_l^N\bigl[F(r_l)\bigr]=\eta_N(v^{(l)})\,r_l^{-N}\,(\mathcal P_N(\theta_l)F)(r_l),
\end{equation}
where $|\eta_N|=1$ and $\mathcal P_N$ is a polynomial of degree at most $c_lN$ depending only on $N$ and $c_l$. At a real place, $r=|v|$ and $\partial_v[v^{-k}G(|v|)]=v^{-k-1}((\theta-k)G)(|v|)$, which gives~\eqref{eq:euler} by induction with $\eta_N=\operatorname{sgn}(v)^N$. At a complex place, $r=|v|^2$, and for a function $G$ of $r$ we have $\partial_x^2G+\partial_y^2G=4G'+4rG''=4r^{-1}\theta^2G$. Applying this to $r^{-k}G$, for which $\theta(r^{-k}G)=r^{-k}(\theta-k)G$, gives $(\partial_x^2+\partial_y^2)(r^{-k}G)=4r^{-k-1}(\theta-k)^2G$. Again~\eqref{eq:euler} follows by induction.

Let $S=\{l:|t^{(l)}|>1\}$, so that $\cH(t)=\prod_{l\in S}|t^{(l)}|^{c_l}$. If $l$ is real then $\mathsf D_l\ee(-tv)=-2\pi it^{(l)}\ee(-tv)$. If $l$ is complex then $\ee(-tv)$ contains the factor $\exp(-4\pi i\Re(t^{(l)}v^{(l)}))$, and $\mathsf D_l\ee(-tv)=-16\pi^2|t^{(l)}|^2\ee(-tv)$. Since $g_T$ has compact support, we may integrate by parts $N$ times with each $\mathsf D_l$, $l\in S$. By~\eqref{eq:euler}, applied successively at the places of $S$ (the operators act on different variables and commute),
\[
|\hat g_T(t)|\ll_N\cH(t)^{-N}\int_{\mathcal S_T}\prod_{l\in S}r_l^{-N}\,dv.
\]
This holds also for $N=0$. It remains to show that the last integral is $O(T^{1-N}(1+|\log T|)^{r-1})$.

Pass to the variables $u_l=\log r_l$, and the angular variables at the complex places. Then $dv\asymp\prod_le^{u_l}\,du_l$ times the angular measure. On $\mathcal S_T$ we have $u_l\leqslant c_l\log2$ for each $l$, and $\log(T/2)\leqslant\sum_lu_l\leqslant\log T$. Consequently each $u_l$ lies in the interval $[\log T-C,C]$ with $C=d\log2+1$, which has length $O(1+|\log T|)$. If $N\geqslant1$, then
\[
\prod_{l\in S}e^{-(N-1)u_l}=e^{-(N-1)\sum_{l\in S}u_l}\leqslant\Bigl(\frac2T\Bigr)^{N-1}\prod_{l\notin S}e^{(N-1)u_l},
\]
so the integrand $\prod_{l\in S}e^{-Nu_l}\prod_le^{u_l}$ is at most $(2/T)^{N-1}\prod_{l\notin S}e^{Nu_l}$. If $N=0$ the integrand is $\prod_le^{u_l}=e^{\sum u_l}\leqslant T$. Now integrate. Choose one variable $u_{l_0}$, with $l_0\in S$ if $S\neq\emptyset$. For fixed values of the others it ranges over an interval of length $\log2$, because of the constraint on $\sum u_l$. The variables $u_l$ with $l\notin S\cup\{l_0\}$ contribute $\int_{-\infty}^Ce^{Nu}\,du\ll1$ if $N\geqslant1$, and $O(1+|\log T|)$ if $N=0$. The variables $u_l$ with $l\in S\setminus\{l_0\}$ each contribute $O(1+|\log T|)$. In all cases the integral is $O(T^{1-N}(1+|\log T|)^{r-1})$.
\end{proof}

\begin{lemma}[The second weight]\label{lem:p2}
Let $k$ be a shell weight, $u_V$ a product cut-off and $0<\rho\leqslant c_{\ast}$, and put
\[
p_\rho(t)=\int_Vu_V(v)\,k(\rho,\Nm v)\,\ee(-tv)\,dv.
\]
Then, for every $N\geqslant0$,
\[
p_\rho(t)\ll_N\frac{(1+|\log\rho|)^{r_1+r_2}}{\rho\,\cH(t)}\min\Bigl\{1,\bigl(\rho\cH(t)\bigr)^{-N}\Bigr\}.
\]
In particular, $p_\rho(t)\ll\rho^{-1}(1+|\log\rho|)^{r_1+r_2}\cH(t)^{-1}$, and $p_\rho(t)\ll_NP^{-N}\cH(t)^{-1}$ if $\cH(t)\geqslant\rho^{-1}P^\varepsilon$.
\end{lemma}

\begin{proof}
Write $k=k_1-k_2$ as in Lemma~\ref{lem:classbasic}. The function $k_1(\rho)$ does not depend on $v$, so its contribution is $k_1(\rho)\hat u_V(t)$. Since $u_V\in C_c^\infty(V)$, we have $\hat u_V(t)\ll_N\prod_l(1+|t^{(l)}|)^{-N}\leqslant\cH(t)^{-N}$, and $k_1(\rho)\ll\rho^{-1}$. This is acceptable.

By definition,
\[
k_2(\rho,\Nm v)=\frac1{\Delta_K}\sum_{\fz}\frac1{\rho\N\fz}\,b\Bigl(\frac{|\Nm v|}{\rho\N\fz}\Bigr),
\]
and the contribution of $k_2$ to $p_\rho(t)$ is $-\Delta_K^{-1}\sum_{\fz}(\rho\N\fz)^{-1}\hat g_{\rho\N\fz}(t)$, with $g_T$ as above. Only $\fz$ with $\rho\N\fz\leqslant2^{d+1}$ occur. Group the ideals $\fz$ according to the dyadic interval $[T,2T)$ containing $\rho\N\fz$. Since $\N\fz\geqslant1$, only $T\geqslant\rho/2$ occur, and each dyadic range contains $O(1+T/\rho)=O(T/\rho)$ ideals. By Lemma~\ref{lem:shell}, the contribution of $k_2$ is therefore
\[
\ll_N\frac{(1+|\log\rho|)^{r-1}}{\rho}\sum_{\substack{\rho/2\leqslant T\leqslant2^{d+1}\\ T\ \text{dyadic}}}T\min\Bigl\{1,(T\cH(t))^{-N-1}\Bigr\}.
\]
The terms with $T\leqslant\cH(t)^{-1}$ contribute $O(\cH(t)^{-1})$, and the terms with $T>\cH(t)^{-1}$ contribute $O(\cH(t)^{-1})$ as well, since they form a geometric series. If moreover $\rho\cH(t)\geqslant1$, all terms have $T\geqslant\rho/2$ and the sum is $O(\rho(\rho\cH(t))^{-N-1})$. In either case the sum is $\ll\cH(t)^{-1}\min\{1,(\rho\cH(t))^{-N}\}$, which proves the lemma; the factor $(1+|\log\rho|)$ accounts for the number of dyadic ranges.
\end{proof}

We turn to the first equation. Fix $w_2\in C_c^\infty(V)$ with $w_2(v)=1$ for $\ip v\leqslant A_0/2$ and $w_2(v)=0$ for $\ip v\geqslant A_0$. Then $w_2(Q_1(\by))=1$ on $\{\ip{\by-\bxi}\leqslant\delta_0\}$, by the choice of $A_0$ in~\eqref{eq:Ibig}. For an element weight $A$ attached to $c_1$, put
\[
q(t)=\int_Vw_2(v)A(v)\ee(-tv)\,dv,\qquad\varpi_l(t)=\min\bigl(1,|t|^{-c_l}\bigr)\quad(t\in K_l).
\]

\begin{lemma}[The first weight]\label{lem:p1}
For every $N\geqslant0$,
\[
q(t)\ll_N\prod_l\bigl(1+|t^{(l)}|\bigr)^{-N}+\prod_l\Bigl(\frac{|c_1^{(l)}|}P\Bigr)^{c_l}\Bigl(1+\frac{|c_1^{(l)}||t^{(l)}|}P\Bigr)^{-N}.
\]
Consequently:
\begin{enumerate}[label=\textup{(\roman*)}]
\item $|q(t)|\ll\prod_l\varpi_l(t^{(l)})$ and $\int_V|q(t)|\,dt\ll1$;
\item $q(t)\ll_NP^{-N}\bigl(\prod_l(1+|t^{(l)}|)^{-2}+\prod_l(|c_1^{(l)}|/P)^{c_l}(1+|c_1^{(l)}||t^{(l)}|/P)^{-2}\bigr)$ unless $|t^{(l)}|\leqslant P^{1+\varepsilon}/|c_1^{(l)}|$ for every $l$.
\end{enumerate}
\end{lemma}

\begin{proof}
The term $\kappa_0w_2$ has Fourier transform $\kappa_0\hat w_2(t)$, which gives the first term. The other term is the Fourier transform of $g(v)=w_2(v)\prod_la_l(|Pv^{(l)}/c_1^{(l)}|)$. This function is supported in $\{|v^{(l)}|\leqslant|c_1^{(l)}|/P\ \forall l\}$, a box of volume $\ll\prod_l(|c_1^{(l)}|/P)^{c_l}$. Its partial derivatives of order $N$ in the coordinates of $v^{(l)}$ are $\ll_N(P/|c_1^{(l)}|)^N$, since $P/|c_1^{(l)}|\geqslant A_0^{-1}$ and the derivatives of $w_2$ are bounded. Integrating by parts $N$ times at each place $l$ with $|c_1^{(l)}||t^{(l)}|/P>1$, in the real coordinate of $v^{(l)}$ in which $t^{(l)}$ is largest, gives the second term.

For~(i), consider a single factor of the second term. If $|t|\leqslant P/|c_1^{(l)}|$, then $(|c_1^{(l)}|/P)^{c_l}\leqslant\min(A_0^{c_l},|t|^{-c_l})$. If $|t|>P/|c_1^{(l)}|$, then, taking $N=c_l$, the factor is at most $(|c_1^{(l)}|/P)^{c_l}(P/(|c_1^{(l)}||t|))^{c_l}=|t|^{-c_l}$, and also at most $A_0^{c_l}$. The first term is trivially $\ll\prod_l\varpi_l$. The bound for $\int|q|$ follows with $N=3$, after the substitution $t^{(l)}\mapsto Pt^{(l)}/|c_1^{(l)}|$ in the second term. For~(ii), suppose that $|t^{(l_0)}|>P^{1+\varepsilon}/|c_1^{(l_0)}|$. Then $|t^{(l_0)}|\geqslant P^\varepsilon/A_0$, and in the first term the factor at $l_0$ is $O(P^{-\varepsilon(N-2)}(1+|t^{(l_0)}|)^{-2})$. In the second term, the factor at $l_0$ is at most $P^{-\varepsilon(N-2)}$ times the same factor with $N=2$. Replacing $N$ by $N/\varepsilon+2$ gives~(ii).
\end{proof}

\begin{lemma}[Fourier inversion]\label{lem:inversion}
Let $I^\natural$ be of type $\mathfrak J_D$, with ingredients $W^\natural,A,k,u_V$, and let $p_{\rho_2}$ be as in Lemma~\ref{lem:p2} for $k$ and $u_V$. For $\bs=(s_1,s_2)\in V^2$ and $\bu\in V^n$ put
\[
\cK^\natural(\bs;\bu)=\int_{V^n}W^\natural(\by)\,\ee\bigl(s_1Q_1(\by)+s_2Q_2(\by)-\bu\cdot\by\bigr)\,d\by.
\]
Then
\begin{equation}\label{eq:fourier}
I^\natural(\bm)=\iint_{V^2}q(t_1)\,p_{\rho_2}(t_2)\,\cK^\natural(t_1,\bl t_2;P\bm)\,dt_1\,dt_2,
\end{equation}
where $\bl t_2$ is the componentwise product, and the integral converges absolutely.
\end{lemma}

\begin{proof}
The functions $w_2A$ and $G(v)=u_V(v)k(\rho_2,\Nm v)$ lie in $C_c^\infty(V)$: the first because each $a_l$ vanishes near $0$, the second by Lemma~\ref{lem:classbasic}(i). By Fourier inversion,
\[
A(Q_1(\by))=w_2(Q_1(\by))A(Q_1(\by))=\int_Vq(t_1)\ee(t_1Q_1(\by))\,dt_1
\]
for $\by$ in the support of $W^\natural$, and
\[
k\bigl(\rho_2,\Nm(\bl Q_2(\by))\bigr)u_V\bigl(\bl Q_2(\by)\bigr)=G(\bl Q_2(\by))=\int_Vp_{\rho_2}(t_2)\ee\bigl(t_2\bl Q_2(\by)\bigr)\,dt_2.
\]
Insert these into~\eqref{eq:typeJ}. The functions $q$ and $p_{\rho_2}$ are integrable by Lemmas~\ref{lem:p1} and~\ref{lem:p2} (using Lemma~\ref{lem:Hint} for the latter), and $W^\natural$ is integrable. So Fubini's theorem applies.
\end{proof}

\subsection{Gaussian integrals}
We now evaluate $\cK^\natural$. It is essentially a product over the places of Gaussian integrals with a quadratic phase.

\begin{lemma}\label{lem:gaussreal}
Let $N\geqslant1$, $L\geqslant1$, let $B$ be a real symmetric $N\times N$ matrix with eigenvalues $\beta_1,\dots,\beta_N$, let $\mathbf w\in\RR^N$ and let $\alpha\in\ZZ_{\geqslant0}^N$. Put
\[
E=\mathbf w^t(L^2+4\pi^2B^2)^{-1}\mathbf w.
\]
Then
\[
\Bigl|\int_{\RR^N}\bz^\alpha e^{-L|\bz|^2}e\bigl(\bz^tB\bz+\mathbf w\cdot\bz\bigr)\,d\bz\Bigr|\ll_{\alpha,N}L^{-|\alpha|/2}\prod_{j=1}^N(L^2+4\pi^2\beta_j^2)^{-1/4}\exp\Bigl(-\frac{\pi^2L}2E\Bigr).
\]
For $\alpha=0$ the left-hand side equals $\pi^{N/2}\prod_j(L^2+4\pi^2\beta_j^2)^{-1/4}\exp(-\pi^2LE)$.
\end{lemma}

\begin{proof}
Let $O$ be an orthogonal matrix with $O^tBO=\operatorname{diag}(\beta_j)$, and substitute $\bz=O\bz'$. Then $\bz^\alpha$ becomes a homogeneous polynomial of degree $|\alpha|$ in $\bz'$ with coefficients bounded in terms of $\alpha$, and $|\bz|^2$, $\mathbf w^t(L^2+4\pi^2B^2)^{-1}\mathbf w$ are unchanged if we replace $B$ by $O^tBO$ and $\mathbf w$ by $O^t\mathbf w$. It therefore suffices to consider $B$ diagonal and a monomial $\bz^\alpha$. Then the integral factorises into one-dimensional integrals
\[
I_a(\beta,\omega)=\int_\RR z^ae^{-Lz^2}e(\beta z^2+\omega z)\,dz=\int_\RR z^ae^{-\mathcal Az^2+bz}\,dz,\qquad\mathcal A=L-2\pi i\beta,\quad b=2\pi i\omega.
\]
For $a=0$ the classical formula gives $I_0=\sqrt{\pi/\mathcal A}\,e^{b^2/(4\mathcal A)}$, with the principal branch of the square root, and
\[
\Re\frac{b^2}{4\mathcal A}=-\pi^2\omega^2\Re\frac1{\mathcal A}=-\frac{\pi^2L\omega^2}{|\mathcal A|^2},\qquad|\mathcal A|^2=L^2+4\pi^2\beta^2.
\]
This proves the case $\alpha=0$. In general $I_a=\partial_b^aI_0$, and by induction $\partial_b^ae^{b^2/(4\mathcal A)}$ is $e^{b^2/(4\mathcal A)}$ times a linear combination of the terms $(b/2\mathcal A)^{a-2i}(2\mathcal A)^{-i}$, $0\leqslant i\leqslant a/2$. Put $x=|\omega|/|\mathcal A|$. Then $|b/(2\mathcal A)|=\pi x$ and $|2\mathcal A|^{-1}\leqslant(2L)^{-1}$, so that
\[
|I_a|\ll_a|\mathcal A|^{-1/2}e^{-\pi^2Lx^2}\sum_{i\leqslant a/2}x^{a-2i}L^{-i}\ll_aL^{-a/2}|\mathcal A|^{-1/2}e^{-\pi^2Lx^2/2},
\]
since $x^je^{-\pi^2Lx^2/2}\ll_jL^{-j/2}$. Multiplying over the coordinates gives the lemma, because $\sum_j\omega_j^2/|\mathcal A_j|^2=E$.
\end{proof}

For a place $l$, $\bs\in K_l^2$ and $\bu\in K_l^n$, and a multi-index $\alpha$ in the real coordinates of $K_l^n$, put
\[
J_{l,\alpha}(\bs;\bu)=\int_{K_l^n}(\by-\bxi^{(l)})^\alpha\exp\bigl(-c_lL|\by-\bxi^{(l)}|^2\bigr)\,\ee_l\bigl(s_1Q_1(\by)+s_2Q_2(\by)-\bu\cdot\by\bigr)\,d\by,
\]
where $\ee_l(z)=\exp(2\pi ic_l\Re z)$, so that $\ee(v)=\prod_l\ee_l(v^{(l)})$. Recall from Lemma~\ref{lem:spectrum} the singular values $\sigma_1(\bs)\leqslant\dots\leqslant\sigma_n(\bs)$ of $M(\bs)$, and put
\[
\Theta_l(\bs)=\prod_{j=1}^n\bigl(L^2+\sigma_j(\bs)^2\bigr)^{-c_l/4},\qquad\bw_l(\bs,\bu)=2M(\bs)\bxi^{(l)}-\bu=s_1\nabla Q_1(\bxi^{(l)})+s_2\nabla Q_2(\bxi^{(l)})-\bu.
\]

\begin{lemma}\label{lem:gaussian}
There is a constant $c_3>0$, depending only on $Q_1$, $Q_2$ and $\bxi$, such that
\[
|J_{l,\alpha}(\bs;\bu)|\ll_\alpha L^{-|\alpha|/2}\,\Theta_l(\bs)\exp\Bigl(-c_3\frac{L|\bw_l(\bs,\bu)|^2}{L^2+|\bs|^2}\Bigr).
\]
\end{lemma}

\begin{proof}
Substitute $\by=\bxi^{(l)}+\bz$. Then $s_1Q_1(\by)+s_2Q_2(\by)-\bu\cdot\by$ equals a constant plus $\bz^tM(\bs)\bz+\bw_l\cdot\bz$.

At a real place this is Lemma~\ref{lem:gaussreal} with $N=n$, $B=M(\bs)$ and $\mathbf w=\bw_l$. Here $|\beta_j|=\sigma_j(\bs)$, and $(L^2+4\pi^2B^2)^{-1}\geqslant(L^2+4\pi^2\|M(\bs)\|^2)^{-1}$ with $\|M(\bs)\|\ll|\bs|$.

At a complex place, identify $\CC^n$ with $\RR^{2n}$. The phase is $2\Re(\bz^tM(\bs)\bz)+2\Re(\bw_l\cdot\bz)$, and the Gaussian is $e^{-2L|\bz|^2}$. By Takagi's factorisation there is a unitary matrix $U$ with $U^tM(\bs)U=\operatorname{diag}(\sigma_j(\bs))$. The substitution $\bz=U\bz'$ is orthogonal on $\RR^{2n}$, and with $z_j'=x_j+iy_j$ the quadratic part becomes $\sum_j2\sigma_j(x_j^2-y_j^2)$. So we may apply Lemma~\ref{lem:gaussreal} with $N=2n$, with $2L$ in place of $L$, and with a real symmetric matrix $B$ with eigenvalues $\pm2\sigma_j(\bs)$; the vector $\mathbf w$ has length $2|\bw_l|$. Since each $\sigma_j$ appears twice, the product of the factors $((2L)^2+16\pi^2\sigma_j^2)^{-1/4}$ is $\asymp\Theta_l(\bs)$.
\end{proof}

\begin{lemma}\label{lem:Kfactor}
Let $W^\natural$ be a Gaussian weight of degree $D$. Then there are coefficients $\pi_\alpha\ll L^D$, indexed by multi-indices $\alpha=(\alpha_l)_l$ with $|\alpha|\leqslant D$, such that for all $\bs\in V^2$ and $\bu\in V^n$,
\[
\cK^\natural(\bs;\bu)=\sum_\alpha\pi_\alpha\prod_lJ_{l,\alpha_l}(\bs^{(l)};\bu^{(l)})+O\bigl(L^{2D}e^{-c_0L}\bigr).
\]
\end{lemma}

\begin{proof}
Write $\Pi(\by-\bxi)=\sum_\alpha\pi_\alpha\prod_l(\by^{(l)}-\bxi^{(l)})^{\alpha_l}$. Since $\|\by-\bxi\|^2=\sum_lc_l|\by^{(l)}-\bxi^{(l)}|^2$, the function $\Pi(\by-\bxi)e^{-L\|\by-\bxi\|^2}$ is a sum of products over the places, and so is the phase. Replacing $w_0$ by $1$ changes the integral by at most $\int_{\ip{\by-\bxi}\geqslant\delta_0/2}|\Pi|e^{-L\|\by-\bxi\|^2}\,d\by\ll L^De^{-c_0L}$, since $\|\by-\bxi\|^2\geqslant\ip{\by-\bxi}^2$. The term $R$ contributes $O(L^De^{-c_0L})$.
\end{proof}

\subsection{Stationary phase on average over the pencil}
The heart of the archimedean analysis is the following estimate. It shows that, after averaging over the pencil with the weights coming from Lemmas~\ref{lem:p2} and~\ref{lem:p1}, the Gaussian integral decays like $|\bu|^{-c_ln/2}$, as stationary phase suggests. We cannot argue pointwise in $\bs$: when $[\bs]$ is close to a root of $\det M(\bs)$, the matrix $M(\bs)$ is badly conditioned and the pointwise bound is weaker. The average takes care of this.

\begin{proposition}\label{prop:statphase}
Let $l$ be a place, let $\lambda\geqslant A_0^{-1}$, let $\bu\in K_l^n$, and let $C_1\geqslant1$. Then
\[
\iint_{\substack{t_1,t_2\in K_l\\ |t_1|,|t_2|\leqslant P^{C_1}}}\varpi_l(t_1)\varpi_l(t_2)\,\bigl|J_{l,\alpha}(t_1,\lambda t_2;\bu)\bigr|\,dt_1\,dt_2\ll_{\alpha,C_1}(\log P)^2\min\bigl(1,|\bu|^{-c_ln/2}\bigr).
\]
\end{proposition}

\begin{proof}
Put $\bs=(s_1,s_2)=(t_1,\lambda t_2)$. Then $dt_2=\lambda^{-c_l}ds_2$, and since $\lambda\gg1$ we have
\[
\varpi_l(t_2)\lambda^{-c_l}=\min\bigl(\lambda^{-c_l},|s_2|^{-c_l}\bigr)\ll\varpi_l(s_2).
\] The range $|t_2|\leqslant P^{C_1}$ becomes $|s_2|\leqslant\lambda P^{C_1}\leqslant P^{C_1+d+3}$. So it suffices to bound
\[
\mathcal I=\iint_{|\bs|\leqslant P^{C_2}}\varpi_l(s_1)\varpi_l(s_2)\,|J_{l,\alpha}(\bs;\bu)|\,d\bs,\qquad C_2=C_1+d+3.
\]
We have $\iint_{|\bs|\leqslant P^{C_2}}\varpi_l(s_1)\varpi_l(s_2)\,d\bs\ll(\log P)^2$. By Lemma~\ref{lem:gaussian} we have $|J_{l,\alpha}|\ll L^{-|\alpha|/2}\Theta_l(\bs)\leqslant L^{-c_ln/2}$, and since $L^{-|\alpha|/2}\leqslant1$ we may ignore this factor from now on.

\emph{Localisation.} Since $\sqrt L=(\log P)^2$, Lemma~\ref{lem:gaussian} shows that $|J_{l,\alpha}(\bs;\bu)|\ll e^{-c_3(\log P)^2}$ unless
\begin{equation}\label{eq:localise}
|\bw_l(\bs,\bu)|\leqslant\frac{L+|\bs|}{\log P}.
\end{equation}
The part of $\mathcal I$ where~\eqref{eq:localise} fails is $O(P^{2c_lC_2}e^{-c_3(\log P)^2})$, which is negligible. Let $T:K_l^2\to K_l^n$ be the linear map $\bs\mapsto s_1\nabla Q_1(\bxi^{(l)})+s_2\nabla Q_2(\bxi^{(l)})$, so that $\bw_l=T\bs-\bu$. Since $\bxi\notin X^\ast(V)$, the two gradients are linearly independent over $K_l$, so $T$ is injective and $|T\bs|\asymp|\bs|$, uniformly for $\bxi$ in a compact subset of $V^n\setminus X^\ast(V)$. In particular $|\bs-\bs'|\ll|T\bs-T\bs'|$.

\emph{Small $\bu$.} If $|\bu|\leqslant C_4L$ for a suitable constant $C_4$, we use $|J_{l,\alpha}|\ll L^{-c_ln/2}\ll\min(1,|\bu|^{-c_ln/2})$ and the bound $(\log P)^2$ for the integral of the weights.

\emph{Large $\bu$.} Suppose $|\bu|>C_4L$. If~\eqref{eq:localise} holds, then $|T\bs-\bu|\leqslant(L+C|T\bs|)/\log P$, which for $P$ large forces $|T\bs|\asymp|\bu|$, and therefore $|\bs|\asymp|\bu|$ and $|T\bs-\bu|\ll|\bu|/\log P$. Fix one $\bs_0$ satisfying~\eqref{eq:localise} (if there is none, $\mathcal I$ is negligible). Then every $\bs$ satisfying~\eqref{eq:localise} lies in the disc $\mathcal D=\{|\bs-\bs_0|\leqslant R\}$ with $R\ll|\bu|/\log P$. In particular all $\bs\in\mathcal D$ have $|\bs|\asymp|\bu|$.

\emph{The determinant factor.} By Lemma~\ref{lem:spectrum}, $\sigma_j(\bs)\asymp|\bs|$ for $j\geqslant2$, and $\sigma_1(\bs)\asymp|\bs|\operatorname{dist}([\bs],Z)$. Hence
\[
\Theta_l(\bs)\ll|\bs|^{-c_l(n-1)/2}\max\bigl(L,\sigma_1(\bs)\bigr)^{-c_l/2}\leqslant|\bs|^{-c_l(n-1)/2}\sigma_1(\bs)^{-c_l/2}.
\]
Let $\cZ_l\subset K_l^2$ be the union of the lines $K_l\bs$ with $[\bs]\in Z$ (only real roots contribute when $l$ is real). For $|\bs|\asymp|\bu|$ we have $|\bs|\operatorname{dist}([\bs],Z)\asymp\operatorname{dist}(\bs,\cZ_l)$, where the right-hand side is the Euclidean distance; for real $\bs$ and non-real roots the distance to $Z$ is bounded below, and those roots may be ignored.

By~\ref{G2}, the set $Z$ contains neither $[1:0]$ nor $[0:1]$. Fix a closed conic neighbourhood $\mathcal N$ of $\cZ_l$ which still avoids the two coordinate axes; on $\mathcal N$ we have $|s_1|\asymp|s_2|\asymp|\bs|$. Since $R/|\bu|\ll1/\log P$, for $P$ large the disc $\mathcal D$ either lies in $\mathcal N$ or does not meet $\frac12\mathcal N$, the conic neighbourhood of half the angular width.

\emph{Case 1: $\mathcal D\cap\frac12\mathcal N=\emptyset$.} Then $\operatorname{dist}([\bs],Z)\gg1$ on $\mathcal D$, so $\sigma_1(\bs)\asymp|\bs|\asymp|\bu|$ and $\Theta_l(\bs)\ll|\bu|^{-c_ln/2}$. The weights integrate to $O((\log P)^2)$, which gives the claim.

\emph{Case 2: $\mathcal D\subset\mathcal N$.} Then $\varpi_l(s_1)\varpi_l(s_2)\ll|\bu|^{-2c_l}$ on $\mathcal D$, and
\[
\mathcal I\ll|\bu|^{-2c_l}|\bu|^{-c_l(n-1)/2}\int_{\mathcal D}\operatorname{dist}(\bs,\cZ_l)^{-c_l/2}\,d\bs+\text{negligible}.
\]
If $l$ is real, $\cZ_l$ is a finite union of lines in $\RR^2$, and $\int_{\mathcal D}\operatorname{dist}(\bs,\ell)^{-1/2}\,d\bs\ll R\int_0^{2R}\delta^{-1/2}\,d\delta\ll R^{3/2}$ for each line $\ell$. If $l$ is complex, $\cZ_l$ is a finite union of complex lines in $\CC^2\cong\RR^4$, each of real codimension two; writing $\bs$ in coordinates tangent and normal to such a line, $\int_{\mathcal D}\operatorname{dist}(\bs,\ell)^{-1}\,d\bs\ll R^2\int_{|\nu|\leqslant2R}|\nu|^{-1}\,d\nu\ll R^3$, the inner integral being over $\nu\in\RR^2$. Since $R\ll|\bu|$, we obtain $\mathcal I\ll|\bu|^{-2-(n-1)/2+3/2}=|\bu|^{-n/2}$ if $l$ is real, and $\mathcal I\ll|\bu|^{-4-(n-1)+3}=|\bu|^{-n}$ if $l$ is complex. In both cases this is $|\bu|^{-c_ln/2}$.
\end{proof}

\begin{remark}
Proposition~\ref{prop:statphase} replaces the estimate of the measure of a set of bad $(t_1,t_2)$ in the proof of~\cite[Lemma~11]{Mun15}. That argument assumes that $M(t_1,t_2)$ is uniformly well conditioned, which fails near the singular members of the pencil. The averaged estimate is robust and works in the same way at real and complex places.
\end{remark}

\subsection{Bounds for the integrals}
\begin{proposition}\label{cor:Ibound}
Let $I^\natural$ be an integral of type $\mathfrak J_D$. Then, for all $\bm\in V^n$,
\[
I^\natural(\bm)\ll L^D\rho_2^{-1}(\log P)^{C}\cH(P\bm)^{-n/2},
\]
where $C$ depends only on $d$. In particular
\[
I_{\fb_1,\fb_2}(\bm)\ll(\log P)^{C}\frac{P^{3d/2}}{\N(\fb_1\fb_2)}\cH(P\bm)^{-n/2}.
\]
\end{proposition}

\begin{proof}
If $\ip\bm\geqslant P^{3C_0}$, then Lemma~\ref{lem:ibp} with $N$ large gives $I^\natural(\bm)\ll\ip\bm^{-N/2}\leqslant(P\ip\bm)^{-dn/2}\leqslant\cH(P\bm)^{-n/2}$. So we may assume that $\ip\bm<P^{3C_0}$, so that $\cH(P\bm)\leqslant P^{O(1)}$.

Start from~\eqref{eq:fourier}. By Lemma~\ref{lem:p1}(ii), the part of the integral where $|t_1^{(l)}|>P^{1+\varepsilon}/|c_1^{(l)}|$ for some $l$ is $\ll P^{-N}\int|p_{\rho_2}|\sup|\cK^\natural|\ll P^{-N}$, since $\int|p_{\rho_2}|\ll\rho_2^{-1}(\log P)^{C}$ by Lemmas~\ref{lem:p2} and~\ref{lem:Hint}, and $|\cK^\natural|\leqslant\int|W^\natural|\ll L^D$. Similarly, by Lemma~\ref{lem:p2} the part where $\cH(t_2)>\rho_2^{-1}P^\varepsilon$ is negligible. In the remaining range all $|t_1^{(l)}|$ and $|t_2^{(l)}|$ are at most $P^{C_1}$ for a suitable $C_1$. There, by Lemmas~\ref{lem:p1}(i) and~\ref{lem:p2},
\[
|q(t_1)p_{\rho_2}(t_2)|\ll\rho_2^{-1}(\log P)^{r_1+r_2}\prod_l\varpi_l(t_1^{(l)})\varpi_l(t_2^{(l)}).
\]
Insert Lemma~\ref{lem:Kfactor}. The error term there contributes $\ll\rho_2^{-1}(\log P)^CL^{2D}e^{-c_0L}$, which is negligible. The main term is a sum of $O(1)$ products over the places, with coefficients $O(L^D)$. Since the weights $\prod_l\varpi_l\varpi_l$ also factorise, Fubini's theorem and Proposition~\ref{prop:statphase} (with $\lambda=\bl^{(l)}$ and $\bu=P\bm^{(l)}$) give
\[
\ll L^D\rho_2^{-1}(\log P)^{r_1+r_2}\prod_l(\log P)^2\min\bigl(1,|P\bm^{(l)}|^{-c_ln/2}\bigr)=L^D\rho_2^{-1}(\log P)^{3(r_1+r_2)}\cH(P\bm)^{-n/2}.
\]
Since $\cH(P\bm)\leqslant P^{O(1)}$, the negligible terms are also bounded by the right-hand side. For the second assertion, apply the first with $D=0$ to each $I_{c_1,\fb_2}$ in~\eqref{eq:Ibig}, and use Lemma~\ref{lem:latticesum} together with $P^{2d}/C_P\ll P^d$:
\[
I_{\fb_1,\fb_2}(\bm)\ll P^d\frac{(\log P)^{r_1+r_2}}{\N\fb_1}\cdot\frac{P^{d/2}}{\N\fb_2}(\log P)^{C}\cH(P\bm)^{-n/2}.\qedhere
\]
\end{proof}

\subsection{Dependence on the second modulus}
In Section~\ref{sec:error} we shall sum $I_{\fb_1,\fb_2}(\bm)$ over the second modulus with $\bm$ moving in a prescribed way. The relevant function is the following. Fix $\fb_1$, an integral ideal $\fc$, an element $\alpha_1\in K^\times$ and $\tv\in V^n$ with $|\tv^{(l)}|>0$ for all $l$. For $\beta\in V^\times$ put
\[
\rho_2(\beta)=\frac{|\Nm\beta|}{\N\fc\,P^{d/2}},\qquad\Phi(\beta)=I_{\fb_1}[\rho_2(\beta)]\Bigl(\frac{\tv}{\alpha_1\beta}\Bigr),
\]
where $I_{\fb_1}[\rho](\bm)$ denotes the right-hand side of~\eqref{eq:Ibig} with $\rho$ in place of $\N\fb_2/P^{d/2}$. When $\beta$ generates $\fb_2\fc$, we have $\rho_2(\beta)=\rho_2$, so that $\Phi(\beta)=I_{\fb_1,\fb_2}(\tv/(\alpha_1\beta))$. Recall the operators $\mathcal D_l$ from~\S\ref{sec:hecke}.

\begin{lemma}\label{lem:smooth}
Let $j\geqslant0$ and $D_1,\dots,D_j\in\bigcup_l\mathcal D_l$. Then $D_1\cdots D_j\Phi(\beta)$ is a linear combination, with coefficients bounded in terms of $j$ and $n$, of at most $5^j$ functions of the form
\[
\beta\mapsto\frac{P^{2d}}{C_P}\sum_{\substack{c_1\in\fb_1\setminus\{0\}\\ \ip{c_1}\leqslant A_0P}}\frac1{|\N c_1|}I^\natural_{c_1}[\rho_2(\beta)]\Bigl(\frac{\tv}{\alpha_1\beta}\Bigr),
\]
where each $I^\natural_{c_1}[\rho]$ is an integral of type $\mathfrak J_{2j}$ with $\rho$ in place of $\rho_2$, whose ingredients come from a finite set depending only on $j$. In particular,
\[
D_1\cdots D_j\Phi(\beta)\ll_jL^{2j}(\log P)^{C}\frac{P^{3d/2}}{\N\fb_1}\cdot\frac{\N\fc}{|\Nm\beta|}\cdot\cH\Bigl(\frac{P\tv}{\alpha_1\beta}\Bigr)^{-n/2}.
\]
\end{lemma}

\begin{proof}
The sum over $c_1$ is finite and does not involve $\beta$, so it suffices to treat a single $c_1$. Put $F(\beta)=I_{c_1}[\rho_2(\beta)](\tv/(\alpha_1\beta))$, where $I_{c_1}[\rho]$ denotes~\eqref{eq:Iint} with $\rho$ in place of $\rho_2$. More generally, let $F^\natural(\beta)=I^\natural[\rho_2(\beta)](\tv/(\alpha_1\beta))$, where $I^\natural$ is of type $\mathfrak J_D$. We show that each $D\in\bigcup_l\mathcal D_l$ maps $F^\natural$ to a linear combination of at most five functions of the same form, of type $\mathfrak J_{D+2}$. The lemma then follows by induction, and the final bound follows from Proposition~\ref{cor:Ibound} and Lemma~\ref{lem:latticesum}, since $\rho_2(\beta)^{-1}=P^{d/2}\N\fc/|\Nm\beta|$.

Substitute $\by=\beta\bz$, the product being taken componentwise in $V$. Then $d\by=|\Nm\beta|^nd\bz$, $Q_i(\by)=\beta^2Q_i(\bz)$, $\Nm(\bl Q_2(\by))=\Nm(\beta)^2\Nm(\bl Q_2(\bz))$, and
\[
\frac{P\tv}{\alpha_1\beta}\cdot\by=\frac{P\tv\cdot\bz}{\alpha_1}
\]
does not depend on $\beta$. Hence
\begin{multline*}
F^\natural(\beta)=|\Nm\beta|^n\int_{V^n}W^\natural(\beta\bz)\,A(\beta^2Q_1(\bz))\,k\bigl(\rho_2(\beta),\Nm(\beta)^2\Nm(\bl Q_2(\bz))\bigr)\\
\times u_V(\bl\beta^2Q_2(\bz))\,\ee\Bigl(-\frac{P\tv\cdot\bz}{\alpha_1}\Bigr)d\bz.
\end{multline*}
The integrand is smooth in $(\beta,\bz)$ and, locally uniformly in $\beta$, compactly supported in $\bz$. So we may differentiate under the integral sign.

Let $D=\beta^{(l)}\partial/\partial\beta^{(l)}$ if $l$ is real, or one of the Wirtinger operators $\beta^{(l)}\partial_{\beta^{(l)}}$, $\bar\beta^{(l)}\partial_{\bar\beta^{(l)}}$ if $l$ is complex. We have $D|\Nm\beta|=|\Nm\beta|$ and $D\Nm(\beta)^2=2\Nm(\beta)^2$. By the product rule, $DF^\natural$ is a sum of five terms, according to which factor is differentiated.
\begin{enumerate}[label=(\arabic*)]
\item $D|\Nm\beta|^n=n|\Nm\beta|^n$.
\item $D[W^\natural(\beta\bz)]=(\mathcal E_lW^\natural)(\beta\bz)$, where $\mathcal E_l=\sum_i y_i^{(l)}\partial_{y_i^{(l)}}$ (with Wirtinger derivatives at a complex place, holomorphic or antiholomorphic according to $D$). If $W^\natural=\Pi(\by-\bxi)w_0e^{-L\|\by-\bxi\|^2}+R$, then $\mathcal E_lW^\natural$ equals
\[
\bigl(\mathcal E_l\Pi-L\,\Pi\cdot\mathcal E_l\|\by-\bxi\|^2\bigr)(\by-\bxi)\,w_0e^{-L\|\by-\bxi\|^2}+\Pi\,(\mathcal E_lw_0)e^{-L\|\by-\bxi\|^2}+\mathcal E_lR.
\]
Writing $\by=(\by-\bxi)+\bxi$, the first bracket is a polynomial in $\by-\bxi$ of degree at most $D+2$ with coefficients $O(L^{D+1})$. The middle term is supported where $\delta_0/2\leqslant\ip{\by-\bxi}\leqslant\delta_0$, where $e^{-L\|\by-\bxi\|^2}\leqslant e^{-L\delta_0^2/4}$; together with $\mathcal E_lR$ it satisfies the conditions on the remainder in Definition~\ref{def:class}(d) with $D+2$ in place of $D$. So $\mathcal E_lW^\natural$ is a Gaussian weight of degree $D+2$.
\item If $A(v)=\kappa_0-\prod_{l'}a_{l'}(|Pv^{(l')}/c_1^{(l')}|)$, then $D[A(\beta^2Q_1(\bz))]=-\prod_{l'}\tilde a_{l'}(\cdots)$, with $\tilde a_{l'}=a_{l'}$ for $l'\neq l$ and $\tilde a_l(x)=c\,xa_l'(x)$, $c\in\{1,2\}$. Indeed $D|\beta^{(l)2}v|=c|\beta^{(l)2}v|$. This is an element weight with $\kappa_0=0$.
\item $D[k(\rho_2(\beta),\Nm(\beta)^2Y)]=(x\partial_xk+2y\partial_yk)(\rho_2(\beta),\Nm(\beta)^2Y)$, which is a shell weight by Lemma~\ref{lem:classbasic}(iv).
\item $D[u_V(\bl\beta^2Q_2(\bz))]$ is obtained by replacing $u_l(x)$ by $c\,xu_l'(x)$, which is again a product cut-off.
\end{enumerate}
In each case, changing variables back to $\by=\beta\bz$ shows that the term is $F^{\natural\prime}(\beta)$ for an integral $I^{\natural\prime}$ of type $\mathfrak J_{D+2}$, multiplied by a constant. This completes the induction.
\end{proof}

\subsection{Archimedean truncation}\label{sec:truncation}
We come to the main new point of the archimedean analysis. The dual variable $\bm$ is negligible as soon as $\cH(P\bm)$ exceeds the expected size, $P^{3d/2}/\N(\fb_1\fb_2)$, \emph{for every number field $K$}.

\begin{theorem}[Archimedean truncation]\label{thm:truncation}
Let $\varepsilon>0$, and let $I^\natural$ be an integral of type $\mathfrak J_D$. Then $I^\natural(\bm)$ is negligible whenever
\[
\cH(P\bm)\geqslant\frac{P^{3d/2+\varepsilon}}{\N(\fb_1\fb_2)}.
\]
In particular this holds for $I_{c_1,\fb_2}(\bm)$, and for $I_{\fb_1,\fb_2}(\bm)$.
\end{theorem}

\begin{proof}
By Lemma~\ref{lem:ibp} we may assume that $\ip\bm\leqslant P^{2C_0}$. Put $\varepsilon'=\varepsilon/(2d+2)$. As in the proof of Proposition~\ref{cor:Ibound}, apart from a negligible error we may restrict the integral~\eqref{eq:fourier} to the range
\begin{equation}\label{eq:trange}
|t_1^{(l)}|\leqslant\frac{P^{1+\varepsilon'}}{|c_1^{(l)}|}\quad(1\leqslant l\leqslant r_1+r_2),\qquad\cH(t_2)\leqslant\rho_2^{-1}P^{\varepsilon'}.
\end{equation}
In this range all $|t_i^{(l)}|$ are $P^{O(1)}$, and we use Lemma~\ref{lem:Kfactor}. Its error term is negligible after integration against $|q(t_1)p_{\rho_2}(t_2)|$. In the main term, Lemma~\ref{lem:gaussian} shows that the product $\prod_lJ_{l,\alpha_l}(\bs^{(l)};P\bm^{(l)})$, where $\bs=(t_1,\bl t_2)$, is $O(e^{-c_3(\log P)^2})$ unless~\eqref{eq:localise} holds at every place, with $\bu=P\bm^{(l)}$. If it does, then
\[
|P\bm^{(l)}|\leqslant|T\bs^{(l)}|+|\bw_l|\ll|t_1^{(l)}|+\bl^{(l)}|t_2^{(l)}|+L\qquad\text{for every }l.
\]
Now we use the archimedean sizes of $\bl$. By~\eqref{eq:trange} and~\eqref{eq:lambdabounds}, $|t_1^{(l)}|\leqslant P^{\varepsilon'}P/|c_1^{(l)}|\leqslant P^{\varepsilon'}\bl^{(l)}$; and $\bl^{(l)}\geqslant A_0^{-1}$. Hence
\[
\max\bigl(1,|P\bm^{(l)}|\bigr)\ll P^{\varepsilon'}L\,\bl^{(l)}\max\bigl(1,|t_2^{(l)}|\bigr)\qquad\text{for every }l.
\]
This is the key inequality: at \emph{each} place separately, the first dual variable is dominated by the second. Raising to the power $c_l$ and multiplying over the places, we find
\[
\cH(P\bm)\ll P^{d\varepsilon'}L^d\,\Nm(\bl)\,\cH(t_2)\leqslant P^{(d+1)\varepsilon'}L^d\cdot\frac{P^d}{\N\fb_1}\cdot\frac{P^{d/2}}{\N\fb_2}\leqslant P^{\varepsilon/2}L^d\,\frac{P^{3d/2}}{\N(\fb_1\fb_2)}.
\]
When $P$ is large, this contradicts the hypothesis. Therefore the main term is $O(e^{-c_3(\log P)^2})$ throughout the range~\eqref{eq:trange}. Integrating against $|q(t_1)p_{\rho_2}(t_2)|$, whose integral is $\ll\rho_2^{-1}(\log P)^C\ll P^{d/2}$, gives a negligible result. For $I_{\fb_1,\fb_2}(\bm)$, note that the hypothesis does not involve $c_1$ and that the sum over $c_1$ in~\eqref{eq:Ibig} has $P^{O(1)}$ terms.
\end{proof}

\begin{remark}\label{rem:truncation}
The proof uses only two features of the weights. First, the Fourier transform of the first weight is confined, place by place, to $|t_1^{(l)}|\lesssim P/|c_1^{(l)}|$. Second, $\bl^{(l)}\geqslant P/|c_1^{(l)}|$ at every place, while $\Nm(\bl)=P^d/\N\fb_1$ only depends on the modulus. If the first equation were detected with the weight $h(\rho_1,|\Nm Q_1(\by)|)$ of Theorem~\ref{thm:BV}, Lemma~\ref{lem:p2} would only give $\cH(t_1)\lesssim\rho_1^{-1}$. The first inequality above would then fail when $t_1$ is concentrated at a single place, and for a scalar $\bl$ one would only obtain $\cH(P\bm)\lesssim\rho_1^{-1}\rho_2^{-1}\bl^{\,d-c_{\min}}$; see Section~\ref{sec:higher}.
\end{remark}
\section{The error term}\label{sec:error}

Let $E(P)$ denote the contribution of the terms with $\bm\neq\mathbf0$ in Proposition~\ref{prop:poisson}:
\begin{equation}\label{eq:EP}
E(P)=\Bigl(\frac{2^{r_2}}{\sqrt{D_K}}\Bigr)^nP^{d(n-3)}\sum_{\fb_1,\fb_2}\frac1{\N(\fb_1\fb_2)^n}\sum_{\substack{\bm\in(\fb_1\fb_2)^{\wedge n}\\ \bm\neq\mathbf0}}C_{\fb_1,\fb_2}(\bm)\,I_{\fb_1,\fb_2}(\bm).
\end{equation}
Here $\fb_1$ runs over integral ideals with $\N\fb_1\leqslant(A_0P)^d$, and $\fb_2$ over all non-zero integral ideals. Indeed, the integrand in~\eqref{eq:Iint} vanishes when $\rho_2>c_{\ast}$, because $|\Nm(\bl Q_2(\by))|\leqslant2^d$ on the support of $U_V(\bl Q_2(\by))$ and $h(x,y)=0$ for $x>\max(1,2|y|)$. Only finitely many $\fb_2$ occur, and for each pair $(\fb_1,\fb_2)$ the sum over $\bm$ converges absolutely by Lemma~\ref{lem:ibp}. So all rearrangements below are justified. The aim of this section is the following estimate.

\begin{proposition}\label{prop:error}
Let $K$ be any number field and let $\Phi(n)$ be as in~\eqref{eq:theta}. Then, for every $\varepsilon>0$,
\[
E(P)\ll_\varepsilon P^{d\Phi(n)+\varepsilon}.
\]
\end{proposition}

The proof follows the structure of~\cite[\S5]{Mun15}. Over $\QQ$ the dual variable is an integer vector and the second modulus is an integer, so one can fix the vector and sum over the modulus. Over $K$ the dual lattice $(\fb_1\fb_2)^\wedge$ moves with $\fb_2$, and we first have to choose coordinates in which the dual variable is fixed while $\fb_2$ varies.

\subsection{Two auxiliary lemmas}
\begin{lemma}\label{lem:unitsum}
Let $\kappa>0$ and let $\bw\in V^n$ with $|\bw^{(l)}|>0$ for all $l$. Then
\[
\sum_{u\in U_K}\cH\bigl(u^{-1}\bw\bigr)^{-\kappa}\ll_\kappa\bigl(1+|\log\cN(\bw)|\bigr)^{r_1+r_2-1}\min\bigl(1,\cN(\bw)^{-\kappa}\bigr).
\]
\end{lemma}

\begin{proof}
Put $r=r_1+r_2$, $a_l=c_l\log|\bw^{(l)}|$ and $e_l(u)=c_l\log|u^{(l)}|$, so that $\sum_le_l(u)=0$ and $\sum_la_l=A:=\log\cN(\bw)$. With $x_l=a_l-e_l(u)$ we have
\[
\log\cH(u^{-1}\bw)=\sum_l\max(0,x_l)=\tfrac12\bigl(\|\mathbf x\|_1+A\bigr),
\]
since $\max(0,x)=\frac12(|x|+x)$ and $\sum_lx_l=A$. Moreover $\|\mathbf x\|_1\geqslant|A|$. By Dirichlet's unit theorem the vectors $(e_l(u))_l$ form a lattice of rank $r-1$ in the hyperplane $\{\sum_le_l=0\}$, and each value is taken by exactly $w_K$ units. Hence the number of $u$ with $\|\mathbf x\|_1\leqslant k$ is $O((1+k)^{r-1})$, uniformly in $\mathbf a$. Summing over $k\geqslant|A|$ in unit steps,
\[
\sum_u\cH(u^{-1}\bw)^{-\kappa}\ll\sum_{k\geqslant|A|}(1+k)^{r-1}e^{-\kappa(k+A)/2}.
\]
If $A\geqslant0$ this is $\ll(1+A)^{r-1}e^{-\kappa A}$, and if $A<0$ it is $\ll(1+|A|)^{r-1}$.
\end{proof}

\begin{lemma}\label{lem:smallprime}
Let $\mathcal C$ be an ideal class and let $\fa$ be a non-zero integral ideal with $\N\fa\leqslant P^{O(1)}$. Then there is a prime ideal $\fq\in\mathcal C$ with $\fq\nmid2D_K\Delta_X\fa$ and $\N\fq\ll(\log P)^2$.
\end{lemma}

\begin{proof}
By the prime ideal theorem for the classes of the ideal class group (equivalently, the Chebotarev density theorem for the Hilbert class field of $K$), the number of prime ideals in $\mathcal C$ of norm at most $x$ is $\gg x/\log x$ for $x\geqslant x_0(K)$. The ideal $2D_K\Delta_X\fa$ has $O(\log P)$ prime divisors. Taking $x=C(\log P)^2$ with $C$ large, there are more primes in $\mathcal C$ of norm at most $x$ than prime divisors of $2D_K\Delta_X\fa$.
\end{proof}

\subsection{Normalising the dual variable}
Fix, once and for all, integral ideals $\fc$ representing the ideal classes of $K$. For each $\fb_1$ and each $\fc$, Lemma~\ref{lem:smallprime} provides a prime $\fq=\fq(\fb_1,\fc)$ in the class of $(\fc^{-1}\fb_1\fd)^{-1}$ with $\fq\nmid2D_K\Delta_X\fb_1$ and $\N\fq\ll(\log P)^2$. Then $\fc^{-1}\fb_1\fd\fq$ is principal, and we fix a generator $\alpha_1=\alpha_1(\fb_1,\fc)$:
\begin{equation}\label{eq:alpha1}
(\alpha_1)=\fc^{-1}\fb_1\fd\fq,\qquad|\Nm\alpha_1|=\frac{D_K\N\fb_1\N\fq}{\N\fc}.
\end{equation}
Let $\mathcal F\subset V^n$ be the set of $\bv$ with $|\bv^{(l)}|>0$ for all $l$ such that the vector $(c_l\log|\bv^{(l)}|-\frac{c_l}d\log\cN(\bv))_l$ lies in a fixed bounded fundamental domain for the lattice of logarithms of units. Every element of $\mathcal F$ is balanced, and every orbit of $U_K$ acting on $\{\bv:|\bv^{(l)}|>0\ \forall l\}$ meets $\mathcal F$ in exactly $w_K$ points.

\begin{lemma}\label{lem:param}
Fix $\fb_1$ and $\fc$, and let $\fq$ and $\alpha_1$ be as above. Consider the pairs $(\fb_2,\bm)$, where $\fb_2$ is an integral ideal in the class of $\fc^{-1}$ and $\bm\in(\fb_1\fb_2)^{\wedge n}\setminus\{\mathbf0\}$. Then
\[
(\beta,\tv)\longmapsto\bigl((\beta)\fc^{-1},\ \tv/(\alpha_1\beta)\bigr)
\]
induces a bijection from $U_K\backslash\{(\beta,\tv)\in(\fc\setminus\{0\})\times(\fq^n\setminus\{\mathbf0\})\}$, where $u$ acts by $(\beta,\tv)\mapsto(u\beta,u\tv)$, onto the set of such pairs. Moreover:
\begin{enumerate}[label=\textup{(\roman*)}]
\item for every prime $\fp\mid\fb_1\fb_2$ with $\fp\neq\fq$, the normalised frequency of $\bm$ at $\fp$, at level $\ord_{\fp}(\fb_1\fb_2)$, is a unit multiple of $\tv$;
\item $\cN(P\bm)=P^d\cN(\tv)/(D_K\N\fq\N(\fb_1\fb_2))$.
\end{enumerate}
Consequently, for any function $\Psi$ for which the sums converge absolutely,
\[
\sum_{\fb_2\in[\fc^{-1}]}\ \sum_{\bm\neq\mathbf0}\Psi(\fb_2,\bm)=\frac1{w_K}\sum_{\tv\in\fq^n\cap\mathcal F}\ \sum_{\beta\in\fc\setminus\{0\}}\Psi\Bigl((\beta)\fc^{-1},\frac{\tv}{\alpha_1\beta}\Bigr).
\]
\end{lemma}

\begin{proof}
If $\beta\in\fc\setminus\{0\}$, then $(\beta)\fc^{-1}$ is integral and lies in the class of $\fc^{-1}$. Conversely, if $\fb_2$ is integral in this class, then $\fb_2\fc=(\beta)$ for some $\beta\in\fc$, unique up to units. By~\eqref{eq:alpha1},
\[
(\alpha_1\beta)(\fb_1\fb_2)^{-1}\fd^{-1}=\fc^{-1}\fb_1\fd\fq\cdot\fb_2\fc\cdot\fb_1^{-1}\fb_2^{-1}\fd^{-1}=\fq.
\]
Hence $\bm\mapsto\alpha_1\beta\bm$ is a bijection from $((\fb_1\fb_2)^{-1}\fd^{-1})^n$ onto $\fq^n$. Replacing $\beta$ by $u\beta$ replaces $\tv=\alpha_1\beta\bm$ by $u\tv$. This proves the bijection. For~(i), let $\fp\mid\fb_1\fb_2$ with $\fp\neq\fq$, and put $k=\ord_{\fp}(\fb_1\fb_2)$. Then $\ord_{\fp}(\alpha_1\beta)=k+\ord_{\fp}\fd$, so that $\alpha_1\beta=\eta\varpi^k\delta_{\fp}$ with $\eta\in\cO_{\fp}^\times$, and $\bm=\eta^{-1}\tv/(\varpi^k\delta_{\fp})$. Part~(ii) follows from $|\Nm(\alpha_1\beta)|=|\Nm\alpha_1|\N\fb_2\N\fc$ and~\eqref{eq:alpha1}. For the last assertion, choose in each orbit the $w_K$ representatives with $\tv\in\mathcal F$; this is possible since every $\tv\in\fo^n\setminus\{\mathbf0\}$ has $|\tv^{(l)}|>0$ for all $l$.
\end{proof}

Applying Lemma~\ref{lem:param} for each $\fb_1$ and each class representative $\fc$, we obtain
\begin{equation}\label{eq:Esplit}
E(P)=\frac1{w_K}\Bigl(\frac{2^{r_2}}{\sqrt{D_K}}\Bigr)^nP^{d(n-3)}\sum_{\fb_1}\frac1{\N\fb_1^{\,n}}\sum_{\fc}\ \sum_{\tv\in\fq^n\cap\mathcal F}\Sigma(\fb_1,\fc,\tv),
\end{equation}
where, writing $\fb_2=(\beta)\fc^{-1}$,
\[
\Sigma(\fb_1,\fc,\tv)=\sum_{\beta\in\fc\setminus\{0\}}\frac{C_{\fb_1,\fb_2}(\tv/(\alpha_1\beta))}{\N\fb_2^{\,n}}\,I_{\fb_1,\fb_2}\Bigl(\frac{\tv}{\alpha_1\beta}\Bigr).
\]

\subsection{Truncation of the dual variable}
Put $M_0=P^{d/2+\varepsilon}$. Recall that elements of $\mathcal F$ are balanced, so that $\cH(\tv)\asymp\cN(\tv)$ for $\tv\in\mathcal F\cap\fo^n$ by Lemma~\ref{lem:balance}.

\begin{lemma}\label{lem:tails}
The contribution to~\eqref{eq:Esplit} of the vectors $\tv$ with $\cH(\tv)>M_0$ is negligible.
\end{lemma}

\begin{proof}
Let $\cH(\tv)>M_0$, let $\beta\in\fc\setminus\{0\}$, and put $\bm=\tv/(\alpha_1\beta)$ and $\fb_2=(\beta)\fc^{-1}$. By Lemma~\ref{lem:param}(ii),
\[
\cH(P\bm)\geqslant\cN(P\bm)\gg\frac{P^d\cH(\tv)}{(\log P)^2\N(\fb_1\fb_2)}\geqslant\frac{P^{3d/2+\varepsilon/2}}{\N(\fb_1\fb_2)}
\]
for $P$ large. So $|I_{\fb_1,\fb_2}(\bm)|\leqslant P^{-A}$ by Theorem~\ref{thm:truncation}, for any fixed $A$. We also have $|I_{\fb_1,\fb_2}(\bm)|\ll P^{O(1)}\cH(P\bm)^{-n/2}$ by Proposition~\ref{cor:Ibound}, and $|I_{\fb_1,\fb_2}(\bm)|\ll_NP^{O(N)}\ip{P\bm}^{-N}$ by Lemma~\ref{lem:ibp}, where $\ip{P\bm}^d\geqslant\cN(P\bm)\gg P^{d-O(1)}\cH(\tv)$. Taking the geometric mean of the three bounds, with $N=3dn$,
\[
|I_{\fb_1,\fb_2}(\bm)|\ll P^{-A/3+O(1)}\cH(P\bm)^{-n/6}\cH(\tv)^{-n}.
\]
Lemma~\ref{lem:unitsum} with $\kappa=n/6$, applied to $\bw=P\tv/(\alpha_1\beta)$, shows that the sum of $\cH(P\bm)^{-n/6}$ over the $\beta$ in a fixed $U_K$-orbit is $O((\log P)^{r_1+r_2})$, since $|\log\cN(P\bm)|\ll\log P$. The number of orbits with $\N\fb_2\leqslant c_{\ast}P^{d/2}$ is $O(P^{d/2})$, and $|C_{\fb_1,\fb_2}(\bm)|\leqslant P^{O(1)}$ by Corollary~\ref{cor:pointwise}. Finally, $\sum_{\tv}\cH(\tv)^{-n}$ over balanced $\tv\in\fo^n\setminus\{\mathbf0\}$ converges, by Lemma~\ref{lem:count-trivial}(iii) and partial summation. Since the sums over $\fb_1$ and $\fc$ have $P^{O(1)}$ terms, the total is $O(P^{-A/3+O(1)})$, and $A$ is arbitrary.
\end{proof}

\subsection{Splitting the second modulus}
From now on we fix $\fb_1$, $\fc$ and $\tv\in\fq^n\cap\mathcal F$ with $\cH(\tv)\leqslant M_0$. Let $\cS=\cS(\fb_1,\fc,\tv)$ be the set of prime ideals dividing
\[
2D_K\Delta_X\,\fc\,\fq\,\fb_1\cdot Q_2^\ast(\tv),
\]
where the factor $Q_2^\ast(\tv)$ is omitted if $Q_2^\ast(\tv)=0$. All of these ideals have norm $P^{O(1)}$. For $\beta\in\fc\setminus\{0\}$ we write $\fb_2=(\beta)\fc^{-1}=\fd_2\fb_2'$, where $\fd_2$ is supported on $\cS$ and $\fb_2'$ is coprime to $\cS$. Put $\fc'=\fc\fd_2$, so that $\beta\in\fc'$ and $\fb_2'=(\beta)\fc'^{-1}$.

\begin{lemma}\label{lem:factor}
In this situation
\[
C_{\fb_1,\fb_2}\Bigl(\frac{\tv}{\alpha_1\beta}\Bigr)=C^{\cS}_{\fb_1,\fd_2}(\tv)\cdot C^\flat(\fb_2'),
\]
where $C^{\cS}_{\fb_1,\fd_2}(\tv)$ depends on $\fb_1$, $\fd_2$ and $\tv$ but not on $\beta$, and where $C^\flat(\fb_2')=\prod_{\fp^\ell\|\fb_2'}\mathcal Q_{\fp}(\ell;\tv)$ satisfies the following.
\begin{enumerate}[label=\textup{(\roman*)}]
\item If $n$ is odd and $Q_2^\ast(\tv)\neq0$, then $C^\flat(\fb_2')=\mu^2(\fb_2')\N\fb_2'^{\,(n+1)/2}\chi_{\tv}(\fb_2')$, with $\chi_{\tv}$ as in Remark~\ref{rem:hecke}.
\item If $n$ is odd and $Q_2^\ast(\tv)=0$, then $C^\flat(\fb_2')=0$ unless $\fb_2'$ is square-full, and $|C^\flat(\fb_2')|\leqslant\N\fb_2'^{\,n/2+1}$.
\item If $n$ is even and $Q_2^\ast(\tv)\neq0$, then $C^\flat(\fb_2')=0$ unless $\fb_2'$ is square-free, and $|C^\flat(\fb_2')|=\N\fb_2'^{\,n/2}$.
\item If $n$ is even and $Q_2^\ast(\tv)=0$, then $|C^\flat(\fb_2')|\leqslant\N\fb_2'^{\,n/2+1}$.
\end{enumerate}
\end{lemma}

\begin{proof}
By Lemma~\ref{lem:mult}, $C_{\fb_1,\fb_2}(\bm)$ is the product of the local sums $C_{\fp}(r_{\fp},\ell_{\fp};\bm)$ over $\fp\mid\fb_1\fb_2$. For $\fp\in\cS$, the normalised frequency of $\bm$ at $\fp$ is a unit multiple of $\tv$, or of $\tv/\varpi_{\fq}$ if $\fp=\fq$ (Lemma~\ref{lem:param}(i), and the same computation at $\fq$). The unit depends on $\beta$, but by Lemma~\ref{lem:unitinv} the local sum does not; and the exponents $r_{\fp},\ell_{\fp}$ only depend on $\fb_1$ and $\fd_2$. This defines $C^{\cS}_{\fb_1,\fd_2}(\tv)$.

Now let $\fp^\ell\|\fb_2'$. Then $\fp\notin\cS$, so $\fp$ is good, $\fp\nmid\fb_1$ (hence $r_{\fp}=0$), $\fp\neq\fq$, and $\fp\nmid Q_2^\ast(\tv)$ if $Q_2^\ast(\tv)\neq0$. The local factor is $C_{\fp}(0,\ell;\bm)=\mathcal Q_{\fp}(\ell;\bm)$, and by Lemma~\ref{lem:unitinv} we may take the normalised frequency to be $\tv$. If $Q_2^\ast(\tv)\neq0$, then $\mathcal Q_{\fp}(\ell;\tv)=0$ for $\ell\geqslant2$ by Lemma~\ref{lem:gauss}(iii), while for $\ell=1$ Lemma~\ref{lem:gauss}(i), (ii) give $q^{(n+1)/2}\chi_{\fp}(\epsilon_nQ_2^\ast(\tv))$ for odd $n$ and a number of absolute value $q^{n/2}$ for even $n$. For odd $n$, $\chi_{\fp}(\epsilon_nQ_2^\ast(\tv))=\chi_{\tv}(\fp)$ since $\fp\nmid2Q_2^\ast(\tv)$. This proves~(i) and~(iii). If $Q_2^\ast(\tv)=0$ and $n$ is odd, then $\mathcal Q_{\fp}(1;\tv)=q^{(n+1)/2}\chi_{\fp}(0)=0$, and $|\mathcal Q_{\fp}(\ell;\tv)|\leqslant q^{\ell(n/2+1)}$ by Lemma~\ref{lem:gauss}(iv); this proves~(ii). If $Q_2^\ast(\tv)=0$ and $n$ is even, Lemma~\ref{lem:gauss}(ii) gives $|\mathcal Q_{\fp}(1;\tv)|=q^{n/2}(q-1)$, and~(iv) follows with Lemma~\ref{lem:gauss}(iv).
\end{proof}

We attach to $\tv$ two exponents $\theta(\tv)$ and $\theta_2(\tv)$, as follows. Call $\tv$ \emph{exceptional} if $n$ is odd, $Q_2^\ast(\tv)\neq0$ and $\epsilon_nQ_2^\ast(\tv)\in\mathcal A_K$ (Remark~\ref{rem:unramified}); equivalently, $\chi_{\tv}\psi=1$ for some character $\psi$ of the class group. Put
\[
(\theta(\tv),\theta_2(\tv))=\begin{cases}
(\theta_K,0)&\text{$n$ odd, $Q_2^\ast(\tv)\neq0$, $\tv$ not exceptional},\\
(0,\frac12)&\text{$n$ odd, and $\tv$ exceptional or $Q_2^\ast(\tv)=0$},\\
(0,0)&\text{$n$ even, $Q_2^\ast(\tv)\neq0$},\\
(0,1)&\text{$n$ even, $Q_2^\ast(\tv)=0$}.
\end{cases}
\]

\subsection{The sum over the second modulus}
\begin{lemma}\label{lem:betasum}
With $\fb_1$, $\fc$, $\tv$ as above and $\fd_2$ supported on $\cS$, put
\[
\Sigma_{\fd_2}=\sum_{\substack{\beta\in\fc'\setminus\{0\}\\ ((\beta)\fc'^{-1},\cS)=1}}\frac{C^\flat\bigl((\beta)\fc'^{-1}\bigr)}{\N\fb_2^{\,n}}\,I_{\fb_1,\fb_2}\Bigl(\frac{\tv}{\alpha_1\beta}\Bigr),\qquad\fb_2=(\beta)\fc^{-1}.
\]
Then
\[
\Sigma_{\fd_2}\ll P^\varepsilon\,\frac1{\N\fd_2^{\,n}}\cdot\frac{P^{3d/2}}{\N\fb_1\N\fd_2}\Bigl(\frac{\N\fb_1\N\fd_2}{P^d\cH(\tv)}\Bigr)^{n/2}\cH(\tv)^{\theta(\tv)}\Bigl(\frac{P^{d/2}}{\N\fd_2}\Bigr)^{\theta_2(\tv)}.
\]
\end{lemma}

\begin{proof}
Fix $\Omega\in C_c^\infty([\frac12,2])$ with $\sum_{k\in\ZZ}\Omega(x/2^k)=1$ for $x>0$. For $X=2^k$ and $\beta\in V^\times$ put
\[
F_X(\beta)=\Bigl(\frac{|\Nm\beta|}{\N\fc}\Bigr)^{-n}\Phi(\beta)\,\Omega\Bigl(\frac{|\Nm\beta|}{\N\fc'\,X}\Bigr),
\]
with $\Phi$ as in Lemma~\ref{lem:smooth}. Then $\Sigma_{\fd_2}=\sum_X\Sigma_{\fd_2,X}$, where $\Sigma_{\fd_2,X}$ is the same sum with $F_X(\beta)$ in place of $\N\fb_2^{-n}I_{\fb_1,\fb_2}(\tv/(\alpha_1\beta))$. On the support of $F_X$ we have $\N\fb_2'\asymp X$ and $\N\fb_2\asymp\N\fd_2X$. Since $\Phi(\beta)=0$ when $\rho_2(\beta)>c_{\ast}$, and $\N\fb_2'\geqslant1$, only $O(\log P)$ values of $X$ occur, all with $X\ll P^{d/2}/\N\fd_2$.

We bound $F_X$ and its derivatives on unit orbits. By Lemma~\ref{lem:smooth}, for $D_1,\dots,D_j\in\bigcup_l\mathcal D_l$,
\[
D_1\cdots D_jF_X(\beta)\ll_jL^{2j}(\log P)^C(\N\fd_2X)^{-n}\frac{P^{3d/2}}{\N\fb_1\N\fd_2X}\cH\Bigl(\frac{P\tv}{\alpha_1\beta}\Bigr)^{-n/2},
\]
the derivatives of the factors $|\Nm\beta|^{-n}$ and $\Omega(\cdots)$ being harmless since these are functions of $|\Nm\beta|$ with bounded logarithmic derivatives. The right-hand side depends on $\beta$ only through $\cH(P\tv/(\alpha_1\beta))$, and the operators $D_i$ commute with multiplication by units. By Lemma~\ref{lem:unitsum} with $\kappa=n/2$ and Lemma~\ref{lem:param}(ii),
\[
\sum_{u\in U_K}\cH\Bigl(\frac{P\tv}{\alpha_1u\beta}\Bigr)^{-n/2}\ll(\log P)^{r_1+r_2}\cN\Bigl(\frac{P\tv}{\alpha_1\beta}\Bigr)^{-n/2}\ll(\log P)^{C}\Bigl(\frac{\N\fb_1\N\fd_2X}{P^d\cH(\tv)}\Bigr)^{n/2},
\]
where we used $\min(1,x^{-n/2})\leqslant x^{-n/2}$ and $\N\fq\ll(\log P)^2$. Hence $F_X^U(\beta)=\sum_uF_X(u\beta)$ converges absolutely, and for every $j$ we have $\sup_\beta|D_1\cdots D_jF_X^U(\beta)|\leqslant P^{j\varepsilon}B$ and $\sup_\beta\sum_u|F_X(u\beta)|\leqslant B$, where
\[
B=P^{\varepsilon/2}(\N\fd_2X)^{-n}\frac{P^{3d/2}}{\N\fb_1\N\fd_2X}\Bigl(\frac{\N\fb_1\N\fd_2X}{P^d\cH(\tv)}\Bigr)^{n/2}.
\]
Note that $X^{n/2+1}B$ does not depend on $X$.

Suppose first that $n$ is odd and $Q_2^\ast(\tv)\neq0$. By Lemma~\ref{lem:factor}(i), $\Sigma_{\fd_2,X}$ is exactly the sum of Proposition~\ref{prop:hecke}, with $\fc'$, $\cS$, $\chi=\chi_{\tv}$ and $F=F_X$. The hypotheses are satisfied: $\cS$ contains the primes dividing $2\fc'$, the conductor $\ff_{\tv}$ divides $4Q_2^\ast(\tv)$ and so is supported on $\cS$, and all norms involved are $P^{O(1)}$. If $\tv$ is not exceptional, the first bound of Proposition~\ref{prop:hecke} gives
\[
\Sigma_{\fd_2,X}\ll P^\varepsilon BX^{\frac{n+2}2}(\N\ff_{\tv})^{\frac14-\delta_{\mathrm{MV}}}.
\]
Since $\tv$ is balanced, $\N\ff_{\tv}\leqslant\N(4Q_2^\ast(\tv))\ll\cN(\tv)^2\ll\cH(\tv)^2$, so that $(\N\ff_{\tv})^{1/4-\delta_{\mathrm{MV}}}\ll\cH(\tv)^{\theta_K}$. If $\tv$ is exceptional, the second bound gives $\Sigma_{\fd_2,X}\ll P^\varepsilon BX^{(n+3)/2}$. In the remaining cases we use Lemma~\ref{lem:betatrivial} together with Lemma~\ref{lem:factor}(ii)--(iv). This gives $\Sigma_{\fd_2,X}\ll BX^{n/2+3/2}$ if $n$ is odd and $Q_2^\ast(\tv)=0$; $\Sigma_{\fd_2,X}\ll BX^{n/2+1}$ if $n$ is even and $Q_2^\ast(\tv)\neq0$; and $\Sigma_{\fd_2,X}\ll BX^{n/2+2}$ if $n$ is even and $Q_2^\ast(\tv)=0$.

In all cases $\Sigma_{\fd_2,X}\ll P^\varepsilon X^{n/2+1}B\cdot\cH(\tv)^{\theta(\tv)}X^{\theta_2(\tv)}$. Summing over the $O(\log P)$ values of $X\ll P^{d/2}/\N\fd_2$ gives the lemma.
\end{proof}

\subsection{The sum over the moduli}
We now sum over $\fb_1$ and $\fd_2$. In addition to $\theta_2(\tv)$, put $\psi_1(\tv)=0$ if $G(\tv)\neq0$ and either $Q_2^\ast(\tv)\neq0$ or $\theta_2(\tv)=1$, and put $\psi_1(\tv)=\frac12$ otherwise.

\begin{lemma}\label{lem:assembly}
Let $\tv\in\fo^n\setminus\{\mathbf0\}$ with $\cH(\tv)\leqslant P^{O(1)}$, and fix $\fc$. Then
\[
\sum_{\N\fb_1\leqslant(A_0P)^d}\ \sum_{\substack{\fd_2\ \text{supported on}\ \cS\\ \N\fd_2\leqslant c_{\ast}P^{d/2}}}\frac{|C^{\cS}_{\fb_1,\fd_2}(\tv)|}{\N(\fb_1\fd_2)^{n/2+1}}\Bigl(\frac{P^{d/2}}{\N\fd_2}\Bigr)^{\theta_2(\tv)}\ll P^{\varepsilon+d\theta_2(\tv)/2+d\psi_1(\tv)},
\]
where $\cS=\cS(\fb_1,\fc,\tv)$ and the sum over $\fb_1$ is restricted to those $\fb_1$ for which $\tv\in\fq(\fb_1,\fc)^n$.
\end{lemma}

\begin{proof}
Let $\fg_0$ be the product of the prime ideals dividing $2D_K\Delta_X\fc\cdot G(\tv)Q_2^\ast(\tv)$, where the factors $G(\tv)$ and $Q_2^\ast(\tv)$ are omitted if they vanish. Then $\N\fg_0\leqslant P^{O(1)}$. Write $\fb_1=\fd_1\fb_1'$, where $\fd_1$ is supported on the primes of $\fg_0$ and $\fb_1'$ is coprime to $\fg_0$. Every prime of $\cS$ divides $\fg_0\fq\fb_1'$, and $\fq\nmid\fb_1$. Write accordingly $\fd_2=\fd_2^0\fd_2^1$, where $\fd_2^0$ is supported on the primes of $\fg_0\fq$ and $\fd_2^1$ on those of $\fb_1'$. By Lemma~\ref{lem:mult},
\[
C^{\cS}_{\fb_1,\fd_2}(\tv)=\prod_{\fp\mid\fg_0\fq}C_{\fp}(r_{\fp},\ell_{\fp})\cdot\prod_{\fp\mid\fb_1'}C_{\fp}(r_{\fp},\ell_{\fp}),
\]
where we suppress the frequency, whose normalisation is a unit multiple of $\tv$ (or of $\tv/\varpi_{\fq}$ at $\fq$). We also use $(P^{d/2}/\N\fd_2)^{\theta_2}\leqslant P^{d\theta_2/2}\N(\fd_2^1)^{-\theta_2}$.

\emph{The primes of $\fg_0\fq$.} By Lemma~\ref{lem:universal}, the first product is at most $\prod_{\fp\mid\fg_0\fq}c_{\fp}\cdot\N(\fd_1\fd_2^0)^{n/2+1}$. The constants $c_{\fp}$ equal $4n$ except at the finitely many bad primes, and $(4n)^{\omega(\fg_0\fq)}\ll P^\varepsilon$ by Lemma~\ref{lem:smoothsupp}. By the same lemma, the number of pairs $(\fd_1,\fd_2^0)$ supported on the primes of $\fg_0\fq$ with norms at most $P^{O(1)}$ is $O(P^\varepsilon)$, uniformly in $\fq$ since $\N\fq\ll(\log P)^2$. So these primes contribute a factor $O(P^{2\varepsilon})$ in total.

\emph{The primes of $\fb_1'$.} Let $\fp^r\|\fb_1'$, so $r\geqslant1$. Then $\fp$ is good, $\fp\neq\fq$, $\fp\nmid G(\tv)$ if $G(\tv)\neq0$, and $\fp\nmid Q_2^\ast(\tv)$ if $Q_2^\ast(\tv)\neq0$. The contribution of $\fp$, summed over the possible exponents $\ell\geqslant0$ of $\fp$ in $\fd_2^1$, is
\[
f_{\fp}(r)=\sum_{\ell\geqslant0}\frac{|C_{\fp}(r,\ell)|}{q^{(r+\ell)(n/2+1)}}q^{-\ell\theta_2}.
\]
For $\ell\geqslant1$ we have $C_{\fp}(r,\ell)=0$ if $Q_2^\ast(\tv)\neq0$, by Lemma~\ref{lem:mixed}. If $Q_2^\ast(\tv)=0$, then $\theta_2\geqslant\frac12$, and Lemma~\ref{lem:universal} bounds the terms with $\ell\geqslant1$ by $\sum_{\ell\geqslant1}4nq^{-\ell\theta_2}\leqslant14nq^{-\theta_2}$. For the term $\ell=0$, Lemmas~\ref{lem:first-generic}, \ref{lem:first-degenerate} and~\ref{lem:universal} give:
\[
\frac{|C_{\fp}(r,0)|}{q^{r(n/2+1)}}\ll\begin{cases}q^{-1}&\text{if $G(\tv)\neq0$ (all $r\geqslant1$)},\\ q^{-1/2}&\text{if $G(\tv)=0$ and $r=1$},\\ 1&\text{if $G(\tv)=0$ and $r\geqslant2$}.\end{cases}
\]
Let $\fs$ denote the product of the primes that divide $\fb_1'$ exactly once. Putting these bounds together, we obtain $f_{\fp}(r)\leqslant C_nq^{-\tau}$ for all $r\geqslant1$ when $G(\tv)\neq0$, with $\tau=1$ if $\psi_1(\tv)=0$ and $\tau=\frac12$ otherwise; and $f_{\fp}(1)\leqslant C_nq^{-1/2}$, $f_{\fp}(r)\leqslant C_n$ for $r\geqslant2$ when $G(\tv)=0$. Hence the total contribution of $\fb_1'$, summed over $\fd_2^1$, is at most
\[
\sum_{\N\fb_1'\leqslant(A_0P)^d}C_n^{\omega(\fb_1')}\times\begin{cases}\N(\rad\fb_1')^{-1}&\text{if }\psi_1(\tv)=0,\\ \N(\fs)^{-1/2}&\text{if }\psi_1(\tv)=\frac12.\end{cases}
\]
In the first case this is $O(P^\varepsilon)$ by Lemma~\ref{lem:smoothsupp}. In the second case write $\fb_1'=\fs\ff$ with $\ff$ square-full and coprime to $\fs$. Since there are $O(y^{1/2})$ square-full ideals of norm at most $y$, and $C_n^{\omega(\fb)}\ll\N\fb^\varepsilon$,
\[
\sum_{\N\fs\ff\leqslant x}C_n^{\omega(\fs\ff)}\N\fs^{-1/2}\ll x^\varepsilon\sum_{\N\ff\leqslant x}\Bigl(\frac x{\N\ff}\Bigr)^{1/2}\ll x^{1/2+2\varepsilon},\qquad x=(A_0P)^d.
\]
Finally, the number of $\fd_1$ is $O(P^\varepsilon)$, as noted above. This proves the lemma.
\end{proof}

\begin{remark}
When $Q_2^\ast(\tv)\neq0$, a prime dividing both moduli must divide $Q_2^\ast(\tv)$, by Lemma~\ref{lem:mixed}; this is the analogue of~\cite[Lemma~9]{Mun15}. When $Q_2^\ast(\tv)=0$, there is no such restriction, and for odd $n$ we lose a factor $P^{d/2}$ in the sum over $\fb_1$. The vectors with $Q_2^\ast(\tv)=0$ are so rare that this loss is harmless; see the table below.
\end{remark}

\subsection{Counting the dual vectors}
\begin{lemma}\label{lem:vcount}
Let $M\geqslant2$. The number of balanced $\tv\in\fo^n\setminus\{\mathbf0\}$ with $\cH(\tv)\leqslant M$ is
\begin{enumerate}[label=\textup{(\roman*)}]
\item $O(M^n)$ in general;
\item $O(M^{n-1/2}\log M)$ if $\tv$ is restricted to be exceptional;
\item $O(M^{n-3/2}\log M)$ if $G(\tv)=0$, or if $Q_2^\ast(\tv)=0$;
\item $O(M^{n-2})$ if $G(\tv)=Q_2^\ast(\tv)=0$.
\end{enumerate}
\end{lemma}

\begin{proof}
A balanced $\tv$ with $\cH(\tv)\leqslant M$ satisfies $\ip\tv\ll M^{1/d}$, and~(i) is Lemma~\ref{lem:count-trivial}(iii). The quadratic form $Q_2^\ast$ has rank $n\geqslant3$, so it is absolutely irreducible and is not a constant multiple of the square of a form. The form $G$ is absolutely irreducible of degree $4(n-2)\geqslant2$ by Lemma~\ref{lem:geometry}\ref{G4}. Part~(iii) therefore follows from the first part of Lemma~\ref{lem:count-BV} with $H\asymp M^{1/d}$. For~(ii), apply the second part of Lemma~\ref{lem:count-BV} to $F=\epsilon_nQ_2^\ast$ and to each of the finitely many classes $a\in\mathcal A_K$. For~(iv), the forms $G$ and $Q_2^\ast$ are distinct and irreducible, since their degrees differ; so they have no common factor, and $\{G=Q_2^\ast=0\}$ has dimension $n-2$. Lemma~\ref{lem:count-trivial}(iii) applies.
\end{proof}

\subsection{Proof of Proposition~\ref{prop:error}}
By Lemma~\ref{lem:tails}, we may restrict~\eqref{eq:Esplit} to $\cH(\tv)\leqslant M_0$. For such $\tv$, Lemma~\ref{lem:factor} gives
\[
\Sigma(\fb_1,\fc,\tv)=\sum_{\fd_2}C^{\cS}_{\fb_1,\fd_2}(\tv)\,\Sigma_{\fd_2},
\]
the sum being over ideals $\fd_2$ supported on $\cS$ with $\N\fd_2\leqslant c_{\ast}P^{d/2}$. Inserting Lemma~\ref{lem:betasum}, we find that $\N\fb_1^{-n}|\Sigma(\fb_1,\fc,\tv)|$ is at most
\[
P^\varepsilon\,P^{\frac{3d}2-\frac{dn}2}\cH(\tv)^{\theta(\tv)-\frac n2}\sum_{\fd_2}\frac{|C^{\cS}_{\fb_1,\fd_2}(\tv)|}{\N(\fb_1\fd_2)^{n/2+1}}\Bigl(\frac{P^{d/2}}{\N\fd_2}\Bigr)^{\theta_2(\tv)},
\]
since $\N\fb_1^{-n}\N\fd_2^{-n}(\N\fb_1\N\fd_2)^{-1+n/2}=\N(\fb_1\fd_2)^{-n/2-1}$. Now interchange the sums over $\fb_1$ and $\tv$ in~\eqref{eq:Esplit}: all terms are non-negative, and we may enlarge the range of $\tv$ to all balanced $\tv\in\fo^n\setminus\{\mathbf0\}$ with $\cH(\tv)\leqslant M_0$. Lemma~\ref{lem:assembly} then gives
\[
E(P)\ll P^{d(\frac n2-\frac32)+\varepsilon}\sum_{\substack{\tv\ \text{balanced}\\ 0<\cH(\tv)\leqslant M_0}}\cH(\tv)^{\theta(\tv)-\frac n2}P^{d\psi(\tv)},\qquad\psi(\tv)=\frac{\theta_2(\tv)}2+\psi_1(\tv).
\]
We split the vectors $\tv$ into the classes listed in the table below. In each class, $\theta$ and $\psi$ are constant, and Lemma~\ref{lem:vcount} bounds the number of $\tv$ with $\cH(\tv)\leqslant M$ by $O(M^{e+\varepsilon})$, with $e$ as in the table. Since $e+\theta-\frac n2>0$ in every class, partial summation gives
\[
\sum_{\cH(\tv)\leqslant M_0}\cH(\tv)^{\theta-\frac n2}\ll M_0^{e+\theta-\frac n2+\varepsilon}=P^{\frac d2(e+\theta-\frac n2)+O(\varepsilon)}.
\]
So the class contributes $O(P^{d\gamma+O(\varepsilon)})$, where
\[
\gamma=\frac n2-\frac32+\psi+\frac12\Bigl(e+\theta-\frac n2\Bigr).
\]

\medskip
\begin{center}
\renewcommand{\arraystretch}{1.3}
\begin{tabular}{l|c|c|c|c|c|c}
class of $\tv$ & $\theta$ & $\theta_2$ & $\psi_1$ & $\psi$ & $e$ & $\gamma$\\\hline
\multicolumn{7}{c}{$n$ odd}\\\hline
$G\neq0$, $Q_2^\ast\neq0$, not exceptional & $\theta_K$ & $0$ & $0$ & $0$ & $n$ & $\frac{3n}4-\frac32+\frac{\theta_K}2$\\
$G\neq0$, exceptional & $0$ & $\frac12$ & $0$ & $\frac14$ & $n-\frac12$ & $\frac{3n}4-\frac32$\\
$G\neq0$, $Q_2^\ast=0$ & $0$ & $\frac12$ & $\frac12$ & $\frac34$ & $n-\frac32$ & $\frac{3n}4-\frac32$\\
$G=0$, $Q_2^\ast\neq0$, not exceptional & $\theta_K$ & $0$ & $\frac12$ & $\frac12$ & $n-\frac32$ & $\frac{3n}4-\frac74+\frac{\theta_K}2$\\
$G=0$, exceptional & $0$ & $\frac12$ & $\frac12$ & $\frac34$ & $n-\frac32$ & $\frac{3n}4-\frac32$\\
$G=0$, $Q_2^\ast=0$ & $0$ & $\frac12$ & $\frac12$ & $\frac34$ & $n-2$ & $\frac{3n}4-\frac74$\\\hline
\multicolumn{7}{c}{$n$ even}\\\hline
$G\neq0$, $Q_2^\ast\neq0$ & $0$ & $0$ & $0$ & $0$ & $n$ & $\frac{3n}4-\frac32$\\
$G\neq0$, $Q_2^\ast=0$ & $0$ & $1$ & $0$ & $\frac12$ & $n-\frac32$ & $\frac{3n}4-\frac74$\\
$G=0$, $Q_2^\ast\neq0$ & $0$ & $0$ & $\frac12$ & $\frac12$ & $n-\frac32$ & $\frac{3n}4-\frac74$\\
$G=0$, $Q_2^\ast=0$ & $0$ & $1$ & $\frac12$ & $1$ & $n-2$ & $\frac{3n}4-\frac32$
\end{tabular}
\end{center}
\medskip

In every row $\gamma\leqslant\Phi(n)$, where we recall that $\theta_K=0$ for even $n$ in~\eqref{eq:theta}. Since $\varepsilon$ is arbitrary, this completes the proof of Proposition~\ref{prop:error}.\qed

\begin{remark}\label{rem:munshi-check}
The first row is the generic term, and it is the only one in which subconvexity enters. For $K=\QQ$ and $n=11$, with the exponent $\theta_\QQ=\frac7{16}$ implicit in~\cite[Lemma~5]{Mun15}, it gives $\frac{33}4-\frac32+\frac7{32}=\frac{223}{32}$, which is the exponent in~\cite[Theorem~1]{Mun15}. With the convexity bound, $\theta_K=\frac12$, the main term $P^{d(n-4)}$ dominates only when $n\geqslant12$.
\end{remark}
\section{The main term}\label{sec:main}

Let $M(P)$ denote the contribution of $\bm=\mathbf0$ in Proposition~\ref{prop:poisson}:
\[
M(P)=\Bigl(\frac{2^{r_2}}{\sqrt{D_K}}\Bigr)^nP^{d(n-3)}\sum_{\fb_1,\fb_2}\frac{C_{\fb_1,\fb_2}(\mathbf0)}{\N(\fb_1\fb_2)^n}\,I_{\fb_1,\fb_2}(\mathbf0).
\]
We define the singular series by
\begin{equation}\label{eq:singularseries}
\fS_K=\sum_{(0)\neq\fb_1,\fb_2\subseteq\fo}\frac{C_{\fb_1,\fb_2}(\mathbf0)}{\N\fb_1^{\,n-1}\N\fb_2^{\,n}}.
\end{equation}

\subsection{The integral at the zero frequency}
\begin{lemma}\label{lem:I0}
\begin{enumerate}[label=\textup{(\roman*)}]
\item For all $\fb_1,\fb_2$ we have $I_{\fb_1,\fb_2}(\mathbf0)\ll(\log P)^{r_1+r_2}$.
\item If $\N\fb_1\leqslant P^{d-\varepsilon}$ and $\N\fb_2\leqslant P^{d/2-\varepsilon}$, then
\[
I_{\fb_1,\fb_2}(\mathbf0)=\Bigl(\frac{\sqrt{D_K}}{2^{r_2}}\Bigr)^2\frac{\N\fb_1}{P^d}\,\sigma_\infty+O_A(P^{-A}).
\]
\end{enumerate}
\end{lemma}

\begin{proof}
Fix $c_1$, and write $A_{c_1}(v)=\bW(c_1/P)-\bW(Pv/c_1)$ and $\bl=\bl(\fb_1,c_1)$. By the definition of $\cW$ as the push-forward of $W(\by)\,d\by$,
\[
I_{c_1,\fb_2}(\mathbf0)=\iint_{V^2}\cW(v_1,v_2)A_{c_1}(v_1)\,h\bigl(\rho_2,\Nm(\bl v_2)\bigr)\,U_V(\bl v_2)\,dv_1\,dv_2.
\]
Substitute $v_2=v_2'/\bl$. Since every $\bl^{(l)}$ is positive, $dv_2=\Nm(\bl)^{-1}dv_2'=\rho_1\,dv_2'$, and this Jacobian does not depend on $c_1$. Thus
\begin{equation}\label{eq:I0a}
I_{c_1,\fb_2}(\mathbf0)=\rho_1\int_VA_{c_1}(v_1)\Psi_{c_1}(v_1)\,dv_1,\qquad\Psi_{c_1}(v_1)=\int_V\cW\Bigl(v_1,\frac{v_2'}{\bl}\Bigr)h(\rho_2,\Nm v_2')U_V(v_2')\,dv_2'.
\end{equation}
Note that $\Psi_{c_1}(v_1)=0$ unless $\ip{v_1}\leqslant A_0/2$.

(i) By Lemma~\ref{lem:hL1} with $R_0=2$, $|\Psi_{c_1}(v_1)|\ll(1+|\log\rho_2|)^{r_1+r_2}\ll(\log P)^{r_1+r_2}$. The proof of Lemma~\ref{lem:Hmain}(ii), with $Q=P$ and $f$ the indicator function of $\{\ip v\leqslant A_0/2\}$, shows that
\[
\sum_{c_1\in\fb_1\setminus\{0\}}\frac1{|\N c_1|}\int_{\ip{v_1}\leqslant A_0/2}|A_{c_1}(v_1)|\,dv_1\ll\frac1{\N\fb_1}.
\]
Hence $\sum_{c_1}|\N c_1|^{-1}|I_{c_1,\fb_2}(\mathbf0)|\ll\rho_1(\log P)^{r_1+r_2}/\N\fb_1=(\log P)^{r_1+r_2}P^{-d}$. Since $P^{2d}/C_P\ll P^d$, (i) follows from~\eqref{eq:Ibig}.

(ii) The function $f(v_2')=\cW(v_1,v_2'/\bl)U_V(v_2')$ is supported in $\{\ip{v_2'}\leqslant2\}$, satisfies $f(0)=\cW(v_1,0)$, and has $\|f\|_{C^k}\ll L^{k/2}$, because $\bl^{(l)}\geqslant A_0^{-1}$ for every $l$. Lemma~\ref{lem:h4} with $x=\rho_2\leqslant P^{-\varepsilon}$ gives
\[
\Psi_{c_1}(v_1)=\frac{\sqrt{D_K}}{2^{r_2}}\cW(v_1,0)+O_N\bigl(P^{-\varepsilon N}L^{k}\bigr),
\]
uniformly in $c_1$ and $v_1$. The cut-off $U_V$ has disappeared from the main term, together with the dependence on $c_1$; this is the analogue of the removal of $U$ in~\cite[(35)]{Mun15}. The error term contributes, after summation over $c_1$ as in~(i), $O(P^d\cdot\rho_1P^{-\varepsilon N}L^k/\N\fb_1)=O(P^{-\varepsilon N}L^k)$ to $I_{\fb_1,\fb_2}(\mathbf0)$, which is negligible. For the main term, observe that for $\ip{v_1}\leqslant A_0/2$ the terms with $\ip{c_1}>A_0P$ vanish in the definition~\eqref{eq:Hb} of $\mathcal H_{\fb_1}(P^2v_1)$ (with $Q=P$). So the sum over $c_1$ in~\eqref{eq:Ibig} produces exactly $\mathcal H_{\fb_1}(P^2v_1)$, and Lemma~\ref{lem:Hmain}(i) with $Q=P$, $f=\cW(\cdot,0)$ and $\rho=\rho_1\leqslant P^{-\varepsilon}$ gives
\begin{align*}
\sum_{c_1}\frac{\rho_1}{|\N c_1|}\int_VA_{c_1}(v_1)\cW(v_1,0)\,dv_1&=\rho_1\int_V\cW(v_1,0)\mathcal H_{\fb_1}(P^2v_1)\,dv_1\\
&=\frac{\rho_1}{P^d}\Bigl(\sigma_\infty\int_V\bW+O_N(P^{-\varepsilon N}L^k)\Bigr),
\end{align*}
since $\cW(0,0)=\sigma_\infty$. Finally, by Theorem~\ref{thm:elementdelta},
\[
\frac{P^{2d}}{C_P}=\frac{\sqrt{D_K}}{2^{r_2}}\,P^d\Bigl(\int_V\bW\Bigr)^{-1}\bigl(1+O_N(P^{-N})\bigr).
\]
Multiplying out gives~(ii).
\end{proof}

\subsection{Evaluation of the main term}
\begin{proposition}\label{prop:main}
Let $n\geqslant9$. Then $\fS_K$ converges absolutely, and
\[
M(P)=\Bigl(\frac{2^{r_2}}{\sqrt{D_K}}\Bigr)^{n-2}\sigma_\infty\fS_KP^{d(n-4)}+O\bigl(P^{d(n-4)-\eta_0}\bigr),\qquad\eta_0=\frac d4(n-8)-\varepsilon.
\]
\end{proposition}

\begin{proof}
By Corollary~\ref{cor:pointwise}, $|C_{\fb_1,\fb_2}(\mathbf0)|\ll\N(\fb_1\fb_2)^{n/2+1+\varepsilon}$. Since the number of ideals of norm at most $x$ is $O(x)$, partial summation gives, for $s>1$ and $Y\geqslant1$,
\begin{equation}\label{eq:tailsum}
\sum_{\N\fb>Y}\N\fb^{-s}\ll_sY^{1-s},\qquad\sum_{\fb}\N\fb^{-s}\ll_s1.
\end{equation}
In particular $\sum|C_{\fb_1,\fb_2}(\mathbf0)|\N\fb_1^{1-n}\N\fb_2^{-n}\ll\sum\N\fb_1^{2-n/2+\varepsilon}\N\fb_2^{1-n/2+\varepsilon}<\infty$, since $n\geqslant7$.

Let $\mathcal R$ be the set of pairs with $\N\fb_1\leqslant P^{d-\varepsilon}$ and $\N\fb_2\leqslant P^{d/2-\varepsilon}$. By Lemma~\ref{lem:I0}(ii) and the identities
\[
\Bigl(\frac{2^{r_2}}{\sqrt{D_K}}\Bigr)^n\Bigl(\frac{\sqrt{D_K}}{2^{r_2}}\Bigr)^2=\Bigl(\frac{2^{r_2}}{\sqrt{D_K}}\Bigr)^{n-2},\qquad P^{d(n-3)}\cdot\frac{\N\fb_1}{P^d}\cdot\frac1{\N\fb_1^{\,n}}=\frac{P^{d(n-4)}}{\N\fb_1^{\,n-1}},
\]
the pairs in $\mathcal R$ contribute $(2^{r_2}/\sqrt{D_K})^{n-2}\sigma_\infty P^{d(n-4)}$ times the part of~\eqref{eq:singularseries} over $\mathcal R$, up to $O(P^{-A})$. By~\eqref{eq:tailsum}, the part of~\eqref{eq:singularseries} outside $\mathcal R$ is
\[
\ll P^{-(d-\varepsilon)(n/2-3)+\varepsilon}+P^{-(d/2-\varepsilon)(n/2-2)+\varepsilon},
\]
and since $\sigma_\infty\ll1$ this is acceptable. The pairs outside $\mathcal R$ contribute to $M(P)$, by Lemma~\ref{lem:I0}(i) and~\eqref{eq:tailsum},
\begin{multline*}
\ll P^{d(n-3)+\varepsilon}\Bigl(\sum_{\N\fb_1>P^{d-\varepsilon}}\N\fb_1^{1-\frac n2+\varepsilon}\sum_{\fb_2}\N\fb_2^{1-\frac n2+\varepsilon}+\sum_{\fb_1}\N\fb_1^{1-\frac n2+\varepsilon}\sum_{\N\fb_2>P^{d/2-\varepsilon}}\N\fb_2^{1-\frac n2+\varepsilon}\Bigr)\\
\ll P^{d(\frac n2-1)+2\varepsilon'}+P^{d(\frac{3n}4-2)+2\varepsilon'},
\end{multline*}
where $\varepsilon'$ depends on $\varepsilon$ and tends to $0$ with it. Compared with $P^{d(n-4)}$, these terms save $P^{d(n/2-3)}$ and $P^{d(n-8)/4}$ respectively, up to $P^{O(\varepsilon)}$. For $n\geqslant9$ the second saving is the smaller one.
\end{proof}

\subsection{The singular series}
For a prime ideal $\fp$, with $q=\N\fp$, put
\[
\sigma_{\fp}=\sum_{r,\ell\geqslant0}\frac{C_{\fp}(r,\ell;\mathbf0)}{q^{r(n-1)+\ell n}}.
\]
For integers $k\geqslant j\geqslant0$ put
\[
N_k(j)=\#\{\bb\bmod\fp^k:\ Q_2(\bb)\equiv0\bmod\fp^k,\ Q_1(\bb)\equiv0\bmod\fp^j\},\quad\nu_k(j)=q^{\,j-k(n-1)}N_k(j).
\]

\begin{lemma}\label{lem:euler}
Let $n\geqslant7$.
\begin{enumerate}[label=\textup{(\roman*)}]
\item The series defining $\sigma_{\fp}$ converges absolutely, and $\fS_K=\prod_{\fp}\sigma_{\fp}$, the product converging absolutely.
\item If $\fp$ is good, then $\sigma_{\fp}=1+O(q^{-2})$, with an implied constant depending only on $n$.
\item For each $j\geqslant0$ the limit $\nu(j)=\lim_{k\to\infty}\nu_k(j)$ exists, and $\sigma_{\fp}=\lim_{j\to\infty}\nu(j)$.
\item If $X(K_{\fp})\neq\emptyset$, then $\sigma_{\fp}>0$.
\end{enumerate}
In particular, if $X$ has a point over every non-archimedean completion of $K$, then $\fS_K>0$.
\end{lemma}

\begin{proof}
(i) This follows from Corollary~\ref{cor:pointwise}, the absolute convergence of~\eqref{eq:singularseries}, and Lemma~\ref{lem:mult} with $\bm=\mathbf0$.

(ii) Let $\fp$ be good. The normalised frequency is $\mathbf0$, and $G(\mathbf0)=Q_2^\ast(\mathbf0)=0$. By Lemma~\ref{lem:first-degenerate}, the term $(r,\ell)=(1,0)$ is at most $3nq^{(n+1)/2-(n-1)}\leqslant3nq^{-2}$. The term $(0,1)$ is $\mathcal Q_{\fp}(1;\mathbf0)q^{-n}$, which vanishes for odd $n$ by Lemma~\ref{lem:gauss}(i), and is at most $q^{n/2+1-n}\leqslant q^{-3}$ for even $n\geqslant8$ by Lemma~\ref{lem:gauss}(ii). For $r+\ell\geqslant2$, Lemma~\ref{lem:universal} bounds the term by $4nq^{r(2-n/2)+\ell(1-n/2)}\leqslant4nq^{-\frac32r-\frac52\ell}$. Summing over $r+\ell\geqslant2$ gives $O(q^{-3})$.

(iii) Let $r\geqslant0$ and $L\geqslant0$, and put $k=r+L$. In Lemma~\ref{lem:dictionary} with $\bm=\mathbf0$, the summand depends only on $\bb$ modulo $\fp^{r+\ell}$ (Lemma~\ref{lem:periodic}), so we may let $\bb$ run modulo $\fp^k$ for every $\ell\leqslant L$, at the cost of a factor $q^{-n(L-\ell)}$. The sum over $a_1$ is the Ramanujan sum $c_{\fp^r}(Q_1(\bb))=\sumstar_{a_1\bmod\fp^r}\ep r(a_1Q_1(\bb))$. Writing $Q_2(\bb)=\varpi^ry$, the sum over $a_2$ is $\sumstar_{a_2\bmod\fp^\ell}\ep\ell(a_2y)$. Every additive character modulo $\fp^L$ is primitive modulo exactly one $\fp^\ell$ with $\ell\leqslant L$, so
\[
\sum_{\ell=0}^L\ \sumstar_{a_2\bmod\fp^\ell}\ep\ell(a_2y)=\sum_{a\bmod\fp^L}\ep L(ay)=q^L\,\mathbf1[\fp^L\mid y].
\]
Consequently
\[
\sum_{\ell=0}^L\frac{C_{\fp}(r,\ell;\mathbf0)}{q^{r(n-1)+\ell n}}=q^{-r(n-1)-L(n-1)}\sum_{\substack{\bb\bmod\fp^k\\ Q_2(\bb)\equiv0\ (\fp^k)}}c_{\fp^r}(Q_1(\bb)).
\]
For $r=0$ the Ramanujan sum is $1$ and the right-hand side is $\nu_L(0)$. For $r\geqslant1$ we have $c_{\fp^r}(x)=q^r\mathbf1[\fp^r\mid x]-q^{r-1}\mathbf1[\fp^{r-1}\mid x]$, and the right-hand side is $\nu_k(r)-\nu_k(r-1)$. Letting $L\to\infty$ and using (i), we see first that $\nu_k(0)$ converges as $k\to\infty$. Then, by induction on $j$, $\nu_k(j)=\nu_k(j-1)+(\nu_k(j)-\nu_k(j-1))$ converges as $k\to\infty$ for every $j$, and
\[
\sum_{\ell\geqslant0}\frac{C_{\fp}(j,\ell;\mathbf0)}{q^{j(n-1)+\ell n}}=\nu(j)-\nu(j-1)\qquad(j\geqslant1).
\]
Summing over $0\leqslant j\leqslant R$ gives $\sum_{r\leqslant R}\sum_\ell C_{\fp}(r,\ell;\mathbf0)q^{-r(n-1)-\ell n}=\nu(R)$, and letting $R\to\infty$ proves (iii).

(iv) Let $\mu$ be the Haar measure on $\cO_{\fp}^n$ with $\mu(\cO_{\fp}^n)=1$. For $k\geqslant j$, the set $\{\bx\in\cO_{\fp}^n:Q_2(\bx)\in\fp^k,\ Q_1(\bx)\in\fp^j\}$ is a union of cosets of $(\fp^k)^n$, so its measure is $q^{-kn}N_k(j)$, and
\[
\nu_k(j)=q^{\,j+k}\mu\bigl(\{\bx\in\cO_{\fp}^n:Q_2(\bx)\in\fp^k,\ Q_1(\bx)\in\fp^j\}\bigr).
\]
Suppose that $X(K_{\fp})\neq\emptyset$, and let $\bx_0\in\cO_{\fp}^n\setminus\{\mathbf0\}$ with $Q_1(\bx_0)=Q_2(\bx_0)=0$. Since $X$ is non-singular, $\nabla Q_1(\bx_0)$ and $\nabla Q_2(\bx_0)$ are linearly independent, by Lemma~\ref{lem:geometry}\ref{G3}. Choose two coordinates, say $x_{n-1},x_n$, for which the Jacobian $J(\bx)=\det\partial(Q_1,Q_2)/\partial(x_{n-1},x_n)$ satisfies $J(\bx_0)\neq0$. By the $\fp$-adic implicit function theorem (see~\cite[Chap.~2]{Igu00}), there are $e\geqslant1$ and an analytic bijection $\Phi$ from $B'\times(\fp^e)^2$, where $B'=\bx_0'+(\fp^e)^{n-2}$ is a ball in the first $n-2$ coordinates, onto an open subset of $\cO_{\fp}^n$, such that $Q_i(\Phi(\bz',y_1,y_2))=y_i$ for $i=1,2$ and $|J(\bx)|_{\fp}=|J(\bx_0)|_{\fp}$ on the image. By the change of variables formula for $\fp$-adic analytic maps~\cite[Chap.~7]{Igu00}, for $k\geqslant j\geqslant e$,
\[
\mu\bigl(\{\bx:Q_2(\bx)\in\fp^k,\ Q_1(\bx)\in\fp^j\}\bigr)\geqslant\mu\bigl(\Phi(B'\times\fp^j\times\fp^k)\bigr)=|J(\bx_0)|_{\fp}^{-1}\,q^{-e(n-2)}q^{-j}q^{-k}.
\]
Hence $\nu_k(j)\geqslant|J(\bx_0)|_{\fp}^{-1}q^{-e(n-2)}>0$ for all $k\geqslant j\geqslant e$, and (iii) shows that $\sigma_{\fp}>0$.

Finally, if every $\sigma_{\fp}$ is positive then so is the absolutely convergent product $\prod_{\fp}\sigma_{\fp}$.
\end{proof}

\begin{remark}
One can show, using~(iii), that $\sigma_{\fp}$ coincides with the classical local density $\lim_{k\to\infty}q^{-k(n-2)}\#\{\bb\bmod\fp^k:Q_1(\bb)\equiv Q_2(\bb)\equiv0\bmod\fp^k\}$. We shall not need this.
\end{remark}

For $n\geqslant9$, the variety $X$ has a point over every non-archimedean completion of $K$, by the results of Demyanov~\cite{Dem56} and Birch, Lewis and Murphy~\cite{BLM62} recalled in the introduction. Hence $\fS_K>0$ under the hypotheses of Theorem~\ref{thm:main}.

\section{Proof of the main results}\label{sec:proofs}

\begin{proof}[Proof of Theorem~\ref{thm:main}]
By Proposition~\ref{prop:poisson}, $N(P)=M(P)+E(P)+O(P^{-A})$. Proposition~\ref{prop:main} evaluates $M(P)$, Lemma~\ref{lem:euler} shows that $\fS_K$ is a convergent product of local densities and is positive, and Proposition~\ref{prop:error} gives $E(P)\ll P^{d\Phi(n)+\varepsilon}$. Hence
\[
N(P)=\Bigl(\frac{2^{r_2}}{\sqrt{D_K}}\Bigr)^{n-2}\fS_K\,\sigma_\infty\,P^{d(n-4)}+O\Bigl(P^{d\Phi(n)+\varepsilon}+P^{d(n-4)-\frac d4(n-8)+\varepsilon}\Bigr).
\]
Since $n-4-\Phi(n)=\frac14(n-10-2\theta_K)$, the error term is $O(P^{d(n-4)-\eta})$ with $\eta=\min\{\frac d4(n-10-2\theta_K),\frac d4(n-8)\}-\varepsilon$. For $n\geqslant11$ this is positive, since $\theta_K<\frac12$ by Theorem~\ref{thm:MV}. All the estimates are uniform for $\bxi$ in a fixed compact subset of $V^n\setminus X^\ast(V)$: the only dependence on $\bxi$ is through $\delta_0$, the constant $c_3$ of Lemma~\ref{lem:gaussian}, and the norm of the inverse of the map $T$ in the proof of Proposition~\ref{prop:statphase}, all of which are bounded on such a set. The statement about $\sigma_\infty$ was proved in~\S\ref{sec:counting}.
\end{proof}

\begin{proof}[Proof of Corollary~\ref{cor:boxes}]
For $\bxi\in V^n$, let $W_{\bxi}$ denote the weight~\eqref{eq:weight}, and $N_{\bxi}(P)$ the corresponding counting function. Let $c_L=\bigl(\int_{V^n}e^{-L\|\bz\|^2}d\bz\bigr)^{-1}\asymp L^{dn/2}$. Let $\mathcal B^+$ and $\mathcal B^-$ be the sets of points at distance at most $(\log P)^{-1}$ from $\mathcal B$, respectively at distance at least $(\log P)^{-1}$ from the complement of $\mathcal B$, and put
\[
g^\pm(\by)=c_L\int_{\mathcal B^\pm}W_{\bxi}(\by)\,d\bxi.
\]
If $\by\notin\mathcal B$, then every $\bxi\in\mathcal B^-$ has $\|\by-\bxi\|\geqslant(\log P)^{-1}$, so $g^-(\by)\ll c_Le^{-L(\log P)^{-2}}=c_Le^{-(\log P)^2}$; and $g^-\leqslant1$ everywhere. If $\by\in\mathcal B$, then $w_0=1$ for all $\bxi$ with $\ip{\by-\bxi}\leqslant\delta_0/2$, and the set of $\bxi$ with $\|\by-\bxi\|\leqslant(\log P)^{-1}$ lies in $\mathcal B^+$; so $g^+(\by)\geqslant1-O(e^{-c(\log P)^2})$. Since the number of $\bx\in\fo^n$ with $\ip\bx\ll P$ is $O(P^{dn})$, we obtain
\[
\begin{split}
\sum_{\bx}g^-(\bx/P)-O(P^{-A})&\leqslant\#\{\bx\in\fo^n\cap P\mathcal B:Q_1(\bx)=Q_2(\bx)=0\}\\
&\leqslant\bigl(1+O(P^{-A})\bigr)\sum_{\bx}g^+(\bx/P),
\end{split}
\]
the sums being over the zeros $\bx\in\fo^n$ of $Q_1=Q_2=0$. Now $\sum_{\bx} g^\pm(\bx/P)=c_L\int_{\mathcal B^\pm}N_{\bxi}(P)\,d\bxi$. Since the closure of $\mathcal B^+$ is a compact subset of $V^n\setminus X^\ast(V)$ for $P$ large, Theorem~\ref{thm:main} applies uniformly for $\bxi\in\mathcal B^\pm$. As $c_L\ll(\log P)^{2dn}$, the error terms contribute $O(P^{d(n-4)-\eta/2})$. In the main term,
\[
c_L\int_{\mathcal B^\pm}\sigma_\infty(W_{\bxi})\,d\bxi=\sigma_\infty(g^\pm),
\]
where $\sigma_\infty(g)$ denotes the singular integral of a weight $g$, that is, the density at $(0,0)$ of the push-forward of $g(\by)\,d\by$ under $(Q_1,Q_2)$. The functions $g^\pm$ differ from the indicator function of $\mathcal B$ by $O(e^{-c(\log P)^2})$, except on the set of points at distance at most $2(\log P)^{-1}$ from the boundary of $\mathcal B$. Near $\mathcal B$, the fibre $\{Q_1=Q_2=0\}$ is smooth, and this set has measure $O((\log P)^{-1})$ for the induced measure on it. Hence $\sigma_\infty(g^\pm)=\sigma_\infty(\mathcal B)+O((\log P)^{-1})$, and the corollary follows.
\end{proof}

\section{Why weights that depend only on norms fail}\label{sec:higher}

In this final section we explain why the nested argument needs the element $\delta$-symbol of Theorem~\ref{thm:elementdelta} when $K$ has more than one archimedean place, and why the ideal $\delta$-symbol of Browning and Vishe suffices for the second equation. The relevant invariant is the number $r_1+r_2$ of archimedean places, not the degree. For instance, an imaginary quadratic field has $d=2$ but only one archimedean place, and for such fields the difficulty described below does not arise.

\subsection{Sizes over $\QQ$ and over $K$}
Over $\QQ$, in Munshi's argument the dual variables satisfy $|t_1|\lesssim(q_1/B)^{-1}$ and $|t_2|\lesssim(q_2/B^{1/2})^{-1}$, and after the rescaling by $\lambda=B/q_1$ that comes from the nesting, the stationary point of the oscillatory integral lies at a frequency of size at most $\max(|t_1|,\lambda|t_2|)\lesssim B^{3/2}/(q_1q_2)$. This is exactly the truncation $|\bm|\lesssim B^{1/2}$ after normalisation.

Over $K$ all sizes are vectors. The weight $h(\rho,|\Nm v|)$ of Theorem~\ref{thm:BV} depends on $v$ only through its norm. By Lemma~\ref{lem:p2}, its Fourier transform is controlled only through $\cH(t)$, and in fact it is genuinely large on the whole hyperbolic region $\cH(t)\lesssim\rho^{-1}$. This region contains points at which all of the size of $t$ is concentrated at a single place. Heuristically, this happens because $\Nm v$ vanishes on every coordinate hyperplane $\{v^{(l)}=0\}$, so that $h(\rho,|\Nm v|)$ has singular behaviour along each of them.

\subsection{The failure of a scalar nesting}
Suppose that both equations were detected by the ideal $\delta$-symbol, with the nesting carried out with a scalar cut-off $U_V(\lambda Q_2(\by))$, $\lambda=\rho_1^{-1/d}$. The proof of Theorem~\ref{thm:truncation} would then only give $|P\bm^{(l)}|\lesssim\max(1,|t_1^{(l)}|,\lambda|t_2^{(l)}|)$, with $\cH(t_1)\lesssim\rho_1^{-1}$ and $\cH(t_2)\lesssim\rho_2^{-1}$. The following lemma shows that this loses a power of $\lambda$ when $r_1+r_2\geqslant2$.

\begin{lemma}\label{lem:profile}
Let $0<\rho_1,\rho_2\leqslant1$ and $\lambda=\rho_1^{-1/d}$, and let $t_1,t_2\in V$ with $\cH(t_1)\leqslant\rho_1^{-1}$ and $\cH(t_2)\leqslant\rho_2^{-1}$. Let $c_{\min}=\min_lc_l$. Then
\[
\prod_l\max\bigl(1,|t_1^{(l)}|,\lambda|t_2^{(l)}|\bigr)^{c_l}\leqslant(\rho_1\rho_2)^{-1}\lambda^{d-c_{\min}}.
\]
If $r_1+r_2\geqslant2$, equality holds up to a constant for suitable $t_1,t_2$. If $r_1+r_2=1$, the left-hand side is at most $(\rho_1\rho_2)^{-1}$.
\end{lemma}

\begin{proof}
Put $A_l=\max(1,|t_1^{(l)}|)^{c_l}$ and $B_l=\lambda^{c_l}\max(1,|t_2^{(l)}|)^{c_l}$, so that $\prod_lA_l\leqslant\lambda^d$ and $\prod_lB_l\leqslant\lambda^d\rho_2^{-1}$. Then
\[
\prod_l\max(A_l,B_l)=\prod_lB_l\prod_l\max\bigl(1,A_lB_l^{-1}\bigr)\leqslant\prod_lB_l\prod_l\max\bigl(1,A_l\lambda^{-c_l}\bigr).
\]
Let $\mathcal L$ be the set of places $l$ with $A_l>\lambda^{c_l}$. Since every $A_l\geqslant1$, we have $\prod_{l\in\mathcal L}A_l\leqslant\prod_lA_l\leqslant\lambda^d$, and therefore
\[
\prod_l\max\bigl(1,A_l\lambda^{-c_l}\bigr)=\prod_{l\in\mathcal L}A_l\lambda^{-c_l}\leqslant\lambda^{d-\sum_{l\in\mathcal L}c_l},
\]
which is $1$ if $\mathcal L=\emptyset$ and at most $\lambda^{d-c_{\min}}$ otherwise. Since $\prod_lB_l\leqslant\lambda^d\rho_2^{-1}=(\rho_1\rho_2)^{-1}$, the inequality follows. For the equality case, take $t_1$ concentrated at a place $l_1$ with $c_{l_1}=c_{\min}$, with $|t_1^{(l_1)}|^{c_{\min}}=\rho_1^{-1}$, and $t_2$ concentrated at another place, with $\cH(t_2)=\rho_2^{-1}$. If there is only one place, then $\max(A,B)\leqslant\lambda^d\rho_2^{-1}$ directly.
\end{proof}

For a real quadratic field, $t_1=(\rho_1^{-1},1)$ and $t_2=(1,\rho_2^{-1})$ give the value $\rho_1^{-3/2}\rho_2^{-1}$, an extra factor $\lambda=\rho_1^{-1/2}$ compared with Theorem~\ref{thm:truncation}. Heuristically, these frequencies are really present. The Fourier transform of $h(\rho_1,|\Nm v|)$ has size about $1$ at such anisotropic points; the corresponding stationary points of the Gaussian integral are non-degenerate; and the anisotropic part of the region carries a volume of $(t_1,t_2)$ comparable to the isotropic part. If this is so, then taking absolute values would lose a factor up to $P^{(d-c_{\min})(n/2+\theta_K)}$ in the sum over $\bm$ for small $\N\fb_1$. For a real quadratic field this exceeds the main term for every $n$. In this sense, a nested $\delta$-method built only on weights that depend on norms fails when $r_1+r_2\geqslant2$.

\subsection{The remedy}
The element $\delta$-symbol of Theorem~\ref{thm:elementdelta} replaces the norm $\Nm a_1$ by the vector $a_1/(c_1P)$. By Lemma~\ref{lem:p1}, its frequency is confined, at each place separately, to $|t_1^{(l)}|\lesssim P/|c_1^{(l)}|$. The nesting vector $\bl(\fb_1,c_1)$ of~\eqref{eq:lambdavec} has the same archimedean sizes as $c_1^{-1}$, up to a scalar, while keeping $\Nm\bl=P^d/\N\fb_1$. Consequently $|t_1^{(l)}|\lesssim\bl^{(l)}$ at every place; in the notation of the proof of Lemma~\ref{lem:profile}, $A_l\lesssim\lambda_l^{c_l}$ for every $l$, so $\mathcal L=\emptyset$ and the extra factor disappears. This is Theorem~\ref{thm:truncation}.

Two features of this construction deserve comment. First, the archimedean information is kept place by place while the arithmetic information (the moduli $\fb_1,\fb_2$, the exponential sums $C_{\fb_1,\fb_2}(\bm)$, and the average over $\fb_2$) is unchanged. This is analogous to the passage from functions of ideals to functions on $K_\infty^\times$ in the work of Blomer and Brumley~\cite{BB11}. The infinite unit group causes no difficulty because $\bW$ has compact support in $V^\times$; the price is the extra summation over $c_1\in\fb_1$, which costs only $(\log P)^{r_1+r_2}/\N\fb_1$ by Lemma~\ref{lem:latticesum}. Second, the choice of $\bl$ adapted to the archimedean sizes of $c_1$ is reminiscent of the systematic choice of generators of principal ideals with controlled archimedean size in the work of Qi~\cite[(2.1)]{Qi24}.

Finally, the second equation may still be detected with the ideal $\delta$-symbol, because no domination at individual places is needed for $t_2$: in the key inequality of Theorem~\ref{thm:truncation}, the second frequency enters only through $\prod_l\max(1,|t_2^{(l)}|)^{c_l}=\cH(t_2)$.

\subsection{Related approaches}
Qi's Dirichlet approximation theorem over number fields~\cite[Lemma~A.2]{Qi19} shows that every $\theta\in V$ admits approximations $\alpha/\beta$, with $\alpha,\beta\in\fo$ and $\beta\neq0$, such that $|\beta|_v\ll Q$ and $|\beta\theta-\alpha|_v\ll Q^{-1}$ at \emph{every} archimedean place $v$, the implied constants depending only on $K$. This isotropic covering is the analogue, in the language of the circle method, of the place-by-place control used here. It would be interesting to compare our construction with an adelic $\delta$-symbol in the spirit of Getz~\cite{Get18}, and with the two-dimensional $\delta$-symbol of Li, Rydin Myerson and Vishe~\cite{LRV24}, whose extension to number fields might lower the range to $n\geqslant10$.

\subsection*{Acknowledgements}
The author is grateful to Professor Ritabrata Munshi for suggesting the problem and for helpful discussions. The author also gratefully acknowledges the Indian Statistical Institute, Kolkata, for its excellent research environment and support through the Research Associate Fellowship.


\end{document}